\documentclass[graybox]{svmult}

\usepackage{mathptmx}       % selects Times Roman as basic font
\usepackage{helvet}         % selects Helvetica as sans-serif font
\usepackage{courier}        % selects Courier as typewriter font
\usepackage{type1cm}        % activate if the above 3 fonts are
\usepackage{makeidx}         % allows index generation
\usepackage{graphicx}        % standard LaTeX graphics tool
\usepackage{multicol}        % used for the two-column index
\usepackage[bottom]{footmisc}% places footnotes at page bottom

\usepackage{amsmath}
\usepackage{amssymb}
\usepackage{amscd}
\usepackage{indentfirst}
\usepackage{hyperref}
\usepackage{ulem}

\usepackage{mathtools}
\mathtoolsset{showonlyrefs=true}

\usepackage{nicefrac}
\usepackage{algorithm,algorithmic}
\usepackage{makecell}
\usepackage{xcolor}

\usepackage[flushleft]{threeparttable} % http://ctan.org/pkg/threeparttable
\usepackage{caption}
\usepackage{multirow}

\newcommand{\Exp}{\mathbb{E}}
\newcommand{\la}{\langle}
\newcommand{\ra}{\rangle}

\newcommand{\lp}{\left(}
\newcommand{\rp}{\right)}
\newcommand{\lf}{\left\{}
\newcommand{\rf}{\right\}}
\newcommand{\ls}{\left[}
\newcommand{\rs}{\right]}

\def\eps{{\varepsilon}}

\def\1{{\bf 1}}

\newtheorem{Def}{Definition}[section]

\newtheorem{assumption}{Assumption}[section]

\newcommand{\EE}{\mathbb{E}}
\newcommand{\R}{{\mathbb R}}

\def\tf{\tilde{f}}
\def\tg{\tilde{g}}

\DeclareMathOperator*{\argmin}{argmin}

\def\eduard#1{{\color{black} #1}}%blue
 
\def\alexander#1{{\color{black} #1}}%orange
\def\pd#1{{\color{black}#1}} %Pavel's corrections magenta 
\def\dk#1{{\color{black} #1}}%green
\def\kesha#1{{\color{black} #1}}%brown
\makeindex             % used for the subject index
\begin{document}

% \addtocontents{toc}{\protect\setcounter{tocdepth}{-10}}

\title*{Recent Theoretical Advances in Non-Convex Optimization}
\titlerunning{Recent theoretical advances in non-convex optimization}
\author{Marina Danilova$^{1,2}$, Pavel Dvurechensky$^{3,4}$, Alexander Gasnikov$^{2,4,5}$,\\ Eduard Gorbunov$^{2,4}$, Sergey Guminov$^{4}$, Dmitry Kamzolov$^{2,6}$,\\ Innokentiy Shibaev$^{2,4}$ }
% \authorrunning{M.~Danilova, P.~Dvurechensky, A.~Gasnikov, E.~Gorbunov, S.~Guminov, D.~Kamzolov, I.~Shibaev}
\authorrunning{Danilova, Dvurechensky, Gasnikov, Gorbunov, Guminov, Kamzolov, Shibaev}
\institute{
$^1$Institute of Control Sciences RAS, Moscow, Russia\\
$^2$Moscow Institute of Physics and Technology, Moscow, Russia\\
$^3$Weierstrass Institute for Applied Analysis and Stochastics, Berlin, Germany\\
$^4$HSE University, Moscow, Russia\\
$^5$Institute for Information Transmission Problems RAS, Moscow, Russia\\
$^6$ Mohamed bin Zayed University of Artificial Intelligence, Abu Dhabi, United Arab Emirates
}

\maketitle

\abstract{Motivated by recent increased interest in optimization algorithms for non-convex optimization in application to training deep neural networks and other optimization problems in data analysis, we give an overview of recent theoretical results on global performance guarantees of optimization algorithms for non-convex optimization. We start with classical arguments showing that general non-convex problems could not be solved efficiently in a reasonable time. Then we give a list of problems \eduard{that} can be solved efficiently to find the global minimizer by exploiting the structure of the problem as much as it is possible. Another way to deal with non-convexity is to relax the goal from finding the global minimum to finding a stationary point or a local minimum. For this setting, we first present known results for the convergence rates of deterministic first-order methods, which are then followed by a general theoretical analysis of optimal stochastic and randomized gradient schemes, and an overview of the stochastic first-order methods. After that, we discuss quite general classes of non-convex problems, such as minimization of $\alpha$-weakly-quasi-convex functions and functions that satisfy Polyak--Łojasiewicz condition, which still allow obtaining theoretical convergence guarantees of first-order methods. Then we consider higher-order and zeroth-order/derivative-free methods and their convergence rates for non-convex optimization problems.}

%%%%%%%%%%%%%%%%%%%%%%%%%%%%%%%%%%%%%%%%%%%%%%%%%%%%%%%%%%%%%%%%%%%%%%%%%%%%%%%%%%%%%%%%%%%
%%%%%%%%%%%%%%%%%%%%%%%%%%%%%%%%%%%%%%%%%%%%%%%%%%%%%%%%%%%%%%%%%%%%%%%%%%%%%%%%%%%%%%%%%%%

\clearpage
\begingroup
\let\cleardoublepage\clearpage
\setcounter{minitocdepth}{2}
\dominitoc
\endgroup

\section{Introduction}

In this \eduard{survey,} we consider non-convex optimization problems in different settings, including stochastic optimization. We are mainly motivated by an increased interest in such problems in connection to applications in machine learning and data analysis, and our main focus is on the methods which possess theoretical guarantees for their global convergence rate or complexity. As we explain first by providing classical examples \cite{murty1987some,nesterov2018lectures}, there is no hope to have any theoretical guarantees for finding a global minimizer in a general non-convex optimization problem in \eduard{a} reasonable time. Despite the \eduard{quite good} practical performance of classical \eduard{general-purpose} methods such as L-BFGS \cite{nocedal2006numerical,floudas2008encyclopedia}\eduard{, and} \eduard{proven} local superlinear convergence\eduard{, their} global complexity is not well understood.

In the last 20 \eduard{years,} theoretical analysis of the global convergence rate or global complexity guarantees has become de facto  a standard in the area of numerical optimization. Since \eduard{the} convexity of the problem allows for such an analysis, many global complexity and convergence results have been obtained in convex optimization \cite{ben-tal2001lectures,bubeck2015convex,nesterov2018lectures,lan2020first,dvurechensky2020advances,dvurechensky2021first-order}. Recent advances in machine learning, which were made possible by the application of neural networks, had lead to the optimization community changing focus to non-convex optimization and, especially to stochastic non-convex optimization. In this \eduard{non-exhaustive survey,} we \eduard{attempt} to highlight existing results on global performance guarantees of large-scale non-convex optimization methods. \eduard{The large} dimension of the decision variable in such problems motivates the use of first-order methods, which possess a cheap iteration. Moreover, the large amount of data motivates to use randomized methods such as stochastic gradient descent, which does not require to look through the whole dataset to make one step of the optimization procedure, thus making the iteration even cheaper. 

Since, in general, non-convex optimization \eduard{problems} cannot be made efficiently \eduard{solved}, we consider several ways to relax this challenging goal. The first relaxation consists \eduard{of} finding problems with hidden convexity or in a convex reformulation of the problem. This requires \eduard{exploitation} of the problem structure as much as it is possible, which limits the generality of the approach, yet leading to a possibility to find a global solution. Another way is to change the goal from finding the global solution to finding a stationary point or a local extremum. In this \eduard{case,} it is possible to obtain polynomial dependence of the complexity of first-order methods on the dimension of the problem and desired accuracy. We consider this approach in the setting of deterministic and stochastic optimization. The third way is to define a class of non-convex problems, which is\eduard{, on} the one \eduard{hand,} quite general, and on the other hand, allows to obtain \eduard{a global} performance guarantees of an algorithm. We consider a class of problems with objective satisfying Polyak--Łojasiewicz condition, which leads to global linear convergence rate, and the class of problems with $\alpha$-weakly-quasi-convex objective, which leads to global sublinear convergence rate. In the above two \eduard{approaches,} we first focus on first-order methods. Then, motivated by several settings in machine learning such as reinforcement learning, black-box adversarial attacks on neural networks, as well as simulation optimization, in which the gradient of the objective is not available, we consider zeroth-order or derivative-free methods and their convergence rates for non-convex optimization problems. By no means we claim that our \eduard{survey} contains all the important results in this area since the literature is huge and we could miss some recent results. We would like to list here some other books \alexander{\cite{polyak1987introduction,conn2009introduction,lan2020first,GasnikovBook}} and \eduard{survey}s  \cite{jain2017non-convex,curtis2017optimization,wright2018optimization,chen2018harnessing,chi2018nonconvex,sun2019optimization,zhang2020from} related to our paper\footnote{See also this webpage with the list of references being updated \url{https://sunju.org/research/nonconvex/}.}.

\section{Preliminaries}

The main challenges in non-convex optimization are caused either by non-convexity of the feasible set or by non-convexity of the objective function. The first case is tightly connected with discrete \eduard{optimization} when the decision variable can take only a discrete set of values. In the second case, yet the variable can take a continuum number of values, the non-convexity of the problem does not allow to hope for finding a global solution in a reasonable amount of time. We start with two particular examples \eduard{that} illustrate the intractability of non-convex optimization in general. This intractability motivates different kinds of relaxations, such as changing the goal to the one consisting \eduard{of} finding an approximate stationary point instead of a global minimum, or introducing additional assumptions on the problem, or heavily using the structure of the problem, which lead to provable convergence to the global minimizer. Next, we present general non-convex optimization problems and some ways to classify them.

\subsection{Global Optimization is NP-hard}
Following \cite{murty1987some}, we consider an example which illustrates that the problem of finding the exact global solution of a non-convex \eduard{problem is} NP-hard. To that end, we consider the minimization problem
\[
\mathop {\min }\limits_{x\in 
{\mathbb{R}}^n} \left\{f\left( x \right):=\sum\limits_{i=1}^n {x_i^4 } -\frac{1}{n}\left( 
{\sum\limits_{i=1}^n {x_i^2 } } \right)^2+\left( {\sum\limits_{i=1}^n {a_i 
x_i } } \right)^4+\left( {1-x_1 } \right)^4 \right\},
\]
where $x_i$ is the $i$-th component of the vector $x$.
Let $A=I-\frac{1}{n}\mathbf{1}\mathbf{1}^{\top}$, where $I$ is the identity matrix of size $n$ and $\mathbf{1}$ is a vector of $n$ ones, and let $[x]^2$ denote a vector with components $[x]^2_i=x_i^2$. In this notation, the objective takes the form
\[
f\left( x \right)=\langle A[x]^2,[x]^2\rangle+\left( {\sum\limits_{i=1}^n {a_i 
x_i } } \right)^4+\left( {1-x_1 } \right)^4.
\]
Since $A$ is a positive semidefinite matrix, $f(x)\geqslant 0$. One may also note that 0 is an eigenvalue of $A$ with multiplicity $1$ and that $\mathbf{1}$ is the corresponding eigenvector. With this in mind, it is not difficult to see that $f(x)=0$ if and only if $x$ satisfies
\[a_1 +\sum\limits_{i=2}^n {a_i x_i } =0,
\quad
x_i =\pm 1,\ i=2,\ldots, n.
\]
The problem of checking whether this equation has a solution is a form of the subset sum problem, which is known to be NP-complete. Since this problem has a solution if and only if the global minimum in the original optimization problem is exactly zero, this implies that the problem of finding even the value of a global minimum for a non-convex objective is NP-hard.

\subsection{Lower Complexity Bound for Global Optimization}
\label{S:lower bound}
Following \cite{nesterov2018lectures}, we now derive a lower bound for the complexity of finding an approximate global minimum of a possibly non-convex objective. Consider the problem 
\[
\mathop{\min }\limits_{x\in \left[ {0,1} \right]^n} f\left( x \right) , 
\] 
where $f$ is possibly non-convex and Lipschitz-continuous function, i.e.,  for some $M>0$ and for all 
% $\forall$
$x,y\in \left[ {0,1} \right]^n$
\[
{|f\left( y \right)-f\left( x \right)|} \leqslant 
M\left\| {y-x} \right\|_\infty .
\]
Such constant exists for all continuous functions $f(x)$ \eduard{on $[0,1]^n$}, so this  assumption is not restrictive.
Let us set the desired accuracy in terms of the objective as $\varepsilon$, i.e., our goal is to find a point $\hat{x}$ such that $f(\hat{x})-f^*\leqslant\varepsilon$, where $f^*$ is the \textit{global} minimum of $f$ on $\left[ {0,1} \right]^n$. For simplicity, we assume $\varepsilon$ to be equal to \eduard{$\nicefrac{1}{N}$} for some $N\in\mathbb{N}$. Consider a family of continuous non-convex objectives $f_k(x)$, $k=1,\ldots, N^n$, constructed as follows: we divide the hypercube $[0,1]^n$ into \eduard{$(\nicefrac{MN}{2})^n$} non-intersecting hypercubes $C_k$ with side length \eduard{$\nicefrac{2}{(NM)}$} and set
\[f_k(x) = \begin{cases} -M\text{dist}_\infty(x,\partial C_k),&\quad x\in C_k, \\
0,&\quad x\notin C_k, \end{cases}\]
where $\partial C_k$ is the boundary of $C_k$ and $\text{dist}_\infty(x,\partial C_k)$ is the distance between $x$ and $\partial C_k$ in the $\|\cdot\|_\infty$-norm. Each $f_k$ has a minimum value of exactly $-\varepsilon$ attained at the center of $C_k$, and the Lipschitz constant of $f_k$ is equal to $M$.

Any minimization method generating its trajectory based on the values of $f(x)$ and its derivatives at the points of the trajectory would need to sample a point from each $C_k$ to find an approximate minimum of each $f_k(x)$. This gives us a lower bound on the number of iterations required: $\Omega((MN)^n)=\Omega(M^n\varepsilon^{-n}).$ And this bound is attained by the algorithm which simply samples the objective values at the vertices of a uniform grid and returns the point with the smallest value. This demonstrates that it is practically impossible to solve a high-dimensional non-convex minimization problem with any reasonable accuracy unless some additional assumptions are introduced.

A similar complexity bound is proved in \cite{nesterov2012make} for finding a point $\hat{x}$ such that $\|\nabla f(\hat{x})\|_{\infty} \leqslant\varepsilon$ and $\|\hat{x}\|_{\infty}\leqslant R$. More precisely, for non-convex functions with Lipschitz continuous
Hessian, such that there exists at least one point $x^\ast$ with $\nabla f(x^\ast) = 0$ and $\|x^\ast\|_{\infty}\leqslant R$, the lower complexity bound is $\Omega\left(\left(\nicefrac{MR^2}{\eps}\right)^{\nicefrac{n}{2}}\right)$.

\subsection{Examples of Non-Convex Problems}

In this \eduard{subsection,} we make a non-extensive overview of non-convex problem formulations and applications where they arise, with a focus on tractable problems. 
One possible way to classify such non-convex problems is to divide them into \eduard{two groups:}
\begin{itemize}
    \item \hyperref[HiddenConvexity]{problems with hidden convexity or analytic solutions;}
    \item \hyperref[ProvableGlobalResults]{problems with provable global solution.}
\end{itemize}
Let us consider formulations of a few concrete problems in each of these classes.

\subsubsection{Problems with Hidden Convexity or Analytic Solutions}\label{HiddenConvexity}

Firstly, it is worth noting a \eduard{broad} class of classical non-convex \eduard{problems} that include linear-fractional programs, geometric programs, problems with two quadratic functions, handling convex equality constraints, convexifying constraint sets. Many such problems are equivalent to convex problems via a simple transformation such as convex relaxation and duality \cite{boyd2004convex}.
        
Next, a wide range of tasks in machine learning and statistics is reduced to eigenproblems. Among these problems are the following principal component analysis, classical multidimensional scaling, and other generalized eigenvalue problems \cite{charisopoulos2020entrywise}. 
        
In the context of non-convex optimization problems, one cannot but mention the class of combinatorial optimization problems as graph problems. Basically, most of these problems are NP-complete, but despite this, there are effective approaches and ways to solve them. Let us consider a closer look at the MAX-CUT problem. This is \eduard{a bright} example of convex reformulations. In some \eduard{problems,} the goal is to find a point with a value as small as possible (or as large as possible in the context of maximization problems), but whether this point is close to the global minimum is not that important. In this case, we can try to approximate the problem with a simpler one and show that the exact solution to the approximate problem corresponds to a good solution of the original problem. We will first illustrate this idea on the MAX-CUT problem
        \[
        \mathop {\max }\limits_{x\in \left\{ {-1,1} 
        \right\}^n} \left\{  f\left( x \right):=\frac{1}{2}\sum\limits_{i,j=1,1}^{n,n} {A_{ij} \left( {x_i 
        -x_j } \right)^2} \right\},
        \]
where $A=\left\| {A_{ij} } \right\|_{i,j=1,1}^{n,n} $ ($A=A^T)$. This is a discrete optimization problem. If we are interested only in the value of the functional and not in the cut itself, we can approximate this problem with a computationally tractable one. \iffalse \textit{Sometimes there exists such threshold, that to find the solution better than this threshold is much more difficult than to find the solution corresponds to this threshold. So let us restrict ourselves only by threshold solution.} \fi Let us introduce \eduard{matrix}
        \[
        L=\mbox{diag}\left\{ {\sum\limits_{j=1}^n {A_{ij} } } \right\}_{i=1}^n -A,
        \] which allows us to write
        \[
        f\left( x \right)=\left\langle {x,Lx} \right\rangle .
        \]
A simple observation: if $\varsigma$ is a random vector uniformly distributed on the Hamming cube $\left\{ {-1,1} \right\}^n$, then
        \[
        \mathbb{E}\left\langle {\varsigma ,L\varsigma } \right\rangle \ge 0.5\mathop {\max 
        }\limits_{x\in \left\{ {-1,1} \right\}^n} \left\langle {x,Lx} \right\rangle 
        .
        \]
In fact, we can do better \eduard{due to the construction of} Goemans and Williamson \cite{goemans1995improved}
    \[\mathop {\max }\limits_{x\in \left\{ {-1,1} \right\}^n} \left\langle {x,Lx} 
    \right\rangle =\mathop {\max }\limits_{x\in \left\{ {-1,1} \right\}^n} 
    \left\langle {L,xx^T} \right\rangle \leqslant \mathop {\max}\limits_{\scriptsize{\begin{array}{c}
     X\in S_+^n \\ 
     X_{ii} =1,\;i=1,...,n \\ 
     \end{array}}} \left\langle {L,X} \right\rangle.\] 
This is an SDP problem. Let $\Sigma $ be the solution of this SDP problem and let
    \[
    \xi \in N\left( {0,\Sigma } \right),
    \quad
    \varsigma =\mbox{sign}\left( \xi \right).
    \]
Then 
    \[
    E\left\langle {\varsigma ,L\varsigma } \right\rangle \ge \alpha_{GW}\mathop {\max 
    }\limits_{x\in \left\{ {-1,1} \right\}^n} \left\langle {x,Lx} \right\rangle 
    ,
    \]
where $\alpha_{GW}\approx0.878567$, and this constant is unimprovable provided that $\text{P}\neq\text{NP}$ and the Unique games conjecture is true \cite{khot2007optimal}.

Further, we would like to highlight the following subclasses of \eduard{non-convex} problems: non-convex proximal operators (Hard-thresholding \cite{blumensath2009iterative}, Potts minimization \cite{kiefer2020iterative}), discrete problems (Binary graph segmentation, Discrete Potts minimization, Nearly optimal $K$-means), infinite-dimensional problems (Smoothing splines, Locally adaptive regression splines, Reproducing kernel Hilbert spaces) and statistical problems.

Another important practical example we would like to mention in this part is Blind Deconvolution. Convolutional models arise in a wide range of problems in image processing and computer vision. The most basic convolutional data model -- blind deconvolution aims to recover a convolution kernel $a_0 \in \R^k$ and signal $x_0 \in \R^m$ from their convolution 
\[
    y = a_0 \circledast x_0,
\]
where $y \in \R^m$ and $\circledast$ is some kind of convolution. This problem is ill-posed in general –- there are infinitely many $(a_0, x_0)$ that convolve to produce $y$. To overcome this issue, some low dimensional priors about $a_0$ and $x_0$ are necessary. As a result, it is essential to use additional constraints and regularization terms. Different priors produce different \eduard{non-convex} optimization problems: Sparse Blind Deconvolution \cite{qu2019nonconvex}, Multi-channel Sparse Blind Deconvolution \cite{shi2020manifold}, Subspace blind deconvolution \cite{li2016identifiability}, Convolutional dictionary learning \cite{papyan2017convolutional}.
    
\eduard{The seen data in many settings in science and engineering are admixtures of several latent sources. Given the observations, we would normally wish to infer the latent sources as well as the admixture distribution. The non-negative matrix factorization (NMF) \cite{ge2015intersecting} mathematical framework offers a natural mathematical framework for modeling numerous mixing problems. In NMF, each row of observation matrix $M \in \R^{n\times m}$ corresponds to a data-point in $\mathbb{R}^m$. Next, the following assumptions are used: 1) there are $r$ latent sources, encoded by the unobserved matrix $W \in \mathbb{R}^{r\times m}$, and 2) each observed data-point can be rewritten as a linear combination of the $r$ sources, the weights of combination are defined via matrix $A \in \mathbb{R}^{n \times r}$. The goal is to find such representation of matrix $M$ that $M=AW$ with the entries of $M$, $A$ and $W$ being non-negative. The number $r$ is called the inner-dimension of the factorization, and the smallest possible $r$ is the nonnegative rank of $M$.}
    
Finally, in the part devoted to problems with Hidden Convexity or Analytical Solution we would like deal with Compressed Sensing and L1-optimization. A vector is said to be $s$-sparse if it has at most $s$ non-zero elements. Consider solving $Ax = b$ for $x$ where $A$ is an $n \times d$ matrix with $n < d$. The set of solutions to $Ax = b$ is a subspace. However, if we restrict ourselves to $s$-sparse solutions, under certain conditions on \eduard{$A$} there is a unique sparse solution \cite{blum2016foundations}. For instance, suppose that there were two $s$-sparse solutions $x_1$ and $x_2$. Then $x_1 -x_2$ would be a $2s$-sparse solution to the homogeneous system $Ax = 0$, which would imply that some $2s$ columns of $A$ are linearly dependent. Unless $A$ has $2$s linearly dependent columns, there can only be one $s$-sparse solution. 
    
There are many areas in which the problem is to find the unique sparse solution to a linear system. One is in plant breeding \cite{blum2016foundations}. Assume we are given a number of apple trees and the strength of some desirable feature of each tree. If we wish to determine which genes are responsible for the feature, we may formulate a system of linear equations $Ax=b$ in which each row of the  matrix $A$ corresponds to a tree and each column corresponds to a position on the genome. The vector $b$ corresponds to the strength of the desired feature in each tree. The solution $x$ tells us the positions on the genome corresponding to the genes that account for the feature.
        
The problem of finding a sparse solution can be stated as the optimization problem \[\mathop {\min }\limits_{Ax=b} \left\| x \right\|_0,\] where $\|x\|_0$ is \eduard{the number} of non-zero coordinates of $x$. This is an NP-hard problem, but it may sometimes be replaced by the convex problem \[\mathop {\min}\limits_{Ax=b} \left\| x \right\|_1.\]
    
What are the sufficient conditions for
    \[
    \mathop {\min }\limits_{Ax=b}  \left\| x \right\|_0 
    \quad
    \Leftrightarrow 
    \quad
    \mathop {\min }\limits_{Ax=b} \left\| x \right\|_1  ?
    \]
A matrix $A$ is said to satisfy \textit{the $s$-restricted isometry property}  if for any $s$-sparse $x$ there exists $\delta_s$ such that
    \[\left( {1-\delta _s } \right)\left\| x \right\|_2^2 \leqslant \left\| {Ax} 
    \right\|_2^2 \leqslant \left( {1+\delta _s } \right)\left\| x \right\|_2^2. \] 
The following theorems give sufficient conditions for the equivalence mentioned above to hold \cite{blum2016foundations,candes2005decoding}.
\begin{theorem}
    Suppose $A$ satisfies the $s$-restricted isometry property with $\delta_{s+1}\leqslant\frac{1}{10\sqrt{s}}.$
    Suppose $x_0$ is $s$-sparse and satisfies $Ax_0 = b$. \eduard{Then,} $x_0$ is the unique minimum 1-norm solution to $Ax = b$.
\end{theorem} 
\begin{theorem}
    Suppose $A$ satisfies the $k$-restricted isometry property for $k\in\{s,2s,3s\}$ with $\delta_{s}+\delta_{2s}+\delta_{3s}\leqslant 1.$
    Suppose $x_0$ is $s$-sparse and satisfies $Ax_0 = b$. \eduard{Then,} $x_0$ is the unique minimum 1-norm solution to $Ax = b$.
\end{theorem} 

Such results demonstrate the importance of matrices satisfying the restricted isometry property for practice. Fortunately, there is an easy way to obtain such matrices \cite{baraniuk2008simple}.
\begin{theorem}
    Suppose $A$ is an $d\times n$ matrix with elements sampled from the Gaussian distribution $\mathcal{N}(0,1/d)$. \eduard{Then,} $A$ satisfies the $s$-restricted isometry property for $s<d$ with $0<\delta_s<1$ with probability $p_s$ satisfying
        \[p_s\geqslant 1-2(12/\delta_s)^s\exp{\left(-\frac{3\delta_s^2-\delta_s^3}{48}d\right)}.\]
\end{theorem}

\subsubsection{Problems with Convergence Results}\label{ProvableGlobalResults}

In this section, we would like to give examples of non-convex optimization problems for which there are methods with proven convergence results. \eduard{We} start with the \textbf{Phase retrieval problem}. The phase retrieval problem has been a topic of study from at least the early 1980s. It is the recovery of a function given the magnitude of its Fourier transform. This problem could be found in various engineering and scientific applications such as optical imaging, electron microscopy, and crystallography, etc. \cite{shechtman2015phase}. We recover a $d$-dimensional signal vector $x^\ast \in \mathbb{C}^d$ from its phaseless measurements
    \[
    y_k=|\langle a_k,x \rangle|^2, \quad k=1,\dots,M,
    \]
    with $a_k$ denoting the measurement vectors. As a result, the phase-retrieval problem can be formulated as the following least squares problem or empirical risk
    minimization
    \[
    \min\limits_{x}\sum\limits_{k=1}^M\left(y_k - |\langle a_k,x \rangle|^2\right)^2.
    \]
    This problem is well-motivated by practical concerns, but unfortunately, this is a non-convex problem, and it is not clear how to find a global minimum even if one exists. In recent literature, there are various approaches to handle this problem \cite{wu2020hadamard, tan2019online, chen2019gradient}, also, algorithms with the provable convergence results were presented in the following papers \cite{candes2015phase, yang2019misspecified}.
    
In the context of non-convex optimization problems with proven convergence result, one cannot but mention \textbf{Low-Rank Matrix Completion}. There are related problems: matrix completion and matrix sensing \cite{bhojanapalli2016dropping}, which are present in big data problems with incompleteness and other machine learning problems. We would like to draw attention to the exact low-rank matrix completion. Given a matrix $Y \in \mathbb{R}^{n \times n}$, partially observed, over a set of indices $\Omega \subseteq \{1,\dots,n\}^2$. Consider the problem of finding the lowest-rank matrix matching $X$ on the observed set
    \[
    \min\limits_{X} \; \mathrm{rank}\left( X\right)
    \]
    \[
    \mathrm{s.t.} \quad X_{ij} = Y_{ij}, \quad (i,j) \in \Omega.
    \]
    This is a non-convex problem having a natural convex relaxation
    \[
    \min\limits_{X} \; \|X\|_{\mathrm{tr}}
    \]
    \[
    \mathrm{s.t.} \quad X_{ij} = Y_{ij}, \quad (i,j) \in \Omega
    \]
    In the paper \cite{jain2013low} the first results of global optimality of alternating minimization were obtained for matrix completion and the related problem of matrix sensing. Proofs of (nearly) linear convergence of gradient descent for Phase retrieval, Matrix completion, Blind deconvolution can be found in the article \cite{ma2018implicit}. Under some assumptions, it can be shown that the solution to the convex problem is exactly equal to the solution to the non-convex problem, with high probability over the sampling model \cite{candes2010power, candes2009exact}. So, this problem can also be attributed to statistical problems with hidden convexity. \eduard{Moreover, we we emphasise another relevant problem called \textbf{Low-Rank Matrix Recovery}. This problem is also known to be non-convex but under some assumptions has no spurious local minima (see \cite{zhang2021sharp,zhang2021general} and references therein).}

    \noindent \textbf{Deep Learning.} In the era of AI, training of the deep neural networks \cite{Goodfellow-et-al-2016} is one of the most popular optimization problems with enormous amount of applications, e.g.,  \cite{kniaz2021adversarial,rezanov2021deep,khritankov2021hidden,kuderov2021planning,surazhevsky2021noise-assisted,ilyuhin2020recognition,demin2021necessary,gorodetskiy2020delta,skrynnik2021forgetful,paparoditis2020wire}. The simplest example of such problem \cite{sun2019optimization} is training fully connected neural network for supervised learning problem
    \begin{equation}
        \min\limits_{\stackrel{W = (W_1,\ldots,W_L)}{W_i \in \R^{n_i\times n_{i-1}}, i=1,\ldots,L}}\left\{f(W) := \frac{1}{m}\sum\limits_{i=1}^m\ell(y_i, f_{x_i}(W))\right\}, \label{eq:fully_connected_NN_problem}
    \end{equation}
    where $\{(x_i,y_i)\}_{i=1}^m$, $x_i\in \R^{n_0}$, $y_i\in\R^{n_y}$ are training data points, $W = (W_1,\ldots, W_L)$ are weights of the model, $L$ is number of fully connected layers, $\ell(\cdot,\cdot)$ is a loss function, e.g., quadratic loss or logistic loss, and
    \begin{equation*}
        f_{x_i}(W) = W_L\phi\left(W_{L-1}\phi \ldots \phi\left(W_2\phi\left(W_1x_i\right)\right)\right),
    \end{equation*}
    where $\phi$ is a scalar\footnote{By $\phi(a)$ where $a = (a_1,\ldots,a_n)^\top \in \R^n$ is multidimensional vector we mean vector $(\phi(a_1),\ldots,\phi(a_n))^\top$.} function called an activation function.
    
    In general, training neural networks is NP-complete problem \cite{blum1989training}. Deep neural networks have bad local minima both for non-smooth activation functions \cite{swirszcz2016local,safran2018spurious} and smooth ones \cite{liang2018understanding,yun2018small} as well as flat saddles \cite{vidal2017mathematics}. Nevertheless, there exist positive results about training neural networks. First of all, under different assumptions it was shown that all local minima are global for 1-layer neural networks \cite{soltanolkotabi2018theoretical,haeffele2017global,feizi2017porcupine}. Next, one can show that GD/SGD converge under some assumptions to global minimum for linear networks \cite{arora2018convergence,ji2019gradient,shin2019effects} and sufficiently wide over-parameterized networks \cite{allen2019convergence}. The detailed summary of recent advances in optimization for deep learning can be found in \cite{sun2019optimization}.

\subsubsection{Geometry of non-convex optimization problems}
In one of the latest \eduard{survey} \cite{zhang2020symmetry}, the authors to distinguish a class of tractable non-convex problems, which have certain properties of symmetry. They highlight non-convex optimization problems with rotational symmetry and discrete symmetry. Problems with rotational symmetry include the previously described phase retrieval and related problems in low-rank matrix factorization and recovery. It turns out that the blind deconvolution and tensor decomposition problems have discrete symmetry.

\section{Deterministic First-Order Methods}\label{sec:deterministic_first_order}

In this section we focus on the following optimization problem 
\begin{equation}
\label{eq:problem_statement}
\mathop {\min }\limits_{x\in Q \subseteq {\mathbb{R}}^n}f\left( x \right), 
\end{equation}
where $Q$ is a simple, closed, convex, set, and $f$ is continuously differentiable function. The simplest method for this kind of problems is projected gradient descent, which can be motivated by a simple continuous-time dynamics. For simplicity we start with the unconstrained case with $Q={\mathbb{R}}^n$. 

\subsection{Unconstrained Minimization}

In the case $Q={\mathbb{R}}^n$, the trajectory of the continuous-time gradient method is the solution to the differential equation $\dot{x}=-\nabla f\left( x(t) \right)$.
It is easy to see that $W\left( x \right)=f\left( x(t) \right)$ is a Lyapunov function for this dynamical system. Indeed, 
\[
\begin{array}{c}
 \frac{dW\left( {x\left( t \right)} \right)}{dt}{\kern 1pt}={\kern 
1pt}\left\langle {\nabla f\left( {x\left( t \right)} \right),\frac{dx\left( 
t \right)}{dt}} \right\rangle {\kern 1pt}=  
\left\langle {\nabla f\left( {x\left( t \right)} 
\right),-\nabla f\left( {x\left( t \right)} \right)} \right\rangle {\kern 
1pt}=-\left\| {\nabla f\left( {x\left( t \right)} \right)} \right\|_2^2 
{\kern 1pt}\leqslant 0. \\ 
 \end{array}
\]
This implies the convergence of the continuous-time gradient descent method to a stationary point.

The classic gradient descent method is then the Euler discretization of the above dynamics and has the form \cite{polyak1987introduction}
\[
x^{k+1}=x^k-h_k\nabla f\left( {x^k} \right),
\]
where $h_k\geq 0$ is the stepsize of the method. One of the main assumptions in this setting is that the function $f$ is $L$-smooth, or, which is the same, its gradient is Lipschitz-continuous, i.e., for some starting point $x^0$, \[\forall x,y\in \left\{ {x\in {\mathbb{R}}^n:\;f\left( x \right)\leqslant f\left( {x^0} \right)} \right\}\quad\left\| {\nabla f\left( y \right)-\nabla f\left( x \right)} \right\|_2 \leqslant 
L\left\| {y-x} \right\|_2.
\]
Then the stepsize $h=1 /L $ guarantees
\[
f( {x^{k+1}} )\leqslant f( {x^k} )-\frac{1}{2L}\left\| 
{\nabla f( {x^k} )} \right\|_2^2.
\]
Summing up these inequalities, we obtain
\[f(x^N)-f(x^0)\leqslant-\frac{1}{2L}\sum_{k=0}^{N-1}\|\nabla f(x^k)\|^2_2\leqslant-\frac{N}{2L}\min\limits_{k=0,...,N-1}\|\nabla f(x^k)\|^2_2.\]
Define $\pd{f_*}=\inf\limits_{x\in\mathbb{R}^n}f(x)$ and assume that this value is finite. Then
\begin{equation}
\label{eq:nonconv_compl_grad}
\min\limits_{k=0,...,N-1}\|\nabla f(x^k)\|^2_2 \leqslant \frac{2L(f(x^0)-\pd{f_*})}{N}.     
\end{equation}
This proves that the complexity of finding an approximate stationary point, i.e. a point $\hat{x}$ such that $\|\nabla f(\hat{x})\|_2 \leqslant\varepsilon$ %\alexander{I guess $\| \| = \| \|_2$?}, 
is $O\left( \frac{L\left( {f\left( {x_0} \right)-\pd{f_*}} 
\right)}{\varepsilon ^2}\right)$. This iteration complexity of finding an $\eps$-stationary point $N\sim \varepsilon ^{-2}$ is unimprovable in terms of its dependence on $\varepsilon$ and $L$ for an arbitrary first-order method applied to minimization of an $L$-smooth objective. 

On the one hand this bound is much better than the exponential in the dimension bound for finding the global minimum, which was derived in Subsection \ref{S:lower bound}. On the other hand we can guarantee only an approximate stationary point, which could be a saddle-point or even a maximum.
This can be illustrated by the example of minimization of the following objective \cite{nesterov2018lectures}
\[
f\left( {x_1,x_2 } \right)=\frac{1}{2}\left( {x_1 } 
\right)^2+\frac{1}{2}\left( {x_2 } \right)^4-\frac{1}{2}\left( {x_2 } 
\right)^2.
\]
If we set $x^0=\left( {1,0} \right)^T$, then 
$x^k$ converges to $ \left( 
{0,0} \right)^T$ as $k \to \infty$, which is a saddle-point.
The good news here is that gradient descent can be perturbed by adding some noise in the iterates in such a way that it converged to a local minimum for almost all initial points and escapes saddle-points \cite{jin2017how}. 

It is important to note that, under additional smoothness assumptions that  higher-order derivatives of the objective are Lipschitz continuous, i.e.
\[
\forall x,y\in \left\{ {x\in {\mathbb{R}}^n:\;f\left( x \right)\leqslant f\left( {x^0} \right)} \right\}\quad\left\| {\nabla^p f\left( y \right)-\nabla^p f\left( x \right)} \right\|_2 \leqslant 
L_p\left\| {y-x} \right\|_2,
\] 
\cite{carmon2019lowerI,carmon2019lowerII} obtain several lower complexity bounds for finding an approximate stationary point. If this inequality holds for $p\in\{1,2\}$, the lower bound becomes $\eps^{-\frac{12}{7}}$, and the additional assumption that the same holds for $p=3$ gives the lower bound to $\eps^{-\frac{8}{5}}$. Surprisingly, Lipschitz continuity of derivatives of order 4 and higher gives the same lower complexity bound.

\subsection{Incorporating Simple Constraints}
It is possible to generalize gradient method for the setting of composite optimization with simple convex constraints, i.e. for the problem
\begin{equation}
\label{eq:PrStateInit}
\min_{x \in Q } \{ F(x) := f(x) + \psi(x)\},
\end{equation}
where $Q$ is a closed convex set, $\psi(x)$ is a simple convex function, e.g. $\|x\|_1$, and $f$ is $L$-smooth function. The standard approach for such problems uses {\it prox-function} $d(x)$ which is continuously differentiable and strongly convex on $Q$, i.e. $d(y)-d(x) -\la \nabla d(x) ,y-x \ra \geq \frac12\|y-x\|^2$ for any $x, y \in Q$.
We define also the corresponding {\it Bregman divergence} $V[z] (x) = d(x) - d(z) - \la d'(z), x - z \ra$, $x,z \in Q$. Then the step of the gradient method from a point $x$ with stepsize $h$ is generalized \cite{nesterov2018lectures,ghadimi2016mini-batch} to
\[
x^+ = \arg \min_{u\in Q}\left\{\la \nabla f(x),u \ra + \frac{1}{h} V[x](u) +\psi(u) 		 \right\}, 
\]
which in the simplest case $\psi(x)\equiv 0$, $d(x)=\frac{1}{2}\|x\|_2^2$, $V[z] (x)=\frac{1}{2}\|x-z\|_2^2$, $Q=\R^n$ coincides with the step of the gradient method. This generalized gradient step leads to a generalized gradient, which is usually referred to as gradient mapping \cite{nesterov2018lectures,ghadimi2016mini-batch}
$g_Q(x)=\frac{1}{h}(x-x^+)$. In this setting, the authors of \cite{ghadimi2016mini-batch} prove that 
\[
\min\limits_{k=0,...,N-1}\|g_Q(x^k)\|^2 \leqslant \frac{2L(F(x^0)-\pd{F_*})}{N} 
\]
if $h=1/L$. Here $\pd{F_*}$ is a lower bound for $F(x)$.
In the described above simple situation this bound coincides with the bound \eqref{eq:nonconv_compl_grad}. The authors of \cite{dang2015stochastic} prove that if  $\|g_Q(x)\|\leqslant \eps$, then $x^+$ is an approximately stationary point of the problem. More precisely, there exist $p \in \partial \psi(x^+)$ such that 
\[
\nabla f(x^+) + p \in -\mathcal{N}_Q(x^+)+\mathcal{B}((1+L(d)\eps),
\]
where $\mathcal{N}_Q(x^+)$ is the normal cone of $Q$ at the point $x^+$,  $B(r)=\{v\in \R^n: \|v\|_* \leqslant r\}$ -- ball in the dual space defined by the conjugate norm, and it is assumed that $d$ is $L(d)$-smooth.
Note that there is no contradiction with the exponential lower bound given in the end of Subsection \ref{S:lower bound} since non-necessarily the obtained point $x^+$ has small norm of the gradient.

This approach was further generalized in \cite{bogolubsky2016learning,dvurechensky2017gradient,gasnikov2018power} for the case of optimization with inexact oracle for the function $f$.
\begin{Def}
\label{D:IO}
We say that a function $f(x)$ is equipped with an \textit{inexact first-order oracle} on a set $X$ if there exists $\delta_u > 0$ and at any point $x \in X$ for any number $\delta_c > 0$ there exists a constant $L(\delta_c) \in (0, +\infty)$ and one can calculate  $\tf(x,\delta_c,\delta_u) \in \R$ and $\tg(x,\delta_c,\delta_u) \in \R^n$ satisfying
\begin{align}
	&|f(x) - \tf(x,\delta_c,\delta_u)| \leqslant\delta_c+\delta_u, 
	\label{eq:dL_or_def_2} \\
	&f(y)-(\tf(x,\delta_c,\delta_u)  - \la\tg(x,\delta_c,\delta_u) ,y-x \ra) \leqslant\frac{L(\delta_c)}{2}\|x-y\|^2 + \delta_c+\delta_u, \quad \forall y \in Q.
  \label{eq:dL_or_def_1}
\end{align}
\end{Def}
In this definition, $\delta_c$ represents the error of the oracle, which we can control and make as small as we would like to. On the opposite, $\delta_u$ represents the error, which we can not control. The proposed for this setting method in \cite{dvurechensky2017gradient} is adaptive to the constant $L$, works under inexact calculation of the point $x^+$, and covers several different settings.  
In particular, smooth functions with H\"older-continuous, i.e. satisfying, for some $\nu\in [0,1]$, $\|\nabla f(x) - \nabla f(y) \|_* \leqslant L_{\nu} \|x-y\|^{\nu}, \forall x,y \in Q$ gradient satisfy this definition with $\delta_u = 0$ and 
\begin{equation}
\label{eq:Lofd}
L(\delta_c) = \left( \frac{1-\nu}{1+\nu} \cdot \frac{2}{\delta_c} \right)^{\frac{1-\nu}{1+\nu}} L_{\nu}^{\frac{2}{1+\nu}}.
\end{equation} 
As a corollary of the general method, \cite{dvurechensky2017gradient} propose a universal method for such problems, which does not require the knowledge of the constants $\nu, L_\nu$ and gives the following convergence rate
\begin{equation}
\min\limits_{k=0,...,N-1}\left\|g_Q(x_k)\right\|^2  \leqslant2^{\frac{1+3\nu}{2\nu}} \left( \frac{1-\nu}{1+\nu} \cdot \frac{40}{\eps} \right)^{\frac{1-\nu}{2\nu}} L_{\nu}^{\frac{1}{\nu}} \left(\frac{F(x^0) - \pd{F_*}}{N}\right) + 
%(4\delta_u + \delta_{pu}) \right) + 
\frac{\eps}{2},
\label{eq:pg_la_2_rate_2}
\end{equation}
or the following complexity estimate $\frac{L_{\nu}^{\frac{1}{\nu}}(F(x^0)-\pd{F_*})}{\eps^{\frac{1+3\nu}{2\nu}}}$ to find $\left\|g_Q(x^k)\right\|\leqslant \eps$. %\alexander{???}
\pd{Inexact oracle models for convex optimization can be useful in non-convex optimization since in some settings a non-convex problem can be considered as a convex problem with inexact oracle \cite{stonyakin2019gradient,stonyakin2020inexact}.}

\subsection{Incorporating Momentum for Acceleration}

The considered above dynamical system $ \dot{x}=-\nabla f\left( x(t) \right)$ does not have any mechanical intuition behind it. In \cite{polyak1964some} the author proposed to consider the following dynamics
\[
\mu \ddot{x}(t) = -\nabla f\left( x(t) \right)-p \dot{x}(t).
\]
One of the ways to discretize it gives the so called heavy-ball method
\[
x^{k+1}=x^k-h\nabla f( {x^k} )+\beta ( {x^k-x^{k-1}} 
),
\]
where $h>0$ is the stepsize and $\beta>0$ is the momentum parameter. Due to the momentum term $\beta \left( {x^k-x^{k-1}} \right)$ the method avoids zigzagging for ill-conditioned problems, which leads to significant efficiency in practice, especially in training neural networks.
Despite practical efficiency, the theoretical guarantee for this method is no better than for the gradient method. In particular, \cite{griewank1981generalized} considers the dynamical system
\[
\mu (t) \ddot{x}(t) = -\nabla f\left( x(t) \right)-p (t) \dot{x}(t),
\]
where $\mu \left( t \right)\sim \left(f\left( {x\left( t \right)} \right)- c \right)$, $c$ is an upper bound on the global minimum of $f(x)$,  and $p\left( t \right)=F\left( {\nabla 
f\left( {x\left( t \right)} \right)} \right)$. With a special choice of 
$F\left(\cdot\right)$, they show that $x\left( t \right)$ converges to 
a local minimizer $x^{loc}$ such that $f\left( {x^{loc}} \right)\leqslant c$ as $t\to +\infty$. 
In \cite{diakonikolas2019generalized} it is shown that for a discretization of a further generalization of the heavy-ball method one may guarantee
\[
\mathop {\min }\limits_{k=1,...,N} \| {\nabla f( {x^k} )} 
\|_2^2 \leqslant \frac{2L(f(x^0)-\pd{f_*})}{N},
\]
which coincides with the bound \eqref{eq:nonconv_compl_grad} for the gradient method.

A different type of momentum was proposed in \cite{nesterov1983method} for convex optimization, which led to the Nesterov's accelerated gradient method
\[
x^1=x^0-h\nabla f\left( {x^0} \right),
\]
\[
x^{k+1}=x^k-h\nabla f( {x^k+\beta _k ( {x^k-x^{k-1}} )} 
)+\beta _k ( {x^k-x^{k-1}} ).
\]
The difference with the heavy-ball method is that the gradient is calculated in the extrapolated point. This idea has been very fruitful and allowed to obtain many accelerated algorithms for convex optimization.
A variant of this method with a special choice of the stepsize $h$ and momentum term $\beta_k$ was shown in \cite{ghadimi2016accelerated} to have the same convergence rate \eqref{eq:nonconv_compl_grad} as the gradient method. This was further extended in \cite{ghadimi2019generalized} for the case of objective with H\"older-continuous gradients to obtain a bound $\frac{L_{\nu}^{\frac{1}{\nu}}(F(x^0)-\pd{F_*})}{\eps^{\frac{1+3\nu}{2\nu}}}$ to find $\left\|g_Q(x^k)\right\|\leqslant \eps$ in the general setting of composite optimization problem \pd{\eqref{eq:PrStateInit}} with simple constraints. Importantly, this method is universal and uniform, which means that it has best possible convergence rates for convex and non-convex problems without knowing whether the problem is convex or not and without knowing its smoothness parameters such as H\"older exponent and H\"older constant.

It is possible to combine this idea with the idea of line-search, i.e. minimization in the direction of the step. The papers \cite{guminov2019accelerated,nesterov2020primal-dual} propose a modification of the accelerated gradient method which is listed as Algorithm \ref{AGMsDR}. Instead of explicitly defining the stepsize $h$ and the momentum term $\beta$, this method uses full one-dimensional relaxation and local information. This makes this method parameter-free and uniform for convex and non-convex smooth optimization by providing optimal complexity bound for the convex and non-convex case. At the same time, inexact line-search is possible and its sufficient accuracy for achieving the desired accuracy is estimated. 
This method shares some similarities with nonlinear conjugate gradient methods which were analyzed in \cite{nemirovski1982orth}.

\begin{algorithm}[!h]
\caption{Accelerated Gradient Method with Small-Dimensional Relaxation (AGMsDR)}
\label{AGMsDR}
\begin{algorithmic}[1]
\ENSURE $x^k$
\STATE Set $k = 0$, $A_0=0$, $x^0 = v^0$, $\psi_0(x) = V[x^0](x)$
\FOR{$k \geqslant 0$}
\STATE \begin{equation}
    \label{eq:beta_k_y_k_def}
    \beta_k = \arg \min_{\beta \in \left[0, 1 \right]} f\left(v^k + \beta (x^k - v^k)\right), \quad y^k = v^k + \beta_k (x^k - v^k).
\end{equation}
\STATE 
Let $(\nabla f(y^k))^{\#}$ be such that $\langle \nabla f(y^k), (\nabla f(y^k))^{\#} \rangle = \|\nabla f(y^k) \|_*^2$ and $\|(\nabla f(y^k))^{\#} \|^2 = 1$.
\begin{equation}
\label{eq:xkp1_opt_b}
h_{k+1} = \arg \min_{h \geqslant 0} f\left(y^k - h(\nabla f(y^k))^{\#}\right), \quad x^{k+1} = y^  k- h_{k+1}(\nabla f(y^k))^{\#}.
\end{equation}
Find $a_{k+1}$ from equation $f(y^k) - \frac{a_{k+1}^2}{2(A_{k} + a_{k+1})} \|\nabla f(y^k) \|_*^2 = f(x^{k+1})$.
\STATE  Set $A_{k+1} = A_{k} + a_{k+1}$.  
\STATE  Set $\psi_{k+1}(x) =  \psi_{k}(x) + a_{k+1}\{f(y^k) + \langle \nabla f(y^k), x - y^k \rangle\}$.
\STATE $v^{k+1} = \arg \min_{x \in \R^n} \psi_{k+1}(x)$,  $k = k + 1$
\ENDFOR
\end{algorithmic}
\end{algorithm}

The above idea was further extended in \cite{guminov2021combination} where an accelerated alternating minimization method was proposed and analyzed for convex and non-convex problems. The main assumption is that the set of coordinates is divided into \pd{$\bar{n}$} disjoint subsets (blocks) $I_p$, $p\in\{1,\ldots,\pd{\bar{n}}\}$ and minimization in each block when the other variables are freezed can be made explicitly. The resulting accelerated alternating minimization algorithm is listed as Algorithm \ref{AAM-2}. 
This method is also parameter-free and uniform for convex and non-convex smooth optimization with optimal complexity bound for the convex and non-convex case.

\begin{algorithm}[H]
\caption{Accelerated Alternating Minimization (AAM)}
\label{AAM-2}
\begin{algorithmic}[1]
   \REQUIRE Starting point $x^0$.
    \ENSURE $x^k$
   \STATE Set $A_0=0$, $x^0 = v^0$.
   \FOR{$k \geqslant 0$}
	\STATE Set $\beta_k = \arg \min\limits_{\beta\in [0,1]} f\left(x^k + \beta (v^k - x^k)\right)$
	\STATE Set $y^k = x^k + \beta_k (v^k - x^k)\quad $
    \STATE Choose $i_k=\arg \max\limits_{i\in\{1,\ldots,\pd{\bar{n}}\}} \|\nabla_i f(y^k)\|_2^2$
\STATE Set $x^{k+1}=\arg \min\limits_{x\in S_{i_k}(y^k)} f(x)$, i.e. minimize $f$ in the corresponding block.
\STATE Find $a_{k+1}$, $A_{k+1} = A_{k} + a_{k+1}$ from \[f(y^k)-\frac{a_{k+1}^2}{2A_{k+1}}\|\nabla f(y^k)\|_2^2=f(x^{k+1})\]\\
\STATE Set $v^{k+1} = v^{k}-a_{k+1}\nabla f(y^k)$ 
\ENDFOR
\end{algorithmic}
\end{algorithm}
\pd{The sequence $y^k$ of this algorithm satisfies
\[
\mathop {\min }\limits_{k=1,...,N} \| {\nabla f( {y^k} )} 
\|_2^2 \leqslant \frac{2\bar{n}L(f(x^0)-f_*)}{N},
\]
i.e. there is an additional multiplier $M$ -- number of blocks.
If the function turns out to be convex, then the same method generates the sequence $x^k$ which gives the decay of the objective similar to accelerated gradient method:
\begin{align*}
    f(x^{k})-f(x^*) \leqslant \frac{2\bar{n}L\|x^0-x^*\|_2^2}{N^2},
\end{align*}
where $x^*$ is the closest to $x^0$ global minimizer.
}

By exploiting the idea of Nesterov's acceleration and combining it with the notion of negative curvature, the authors of \cite{carmon2017convex} manage to accelerate first-order methods for non-convex optimization under additional assumptions that second and third derivatives are Lipschitz continuous. More precisely, if $L$-smooth function has also Lipschitz continuous Hessian, they obtain complexity $O\left( \varepsilon ^{-7/4}\log(1/\eps)\right)$ to find a point $\hat{x}$ such that $\|\nabla f(\hat{x})\|_2 \leqslant\varepsilon$. Assuming additionally that the third derivative is Lipschitz, this bound is improved to  $O\left( \varepsilon ^{-5/3}\log(1/\eps)\right)$.

\section{Stochastic First-Order Methods}\label{sec:stoch_methods}
In this section, we consider the same problem as in Section~\ref{sec:deterministic_first_order}:
\begin{equation}
    \min\limits_{x\in\R^n} f(x), \label{eq:main_pr_non_cvx}
\end{equation}
where function $f$ is a general non-convex $L$-smooth function with the uniform lower bound $f_*$, i.e., it is differentiable and
\begin{eqnarray}
    f(x) &\ge& f_* \quad \forall x\in\R^n,\label{eq:unif_lower_bound}\\
    \|\nabla f(x) - \nabla f(y)\|_2 &\le& L\|x-y\|_2\quad \forall x,y\in\R^n.\label{eq:L_smoothness}
\end{eqnarray}
We are interested in two particular cases: expectation minimization
\begin{equation}
    f(x) = \EE_\xi[f(x,\xi)], \label{eq:f_expectation}
\end{equation}
and finite-sum minimization
\begin{equation}
    f(x) = \frac{1}{m}\sum\limits_{i=1}^m f_i(x). \label{eq:f_finite_sum}
\end{equation}
Such problems usually arise in applications of (deep) machine learning \cite{Goodfellow-et-al-2016,sun2019optimization} and mathematical statistics \cite{spokoiny2012parametric}, and typically they are solved via stochastic first-order methods.

In general, the best one can expect to achieve is an approximate stationary point \cite{vavasis1993black,arjevani2019lower}. To be specific, for this class of problems stochastic first-order methods in the worst case can only find such point $\hat x$ that
\begin{equation}
    \EE\left[\|\nabla f(\hat x)\|_2^2\right] \le \varepsilon^2 \label{eq:epsilon_stationary_point}
\end{equation}
For simplicity, we will call the point $\hat x$ as $\varepsilon$-stationary point, but mean by this that inequality \eqref{eq:epsilon_stationary_point} holds.

Below we summarize recent results about finding $\varepsilon$-stationary point using stochastic first-order methods. We start with presenting the general and unified approach to analyze optimal deterministic and stochastic first-order methods for objectives of types \eqref{eq:f_expectation} and \eqref{eq:f_finite_sum} in the general settings. After that, we consider $3$ big classes of stochastic first-order methods with convergence guarantees: SGD and its variants, variance reduced methods, and adaptive stochastic methods. 

\subsection{General View on Optimal Deterministic and Stochastic First-Order Methods for Non-Convex Optimization}\label{sec:general_view}
Assume that at each point $x$, we have access to the estimator $g(x)$ of the gradient $\nabla f(x)$. For now, it is not important to specify what properties $g(x)$ satisfies. In these settings one can use Algorithm~\ref{alg:opt_non_cvx_methods} in order to find $\varepsilon$-stationary point.
\begin{algorithm}[h!]
    \caption{General scheme of the optimal first-order method for non-convex optimization}
    \label{alg:opt_non_cvx_methods}
    \begin{algorithmic}[1]
        \REQUIRE learning rates $\{h_k\}_{k\ge 0}$ satisfying $h_k \le \frac{1}{2L}$, starting point $x^0\in\R^n$, stopping criterion $C$
        \FOR{ $k=0,1,2,\ldots$ }
        \STATE Get $g^k = g(x^k)$
        \IF{$C$ holds}
        \STATE $x^N = x^k$
        \STATE break
        \ELSE
        \STATE $x^{k+1} = x^k - h_kg^k$
        \ENDIF
        \ENDFOR
        \RETURN $x^{N}$
    \end{algorithmic}
\end{algorithm}

Below we derive preliminary inequalities playing the central role in the analysis of optimal (stochastic) first-order algorithms. From $L$-smoothness of $f$ we have
\begin{eqnarray}
    f(x^{k+1}) &\le& f(x^k) + \langle \nabla f(x^k), x^{k+1} - x^k\rangle + \frac{L}{2}\|x^{k+1}-x^k\|_2^2\notag\\
    &=& f(x^k) + \langle g^k, x^{k+1} - x^k\rangle + \langle \nabla f(x^k) - g^k, x^{k+1} - x^k\rangle + \frac{L}{2}\|x^{k+1}-x^k\|_2^2\notag\\
    &\le& f(x^k) - h_k\|g^k\|_2^2 + h_k\|\nabla f(x^k) - g^k\|_2^2 + \left(\frac{1}{4h_k}+\frac{L}{2}\right)\|x^{k+1}-x^k\|_2^2, \notag
\end{eqnarray}
where in the last inequality we use Fenchel--Young inequality: $\langle a,b \rangle \le \frac{1}{2\alpha}\|a\|_2^2 + \frac{\alpha}{2}\|b\|_2^2$ with $a = \nabla f(x^k) - g^k$, $b = x^{k+1} - x^k$ and $\alpha = \frac{1}{2h_k}$. Since $h_k \le \frac{1}{2L}$ and $x^{k+1} = x^k - h_k g^k$ we can continue our derivations:
\begin{eqnarray}
    f(x^{k+1}) &\le& f(x^k) - \frac{h_k}{2}\|g^k\|_2^2 + h_k\|\nabla f(x^k) - g^k\|_2^2. \label{eq:main_ineq_non_cvx_opt_methods}
\end{eqnarray}
Now it is crucial to specify what we need to assume about $g(x)$. We emphasize that all $3$ cases considered below are based on the tight bounds for $\|\nabla f(x^k) - g^k\|_2^2$ or its expectation.

\subsubsection{Deterministic Case}
In this case we assume that for all $x\in \R^n$ we have an access to such $g(x)$ that
\begin{equation}
    \|g(x) - \nabla f(x)\|_2^2 \le \frac{\varepsilon^2}{10}. \label{eq:deterministic_assumption}
\end{equation}
In other words, $g(x)$ is good enough approximation of $\nabla f(x)$. Consider the stopping criterion $C =\left\{\|g^k\|_2^2 \le \frac{2\varepsilon^2}{5}\right\}$ and let $h_k = \frac{1}{2L}$ for all $k\ge 0$. First of all, if Algorithm~\ref{alg:opt_non_cvx_methods} stops, then $\|g^N\|_2 \le \frac{4\varepsilon^2}{10}$ and $x^N$ satisfies
\begin{equation*}
    \|\nabla f(x^N)\|_2^2 = \|\nabla f(x^N) - g^N + g^N\|_2^2 \le 2\|\nabla f(x^N) - g^N\|_2^2 + 2\|g^N\|_2^2 \overset{\eqref{eq:deterministic_assumption}}{\le} \varepsilon^2.
\end{equation*}
Next, we derive an upper bound for such $N$ that Algorithm~\ref{alg:opt_non_cvx_methods} stops after $N$ iterations. Assume that, after $N$ iterations the method has not stopped. Then for all $k = 0,1,\ldots,T$ we have
\begin{equation*}
    f(x^{k+1}) \overset{\eqref{eq:main_ineq_non_cvx_opt_methods},\eqref{eq:deterministic_assumption}}{\le} f(x^k) - \frac{4h_k\varepsilon^2}{10} + \frac{h_k\varepsilon^2}{10} = f(x^k) - \frac{3\varepsilon^2}{20L}.
\end{equation*}
Unrolling the recurrence we obtain: 
\begin{eqnarray*}
    f(x^{N+1}) &\le& f(x^0) - \frac{3\varepsilon^2}{20L}(N+1)\\
    &\Big\Downarrow& \\
    N &\le& \frac{20L(f(x^0) - f(x^{N+1}))}{3\varepsilon^2} - 1 \le \frac{20L(f(x^0) - f_*)}{3\varepsilon^2} - 1.
\end{eqnarray*}
Therefore, the methods stops after
\begin{equation}
    N \le \frac{20L(f(x^0) - f_*)}{3\varepsilon^2}\label{eq:deterministic_bound}
\end{equation}
iterations. This bound is optimal up to constant factors \cite{carmon2019lowerI}.

\subsubsection{Stochastic Case: Uniformly Bounded Variance}\label{sec:sgd_ubv}
In this case, we assume that for all $x\in\R^n$ we have
\begin{equation}
    \EE\left[g(x)\mid x\right] = \nabla f(x),\quad \EE\left[\left\|g(x) - \nabla f(x)\right\|_2^2\mid x\right] \le \frac{\varepsilon^2}{2}. \label{eq:stoch_assumption}
\end{equation}
For example, this situation appears when
$$
f(x) = \EE_\xi[f(x,\xi)]
$$
where $\xi$ is a random variable with distribution $\cal{D}$ and $g(x)$ is formed as
\begin{equation}
\label{eq:mini-batch_stoch_grad}
g(x) = \frac{1}{r}\sum\limits_{i=1}^r \nabla f_i(x,\xi_i)
\end{equation}
where $\xi_1,\ldots,\xi_r$ are i.i.d.\ samples from $\cal{D}$ and
\begin{equation}
    \EE_\xi[\nabla f(x,\xi)] = \nabla f(x),\quad \EE_\xi\left[\|\nabla f(x,\xi) - \nabla f(x)\|_2^2\right] \le \sigma^2. \label{eq:bounded_var}
\end{equation}
Indeed, if we choose $r = \max\left\{1,\frac{2\sigma^2}{\varepsilon^2}\right\}$, then due to independence of $\xi_1,\ldots,\xi_r$ we have:
\begin{equation*}
    \EE\left[\|g(x) - \nabla f(x)\|_2^2\mid x\right] = \frac{1}{r^2}\sum\limits_{i=1}^r\EE_{\xi_i}\left[\|\nabla f(x,\xi_i) - \nabla f(x)\|_2^2\right] \le \frac{\sigma^2}{r} \le \frac{\varepsilon^2}{2}.
\end{equation*}
Then, taking conditional expectation $\EE\left[\cdot\mid x^k\right]$ from the both sides of \eqref{eq:main_ineq_non_cvx_opt_methods} we derive
\begin{eqnarray*}
    \EE\left[f(x^{k+1})\mid x^k\right] &\le& f(x^k) - \frac{h_k}{2}\EE\left[\|g^k\|_2^2\mid x^k\right] + h_k\EE\left[\|g^k - \nabla f(x^k)\|_2^2\mid x^k\right]\\
    &=& f(x^k) - \frac{h_k}{2}\|\nabla f(x^k)\|_2^2 - \frac{h_k}{2}\EE\left[\|g^k-\nabla f(x^k)\|_2^2\mid x^k\right]\\
    &&\quad + h_k\EE\left[\|g^k - \nabla f(x^k)\|_2^2\mid x^k\right]\\
    &=& f(x^k) - \frac{h_k}{2}\|\nabla f(x^k)\|_2^2 + \frac{h_k}{2}\EE\left[\|g^k - \nabla f(x^k)\|_2^2\mid x^k\right]\\
    &\overset{\eqref{eq:stoch_assumption}}{\le}& f(x^k) - \frac{h_k}{2}\|\nabla f(x^k)\|_2^2 + \frac{h_k\varepsilon^2}{4}.
\end{eqnarray*}
After that, we take the full expectation from the both sides of the previous inequality, choose $h_k \equiv \frac{1}{2L}$ and sum up the result for $k=0,1,\ldots, N-1$:
\begin{eqnarray*}
    \frac{1}{N}\sum\limits_{k=0}^{N-1}\EE\left[\|\nabla f(x^k)\|_2^2\right] &\le& \frac{4L}{N}\sum\limits_{k=0}^{N-1}\left(\EE[f(x^k)] - \EE[f(x^{k+1})]\right) + \frac{\varepsilon^2}{2}\\
    &=& \frac{4L\left(f(x^0) - \EE[f(x^N)]\right)}{N} + \frac{\varepsilon^2}{2}\\
    &\le& \frac{4L\left(f(x^0) - f_*\right)}{N} + \frac{\varepsilon^2}{2}.
\end{eqnarray*}
Finally, we choose the output of the method $\hat x^N$ uniformly at random from $x^0,x^1,\ldots,x^{N-1}$ which implies
\begin{eqnarray*}
    \EE\left[\|\nabla f(\hat x^N)\|_2^2\right] &\le& \frac{4L\left(f(x^0) - f_*\right)}{N} + \frac{\varepsilon^2}{2}.
\end{eqnarray*}
Taking $N = \frac{8L\left(f(x^0)-f_*\right)}{\varepsilon^2}$ we obtain $\EE\left[\|\nabla f(\hat x^N)\|_2^2\right] \le \varepsilon^2$. Moreover, the total number of stochastic oracle calls (number of $\nabla f(x,\xi)$-calculations) is
\begin{equation*}
    \sum\limits_{k=0}^{N-1}r_k = \max\left\{\frac{8L\left(f(x^0)-f_*\right)}{\varepsilon^2},\frac{16L\left(f(x^0)-f_*\right)\sigma^2}{\varepsilon^4}\right\}.
\end{equation*}
This bound is optimal up to constant factors for the case when the variance is uniformly upper bounded \cite{arjevani2019lower}.

\subsubsection{Stochastic Case: Finite Sum Minimization}\label{sec:finite_sums}
In this case we assume that the objective function has a finite sum structure \eqref{eq:f_finite_sum} with $L$-smooth summands. In fact, this smoothness constant $L$ can be significantly larger than the smoothness constant of $f$. It is essential for providing a fair comparison of different complexity results. It is possible to improve the dependence on $L$ in the final complexity bounds \cite{li2020page} using average smoothness assumption, but for simplicity we consider the case when all summands are $L$-smooth. Moreover, we assume that there exists constant $\sigma^2$ (possibly infinite) such that for $\xi$ taken uniformly at random from $\{1,\ldots,m\}$ and for all $x\in\R^n$
\begin{equation}
    \EE_\xi\left[\|\nabla f_\xi(x) - \nabla f(x)\|_2^2\right] \le \sigma^2. \label{eq:bounded_variance_spider}
\end{equation}
We define $r_k$ and $g^k$ in the following way:
\begin{eqnarray}
    r_k = r &=& \max\left\{1,\frac{20\sigma^2}{\varepsilon^2}\right\}, \label{eq:batch_size_spider}\\
    q &=& \min\left\{r,m\right\}, \label{eq:q_def} \\
    g^k &=& \begin{cases}\frac{1}{r}\sum\limits_{j=1}^r\nabla f_{\xi_{k,j}}(x^k),& \text{if } r < m \text{ and } r \text{ divides } k,\\ \nabla f(x^k),& \text{if } m \le r \text{ and } m \text{ divides } k,\\ \nabla f_{\xi_k}(x^k) - \nabla f_{\xi_k}(x^{k-1}) + g^{k-1},& \text{otherwise}\end{cases} \label{eq:stoch_grad_spider}\\
    h_k = h &=& \frac{1}{10L\sqrt{q}}.\label{eq:stepsize_spider}
\end{eqnarray}
Here, at iteration $k$ random index $\xi_k$ is sampled uniformly at random from $\{1,\ldots,m\}$ if $k$ is not divisible by $q$ and random indices $\xi_{k,1},\ldots,\xi_{k,r}$ are i.i.d.\ samples from uniform distribution on $\{1,\ldots,m\}$ if $q = r$ and $r$ divides $k$. As the result, we obtain the variant of SPIDER \cite{fang2018spider}. We notice that for $k = aq+p$, $p\in\{0,1,\ldots,q-1\}$ iteration $k$ requires $2$ calculations of $\nabla f_\xi(x)$ when $p\neq 0$ and $q$ calculations of $\nabla f_\xi(x)$ when $p = 0$. This implies that $q$ iterations of the method requires only $3q$ calculations of $\nabla f_\xi(x)$, so, if $k \ge q$, then the number of stochastic first-order oracle coincides with the number of iterations up to a constant factor $3$.

Below we present a simplified approach to analyze SPIDER. As before, our goal is to show that $\EE\left[\|g^k - \nabla f(x^k)\|_2^2\right]$ can be upper-bounded by either something small or something that can be controlled by other terms in \eqref{eq:main_ineq_non_cvx_opt_methods}. First of all, if $k = aq$, then
\begin{eqnarray}
    \EE\left[\|g^k - \nabla f(x^k)\|_2^2\right] &=& \begin{cases}0,&\text{if } q = m,\\ \EE\left[\left\|\frac{1}{r}\sum\limits_{j=1}^r\nabla f_{\xi_{k,j}}(x^k) - \nabla f(x^k)\right\|_2^2\right],&\text{if } q=r\end{cases}\notag\\
    &=& \begin{cases}0,&\text{if } q = m,\\ \frac{1}{r^2}\sum\limits_{j=1}^r\EE\left[\left\|\nabla f_{\xi_{k,j}}(x^k) - \nabla f(x^k)\right\|_2^2\right],&\text{if } q=r\end{cases}\notag\\
    &\overset{\eqref{eq:bounded_variance_spider}}{\le}& \begin{cases}0,&\text{if } q = m,\\ \frac{\sigma^2}{r},&\text{if } q=r,\end{cases} \overset{\eqref{eq:batch_size_spider}}{\le} \begin{cases}0,&\text{if } q = m,\\ \frac{\varepsilon^2}{20},&\text{if } q=r,\end{cases} \label{eq:var_bound_k_aq}
\end{eqnarray}
where $\EE\left[\left\|\frac{1}{r}\sum\limits_{j=1}^r\nabla f_{\xi_{k,j}}(x^k) - \nabla f(x^k)\right\|_2^2\right] = \frac{1}{r^2}\sum\limits_{j=1}^r\EE\left[\left\|\nabla f_{\xi_{k,j}}(x^k) - \nabla f(x^k)\right\|_2^2\right]$ due to independence of $\xi_{k,1},\ldots,\xi_{k,r}$ and in the third inequality we applied the tower property: $\EE[\cdot] = \EE\left[\EE\left[\cdot\mid x^k\right]\right]$. Secondly, if $k = aq + p$ with $p\in\{1,\ldots,q-1\}$ we have
\begin{eqnarray}
    \EE\left[\left\|g^k - \nabla f(x^k)\right\|_2^2\right] &\overset{\eqref{eq:stoch_grad_spider}}{=}& \EE\left[\left\|\nabla f_{\xi_k}(x^k) - \nabla f_{\xi_k}(x^{k-1}) + g^{k-1} - \nabla f(x^k)\right\|_2^2\right] \notag\\
    &=& \EE\left[\left\|\nabla f_{\xi_k}(x^k) - \nabla f_{\xi_k}(x^{k-1}) - \nabla f(x^k) + \nabla f(x^{k-1})\right\|_2^2\right]\notag\\
    &&\quad + \EE\left[\left\|g^{k-1} - \nabla f(x^{k-1})\right\|_2^2\right]\notag
\end{eqnarray}
where we use the variance decomposition\footnote{Here $\EE_{\xi_k}[\cdot]$ is a mathematical expectation conditioned on everything despite $\xi_k$, i.e.\ expectation is taken w.r.t.\ the randomness coming only from $\xi_k$.} $\EE_{\xi_k}\left[\|\eta\|_2^2\right] = \EE_{\xi_k}\left[\|\eta-\EE_{\xi_k}[\eta]\|_2^2\right] + \left\|\EE_{\xi_k}\left[\eta\right]\right\|$ for random vector $\eta = \nabla f_{\xi_k}(x^k) - \nabla f_{\xi_k}(x^{k-1}) + g^{k-1} - \nabla f(x^k)$ together with the tower property $\EE[\cdot] = \EE\left[\EE_{\xi_k}\left[\cdot\right]\right]$. Using the inequality above together with $\|a+b\|_2^2 \le 2\|a\|_2^2 + 2\|b\|_2^2$, $a,b\in\R^n$ and $L$-smoothness of $f_1,\ldots,f_m, f$ we get
\begin{eqnarray}
    \EE\left[\left\|g^k - \nabla f(x^k)\right\|_2^2\right] &\le& 2\EE\left[\left\|\nabla f_{\xi_k}(x^k) - \nabla f_{\xi_k}(x^{k-1})\right\|_2^2\right] + 2\EE\left[\left\|\nabla f(x^k) - \nabla f(x^{k-1})\right\|_2^2\right]\notag\\
    &&\quad + \EE\left[\left\|g^{k-1} - \nabla f(x^{k-1})\right\|_2^2\right]\notag\\
    &\le& 4L^2\EE\left[\|x^k - x^{k-1}\|_2^2\right] + \EE\left[\left\|g^{k-1} - \nabla f(x^{k-1})\right\|_2^2\right]\notag\\
    &=& 4L^2h^2\EE\left[\|g^{k-1}\|_2^2\right]+ \EE\left[\left\|g^{k-1} - \nabla f(x^{k-1})\right\|_2^2\right]\notag\\
    &\le& 8L^2h^2\EE\left[\|\nabla f(x^{k-1})\|_2^2\right] + \left(1+8L^2h^2\right)\EE\left[\left\|g^{k-1} - \nabla f(x^{k-1})\right\|_2^2\right].\notag
\end{eqnarray}
Unrolling the recurrence we derive
\begin{eqnarray}
    \EE\left[\left\|g^k - \nabla f(x^k)\right\|_2^2\right] &\le& 8L^2h^2\sum\limits_{l=1}^{p}\left(1+8L^2h^2\right)^{l-1}\EE\left[\|\nabla f(x^{k-l})\|_2^2\right] \notag\\
    &&\quad+ \left(1+8L^2h^2\right)^p\EE\left[\left\|g^{aq} - \nabla f(x^{aq})\right\|_2^2\right]\notag\\
    &\overset{\eqref{eq:var_bound_k_aq},p\le q}{\le}&\left(1+8L^2h^2\right)^q\sum\limits_{l=1}^p 8L^2h^2\EE\left[\|\nabla f(x^{aq+l})\|_2^2\right]\notag\\
    &&\quad+ \left(1+8L^2h^2\right)^q\begin{cases}0,&\text{if } q = m,\\ \frac{\varepsilon^2}{20},&\text{if } q=r\end{cases}\notag\\
    &\overset{(1+x)^q \le e^{qx}}{\le}& \exp\left(8L^2h^2q\right)\sum\limits_{l=1}^p 8L^2h^2\EE\left[\|\nabla f(x^{aq+l})\|_2^2\right]\notag\\
    &&\quad+ \exp\left(8L^2h^2q\right)\begin{cases}0,&\text{if } q = m,\\ \frac{\varepsilon^2}{20},&\text{if } q=r.\end{cases}\notag
\end{eqnarray}
Next, using the choice of the stepsize $h = \nicefrac{1}{(10L\sqrt{q})}$ we obtain
\begin{eqnarray}
    \EE\left[\left\|g^k - \nabla f(x^k)\right\|_2^2\right] &\le& \sum\limits_{l=1}^p 9L^2h^2\EE\left[\|\nabla f(x^{aq+l})\|_2^2\right] + \frac{11\varepsilon^2}{200}.\label{eq:spider_var_bound}
\end{eqnarray}
Finally, we put all the inequalities together. We start with modifying \eqref{eq:main_ineq_non_cvx_opt_methods}:
\begin{eqnarray*}
    f(x^{k+1}) &\le& f(x^k) - \frac{h_k}{2}\|g^k\|_2^2 + h_k\|\nabla f(x^k) - g^k\|_2^2\\
    &\le& f(x^k) - \frac{h}{4}\|\nabla f(x^k)\|_2^2 + \frac{3h}{2}\|\nabla f(x^k) - g^k\|_2^2,
\end{eqnarray*}
where we used that inequality $\|a+b\|_2^2 \ge \frac{1}{2}\|a\|_2^2 - \|b\|_2^2$ holds for all $a,b\in\R^n$ (in particular, we use $a = \nabla f(x^k)$ and $b = g^k - \nabla f(x^k)$). Next, we take the full mathematical expectation from the both sides of previous inequality (taking into account that $k = aq+p$):
\begin{eqnarray*}
    \EE[f(x^{aq+p+1})] &\le& \EE[f(x^{aq+p})] - \frac{h}{4}\EE\left[\|\nabla f(x^{aq+p})\|_2^2\right] + \frac{3h}{2}\EE\left[\|g^{aq+p} - \nabla f(x^{aq+p})\|_2^2\right]\\
    &\le& \EE[f(x^{aq+p})] - \frac{h}{4}\EE\left[\|\nabla f(x^{aq+p})\|_2^2\right] + \frac{3h}{2}\sum\limits_{l=1}^p 9L^2h^2\EE\left[\|\nabla f(x^{aq+l})\|_2^2\right]\\
    &&\quad + \frac{33h\varepsilon^2}{400}.
\end{eqnarray*}
We notice that this inequality holds for all integers $a\ge 0$ and $p \in\{0,\ldots,q-1\}$. Summing up these inequalities for $p = 0,\ldots,P$ and taking $a = A$ where $N = Aq + P$, $P\in\{0,\ldots,q-1\}$  we get
\begin{eqnarray*}
    0 &\le& \sum\limits_{p=0}^{P}\left(\EE[f(x^{Aq+p})] - \EE[f(x^{Aq+p+1})]\right) - \frac{h}{4}\sum\limits_{p=0}^{P}\EE\left[\|\nabla f(x^{Aq+p})\|_2^2\right] \\
    &&\quad+ \frac{27L^2h^3}{2}\sum\limits_{p=0}^P\sum\limits_{l=1}^p\EE\left[\|\nabla f(x^{Aq+l})\|_2^2\right] + \frac{33h\varepsilon^2 (P+1)}{400}\\
    &\overset{P\le q-1}{\le}& \EE\left[f(x^{Aq})\right] - \EE\left[f(x^{Aq+P+1})\right] - h\left(\frac{1}{4} - \frac{27L^2h^2q}{2}\right)\sum\limits_{p=0}^{P}\EE\left[\|\nabla f(x^{Aq+p})\|_2^2\right]\\
    &&\quad +\frac{33h\varepsilon^2 (P+1)}{400} \\
    &\overset{\eqref{eq:stepsize_spider}}{=}& \EE\left[f(x^{Aq})\right] - \EE\left[f(x^{Aq+P+1})\right] - \frac{23h}{200}\sum\limits_{p=0}^{P}\EE\left[\|\nabla f(x^{Aq+p})\|_2^2\right] +\frac{33h\varepsilon^2 (P+1)}{400},
\end{eqnarray*}
hence 
\begin{eqnarray*}
    \frac{23h}{200}\sum\limits_{p=0}^{P}\EE\left[\|\nabla f(x^{Aq+p})\|_2^2\right] &\le& \EE\left[f(x^{Aq})\right] - \EE\left[f(x^{Aq+P+1})\right] +\frac{33h\varepsilon^2 (P+1)}{400}.
\end{eqnarray*}
These inequalities hold for all $A$ and $P$. Then we can sum up these inequalities for $(A,P) = (0,q-1),(1,q-1),\ldots, (\hat{A},\hat{P})$ and get that for $\hat N = \hat{A}q + \hat{P}$ and divide the result by $\frac{23h(\hat{N}+1)}{200}$ and get
\begin{eqnarray*}
    \frac{1}{\hat N+1}\sum\limits_{k=0}^{\hat{N}}\EE\left[\|\nabla f(x^k)\|_2^2\right] &\le& \frac{200\left(f(x^0) - \EE\left[f(x^{\hat N + 1})\right]\right)}{23h(\hat{N}+1)} + \frac{33\varepsilon^2}{46}\\
    &\overset{\eqref{eq:stepsize_spider}}{\le}& \frac{2000L\sqrt{q}\left(f(x^0) - f_*\right)}{23(\hat{N}+1)} + \frac{33\varepsilon^2}{46}.
\end{eqnarray*}
Finally, taking $\hat x^{\hat{N}}$ uniformly at random from $x^0,\ldots,x^{\hat{N}}$ we get
\begin{equation*}
    \EE\left[\|\nabla f(\hat x^{\hat N})\|_2^2\right] \le \frac{2000L\sqrt{q}\left(f(x^0) - f_*\right)}{23(\hat{N}+1)} + \frac{33\varepsilon^2}{46}.
\end{equation*}
This implies that after
\begin{eqnarray*}
    \hat{N} &=& \frac{4000L\sqrt{q}(f(x^0) - f(x^\ast))}{13\varepsilon^2}\\ &\overset{\eqref{eq:batch_size_spider},\eqref{eq:q_def}}{=}& \frac{4000L(f(x^0) - f(x^\ast))}{13\varepsilon^2}\min\left\{\sqrt{m},\max\left\{1, \frac{\sqrt{20}\sigma}{\varepsilon}\right\}\right\}
\end{eqnarray*}
iterations we reach $\EE\left[\|\nabla f(\hat x^{\hat N})\|_2^2\right] \le \varepsilon^2$. Moreover, it requires 
\begin{equation*}
    O\left(\frac{L(f(x^0) - f(x^\ast))}{\varepsilon^2}\min\left\{\sqrt{m},\max\left\{1, \frac{\sigma}{\varepsilon}\right\}\right\} + \min\left\{m,\max\left\{1, \frac{\sigma^2}{\varepsilon^2}\right\}\right\}\right)
\end{equation*}
calculations of $\nabla f_\xi(x)$ which is optimal up to constant factors \cite{fang2018spider}.

\subsection{SGD and Its Variants}\label{sec:sgd_like}
As it was shown in the previous section, SGD 
\begin{equation}
    x^{k+1} = x^k - h_kg(x^k),\quad \EE[g(x)] = \nabla f(x) \label{eq:SGD_dynamics}
\end{equation}
in the settings of Section~\ref{sec:sgd_ubv} requires $O\left(\frac{L(f(x^0) - f_*}{\epsilon^2}\right)$ iterations with \kesha{batch size} $r = \Theta\left(\max\left\{1,\frac{\sigma^2}{\epsilon^2}\right\}\right)$ to find an $\epsilon$-stationary point in expectation. The total number of stochastic first-order oracle calls equals 
\begin{equation}
    O\left(\frac{L(f(x^0)-f_*)}{\epsilon^2}\max\left\{1,\frac{\sigma^2}{\epsilon^2}\right\}\right). \label{eq:SGD_non_cvx_bound}
\end{equation}
We emphasize that we use large \kesha{batch size} for the sake of simplicity and unification of the results in 3 different cases. In fact, it is possible to obtain the bound \eqref{eq:SGD_non_cvx_bound} using smaller stepsizes and constant batch sizes of the order $O(1)$ \cite{ghadimi2013stochastic}.

\subsubsection{Assumptions on the Stochastic Gradient}\label{sec:SGD_assumptions}
In addition to assumption \eqref{eq:bounded_var}, which is quite restrictive, there exist several other assumptions on the stochastic gradient studied in the literature. Recently 
% Khaled and Richt\'{a}rik 
in \cite{khaled2020better} it was proposed a simple and unified way to cover the most popular ones.
\begin{assumption}[Expected Smoothness; Assumption 2 from \cite{khaled2020better}]\label{as:expected_smoothness}
    The second moment of stochastic gradients satisfies
    \begin{equation}
        \EE\left[\|g(x)\|_2^2\right] \le 2A\left(f(x) - f_*\right) + B\|\nabla f(x)\|_2^2 + C\label{eq:expected_smoothness}
    \end{equation}
    for some $A, B, C \ge 0$ and for all $x\in\R^n$.
\end{assumption}
This assumption generalizes the notion of expected smoothness introduced and adjusted for convex problems in \cite{gower2019sgd}. Moreover, the following assumptions are stronger than Assumption~\ref{as:expected_smoothness} or can be seen as special cases of Assumption~\ref{as:expected_smoothness} (see more details and formal proofs in \cite{khaled2020better}).

\noindent \textbf{Uniformly upper-bounded variance (UV) assumption.} Indeed, if $A = 0$, $B=1$ and $C = \sigma^2$, then using variance decomposition inequality \eqref{eq:expected_smoothness} implies \eqref{eq:bounded_var}:
\begin{equation*}
    \EE\left[\|g(x) - \nabla f(x)\|_2^2\right] = \EE\left[\|g(x)\|_2^2\right] - \|\nabla f(x)\|_2^2 \overset{\eqref{eq:expected_smoothness}}{\le} \sigma^2.
\end{equation*}

\noindent \textbf{Expected strong growth condition (E-SG).} When $A = C = 0$ and $B = \alpha \ge 1$ inequality \eqref{eq:expected_smoothness} transforms into so-called expected strong growth condition \cite{solodov1998incremental,vaswani2018fast}:
\begin{equation}
    \EE\left[\|g(x)\|_2^2\right] \le \alpha\|\nabla f(x)\|_2^2.\label{eq:expected_strong_growth}
\end{equation}

\noindent \textbf{Maximal strong growth condition (M-SG)} \cite{tseng1998incremental,schmidt2013fast} states that there exists such $\alpha > 0$ that
\begin{equation}
    \|g(x)\|_2^2 \le \alpha \|\nabla f(x)\|_2^2 \text{ almost surely for all } x\in\R^n. \label{eq:maximal_strong_growth}
\end{equation}
This condition implies E-SG \eqref{eq:expected_strong_growth} while known convergence results in expectation under M-SG assumption have no advantage in comparison with their counterparts under E-SG.

\noindent \textbf{Relaxed growth condition (RG)} \cite{bottou2018optimization} can be seen as another special case of Assumption~\ref{as:expected_smoothness} with $A = 0$, $B = \alpha \ge 1$ and $C = \beta\ge 0$ or as an extension of E-SG:
\begin{equation}
    \EE\left[\|g(x)\|_2^2\right] \le \alpha\|\nabla f(x)\|_2^2 + \beta.\label{eq:relaxed_growth}
\end{equation}
However, there exist simple problems of type \eqref{eq:main_pr_non_cvx}+\eqref{eq:f_expectation} that fit the settings we are interested in but do not satisfy \eqref{eq:relaxed_growth} (see Proposition~1 from \cite{khaled2020better}).

\noindent \textbf{Gradient confusion condition (GC)} \cite{sankararaman2019impact} was developed for the finite-sum case \eqref{eq:f_finite_sum}. In particular, it states that there exists such $\eta > 0$ that for all $i,j=1,\ldots,m$ and for all $x\in\R^n$
\begin{equation}
    \la\nabla f_i(x), \nabla f_j(x)\ra \ge -\eta. \label{eq:gradient_confusion}
\end{equation}
One can show (see Theorem~1, \cite{khaled2020better}) that inequality \eqref{eq:gradient_confusion} implies \eqref{eq:relaxed_growth} with $\alpha = m$ and $\beta = \eta (m-1)$, and, as a consequence, it is a special case of Assumption~\ref{as:expected_smoothness} with $A = 0$, $B = m$, and $C = \eta (m-1)$.

\noindent \textbf{Sure-smoothness condition (SS)} \cite{lei2019stochastic} is defined for the case when the objective is represented as an expectation \eqref{eq:f_expectation} and $g(x) = \nabla f(x,\xi)$ where $\xi$ is sampled independently at each iteration of SGD. That is, sure-smoothness condition means that\footnote{In the original paper \cite{lei2019stochastic}, authors considered more general situation when stochastic realizations $f(x,\xi)$ have H\"{o}lder-continuous gradients.} for all $x,y\in\R^n$
\begin{equation}
    \|\nabla f(x,\xi) - \nabla f(y,\xi)\|_2 \le L\|x-y\|_2 \;\;\;\text{and}\;\;\; f(x,\xi) \ge 0\;\;\; \text{almost surely in } \xi. \label{eq:sure_smoothness}
\end{equation}
Applying classical corollaries of $L$-smoothness one can derive inequality \eqref{eq:expected_smoothness} with $A = 2L$, $B=0$, and $C = 2Lf_*$ from \eqref{eq:sure_smoothness}.

Next, Assumption~\ref{as:expected_smoothness} covers \textbf{arbitrary sampling} setup and distributed setup with quantization\footnote{This technique is applied in distributed optimization to reduce the overall communication cost (e.g., see \cite{alistarh2017qsgd,beznosikov2020biased,pmlr-v139-gorbunov21a}). However, methods for distributed optimization are out of scope of our survey.}. For simplicity, we mention only \textbf{sampling with replacement} as a special case of arbitrary sampling (see more examples in \cite{khaled2020better}). In particular, consider the finite-sum optimization problem \eqref{eq:main_pr_non_cvx}+\eqref{eq:f_finite_sum} and assume that $f_i$ is $L_i$-smooth and bounded from below by $f_{i,*}$ for all $i=1,\ldots,m$. Moreover, assume that $g(x) = \nabla f_j(x)$ where $j = i$ with probability $p_i\ge 0$, $i=1,\ldots,m$, $\sum_{i=1}^m p_i =1$. Then, one can prove \cite{khaled2020better} that Assumption~\ref{as:expected_smoothness} is satisfied in this case with $A = \max_i\frac{L_i}{mp_i}$, $B = 0$, and $C = 2A\Delta_* = \frac{2A}{m}\sum_{i=1}^m\left(f_* - f_{i,*}\right)$. That is, if we apply \textbf{uniform sampling}, i.e., $p_i = \frac{1}{m}$ for all $i=1,\ldots,m$, then we get $A = \max_i L_i$, $B=0$, $C = 2\max_i L_i\Delta_*$, and if \textbf{importance sampling} with $p_i = \frac{L_i}{\sum_{l=1}^mL_l}$ is applied, then Assumption~\ref{as:expected_smoothness} holds with $A = \overline{L} = \frac{1}{m}\sum_{i=1}^mL_i$, $B=0$, and $C = 2\overline{L}\Delta_*$.

Finally, under Assumption~\ref{as:expected_smoothness} Khaled and Richt\'{a}rik \cite{khaled2020better} derived the following complexity bound: if $h = \min\left\{\frac{1}{\sqrt{LAN}},\frac{1}{LB},\frac{\varepsilon}{2LC}\right\}$, then inequality
\begin{equation}
    \min\limits_{0\le k \le N-1}\EE\left[\|\nabla f(x^k)\|_2\right] \le \varepsilon \label{eq:min_norm_bound}
\end{equation}
is satisfied after
\begin{equation}
    N = O\left(\frac{L(f(x^0)-f_*)}{\varepsilon^2}\max\left\{B, \frac{A(f(x^0)-f(x^*))}{\varepsilon^2}, \frac{C}{\varepsilon^2}\right\}\right) \label{eq:khaled_bound}
\end{equation}
iterations of SGD. It is worth to mention that this bound gives the sharpest rates for all known special cases. We summarize some of them in Table~\ref{tab:sgd_bounds}. We notice that \eqref{eq:min_norm_bound} is weaker than \eqref{eq:epsilon_stationary_point}, but it is easy to obtain the same bound \eqref{eq:khaled_bound} guaranteeing \eqref{eq:epsilon_stationary_point} instead of \eqref{eq:min_norm_bound} based on the analysis given in \cite{khaled2020better}.
\begin{table}[h!]
    \normalsize
    \centering
    \begin{tabular}{|c|c|c|c|}
        \hline
        Problem & Settings & Citation & Complexity \\
        \hline
        \eqref{eq:main_pr_non_cvx}+\eqref{eq:f_expectation} & \textbf{UV} \eqref{eq:bounded_var} & \cite{ghadimi2013stochastic} & $\frac{L\Delta_0}{\varepsilon^2}\max\left\{1,\frac{\sigma^2}{\varepsilon^2}\right\}$\\
         \hline
         \eqref{eq:main_pr_non_cvx}+\eqref{eq:f_expectation}/\eqref{eq:f_finite_sum} & \textbf{RG} \eqref{eq:relaxed_growth} & \cite{bottou2018optimization,vaswani2018fast} & $\frac{L\Delta_0}{\varepsilon^2}\max\left\{\alpha,\frac{\beta}{\varepsilon^2}\right\}$\\
         \hline
         \eqref{eq:main_pr_non_cvx}+\eqref{eq:f_finite_sum} & \textbf{GC} \eqref{eq:gradient_confusion} & \cite{sankararaman2019impact} & $\frac{L\Delta_0}{\varepsilon^2}\max\left\{m,\frac{\eta(m-1)}{\varepsilon^2}\right\}$\\
         \hline
         \eqref{eq:main_pr_non_cvx}+\eqref{eq:f_finite_sum} & \textbf{Uniform Sampling} & \cite{khaled2020better} & $\frac{L\max_i L_i\Delta_0}{\varepsilon^4}\max\left\{\Delta_0,\Delta_*\right\}$\\
         \hline
         \eqref{eq:main_pr_non_cvx}+\eqref{eq:f_finite_sum} & \textbf{Importance Sampling} & \cite{khaled2020better} & $\frac{L\overline{L}\Delta_0}{\varepsilon^4}\max\left\{\Delta_0,\Delta_*\right\}$\\
         \hline
    \end{tabular}
    \caption{Summary of the complexity results for SGD under different assumptions on the stochastic gradient. The column ``Complexity'' contains an overall number of stochastic first-order oracle calls needed to find $\varepsilon$-stationary point neglecting constant factors. Notation: $\Delta_0 = f(x^0) - \eduard{f_*}$, $\sigma^2 =$ a uniform bound for the variance of the stochastic gradient \eqref{eq:bounded_var}, $\alpha,\beta$ = relaxed growth condition parameters, $\eta$ = gradient confusion parameter, $\Delta_* = \frac{1}{m}\sum_{i=1}^m(f_* - f_{i,*})$, $\max_i L_i$ = maximal smoothness constant of $f_i$ in \eqref{eq:f_finite_sum}, $\overline{L}$ = averaged smoothness constant of $f_i$ in \eqref{eq:f_finite_sum}.}
    \label{tab:sgd_bounds}
\end{table}

\subsubsection{The Choice of the Stepsize}
In practice, instead of using the constant stepsize for SGD it is popular to periodically decrease the stepsize by some factor \cite{bottou2010large,krizhevsky2009learning,he2016deep} even for non-convex problems. For strongly convex problems such a choice is natural: it is well-known \cite{gorbunov2020unified} that if the stepsize equals $h$ and strong convexity parameter equals $\mu$, then SGD converges with linear rate $\widetilde{O}((h\mu)^{-1})$ to the neighborhood of the solution with size proportional to $h$. Surprisingly, SGD enjoys similar behaviour even for non-convex problems which was recently shown in \cite{shi2020learning}. 

In the neural networks training, ``warmup'' \cite{goyal2017accurate, gotmare2018closer} and cyclical stepsize \cite{smith2017cyclical,loshchilov2016sgdr} schedules are also very popular and useful. The first one refers to the strategy when, during several epochs of training, tiny stepsizes are used, and then they are increased. This technique was successfully applied for several deep learning problems like ResNet \cite{he2016deep}, large-batch training of Imagenet \cite{goyal2017accurate} and natural language problems \cite{vaswani2017attention, devlin2018bert}. 

Cyclical stepsize schedule means that the stepsize is changing between some lower and upper bounds. There are different modification of this technique including gradual decrease and increase during one epoch \cite{smith2017cyclical} and gradual decrease of the stepsize followed by the sudden increase \cite{loshchilov2016sgdr}. However, the theoretical understanding of the success of ``warmup'' and cyclical schedules is very limited.

We also discuss different stepsize policies including adaptive ones (Section~\ref{sec:adaptive_methods}), Armijo line-search under expected strong growth assumption and stochastic Polyak stepsizes under relaxed growth assumption (Section~\ref{sec:over_param}) in the following subsections.

\subsubsection{Over-Parameterized Models}\label{sec:over_param}
In Section~\ref{ProvableGlobalResults}, we mentioned that over-parameterization \cite{livni2014computational,neyshabur2017exploring,zhang2016understanding,nguyen2018loss,li2018over,allen2019convergence, allen2019on}, meaning that the last layer has more neurons than the number of samples in the training set, is a good property for neural networks from the optimization and generalization \cite{ma2018power, allen2019learning, allen2019can} point perspectives, but not a panacea: over-parameterized neural networks have no spurious valleys, but still can have bad local minima \cite{ding2019spurious}.

In the papers, focusing mostly on the optimization aspects of over-parameterized models, it was shown that SGD converges with the same (up to the difference in the smoothness constants) rate as GD in terms of the iteration complexity in convex and strongly convex cases \cite{vaswani2018fast,vaswani2019painless,loizou2020stochastic} under interpolation condition: for the finite-sum optimization problem \eqref{eq:main_pr_non_cvx}+\eqref{eq:f_finite_sum} there exists such point $x^*\in\R^n$ that
\begin{equation}
    \min\limits_{x\in\R^n}f_i(x) = f_i(x^*)\quad \forall i = 1,\ldots, m. \label{eq:interpolation_cond}
\end{equation}
Furthermore, in this setting SGD converges with Armijo line-search \cite{vaswani2019painless}, with stochastic Polyak stepsizes \cite{loizou2020stochastic}, and, if additionally expected strong growth condition \eqref{eq:expected_strong_growth} holds, SGD can be accelerated \cite{vaswani2018fast} and the accelerated version converges as good as Nesterov's method \cite{nesterov1983method} in terms of iteration complexity up to expected strong growth multiplicative factor $\alpha$ from \eqref{eq:expected_strong_growth}.

In the general non-convex case, the following results exist.\newline
\textbf{Constant stepsizes.} In \cite{vaswani2018fast}, it was shown that SGD with constant stepsize $h = \nicefrac{1}{\alpha L}$ finds $\varepsilon$-stationary point under expected strong growth condition \eqref{eq:expected_strong_growth} with the rate $O\left(\nicefrac{\alpha L (f(x^0)-f_*)}{\varepsilon^2}\right)$ matching the iteration complexity of GD up to the factor $\alpha$.\newline
\textbf{Armijo line-search.} The idea that under interpolation condition/expected strong growth condition SGD and GD have similar properties was then strengthen in \cite{vaswani2019painless}, where authors showed that SGD with Armijo line-search converges in these settings. In particular, the authors of \cite{vaswani2019painless} considered such stepsizes $h_k$ that
\begin{equation}
    f_{i_k}(x^k - h_k\nabla f_{i_k}(x^k)) \le f_{i_k}(x^k) - ch_k\|\nabla f_{i_k}(x^k)\|_2^2, \label{eq:armijo_line_search}
\end{equation}
where the index $i_k$ is sampled uniformly at random from the set $\{1,\ldots, m\}$, the stochastic gradient $g^k$ is defined as $g^k = \nabla f_{i_k}(x^k)$, and $c > 0$ is a hyper-parameter. Moreover, it is assumed that $h_k \in (0,h_{\max}]$ for all $k\ge 0$. Then SGD with Armijo line-search \eqref{eq:armijo_line_search} with $c > 1 - \nicefrac{L_{\max}}{(\alpha L)}$ and $h_{\max} \le \nicefrac{2}{\alpha L}$ finds $\varepsilon$-stationary point under expected strong growth condition \eqref{eq:expected_strong_growth} with the rate $O\left(\nicefrac{(f(x^0)-f_*)}{(\delta\varepsilon^2)}\right)$, where $\delta = \left(h_{\max} + \nicefrac{2(1-c)}{L_{\max}}\right) - \alpha\left(h_{\max} - \frac{2(1-c)}{L_{\max}} + Lh_{\max}^2\right)$, $L_{\max}$ is the maximal smoothness constant of summands $f_i$, and $f$ is the smoothness constant of $f$. Authors of \cite{vaswani2019painless} also considered the version with samples used for backtracking \eqref{eq:armijo_line_search} independent from those used for determining the stochastic gradient, and the version with non-increasing stepsizes under additional assumption that the iterates lie in some ball with radius $D$. The rates are $O\left(\nicefrac{\max\{L_{\max},\alpha L\}(f(x^0)-f_*)}{\varepsilon^2}\right)$ and $O\left(\nicefrac{\max\{L_{\max},\alpha L\}LD^2}{\varepsilon^2}\right)$ respectively, and both complexity bounds hold with $c = \nicefrac{1}{2}$ and $h_{\max} = \nicefrac{1}{(\alpha L)}$. Finally, in the numerical experiments from \cite{vaswani2019painless} the authors observed that the method's performance is robust to the choic of $c$ and $h_{\max}$.\newline
\textbf{Stochastic Polyak stepsizes.} Next, SGD under expected strong growth condition converges with stochastic Polyak stepsizes introduced and analyzed in \cite{loizou2020stochastic}:
\begin{equation}
    h_k = \min\left\{\frac{f_{i_k}(x^k) - f_{i_k,*}}{c\|\nabla f_{i_k}(x^k)\|_2},h_b\right\},\label{eq:stochastic_polyak_stepsize}
\end{equation}
where the index $i_k$ is sampled uniformly at random from the set $\{1,\ldots, m\}$, the stochastic gradient $g^k$ is defined as $g^k = \nabla f_{i_k}(x^k)$, $f_{i,*}$ is uniform lower bound for $f_i(x)$, and $c > 0$ is a hyper-parameter. In particular, one can show \cite{loizou2020stochastic} that SGD in these settings with $c > \nicefrac{\alpha L}{4L_{\max}}$ and $h_b \le \max\left\{\nicefrac{2}{(\alpha L)},\overline{h}_b\right\}$ finds $\varepsilon$-stationary point under expected strong growth condition \eqref{eq:expected_strong_growth} with the rate $O\left(\nicefrac{(f(x^0)-f_*)}{(\delta\varepsilon^2)}\right)$, where $\delta = \left(h_{b} + \beta\right) - \alpha\left(h_{b} - \beta + Lh_{b}^2\right)$, $\beta = \min\left\{\nicefrac{1}{(2cL_{\max})}, h_b\right\}$, and
\begin{equation*}
    \overline{h}_b = \frac{-(\alpha-1) + \sqrt{(\alpha-1)^2 + \frac{4L\alpha(\alpha+1)}{2cL_{\max}}}}{2L\alpha}.
\end{equation*}

\subsubsection{Proximal Variants}
\eduard{In the previous subsections, all complexity results rely on the smoothness of the objective function. The natural question arises: is it possible to generalize these results to the non-smooth case? In the recent work \cite{kornowski2021oracle}, the authors give a negative answer to this question for \textit{generally non-smooth} non-convex functions, i.e., one cannot find efficiently via first-order methods near $\varepsilon$-stationary points. However, m}any complexity results that we mentioned before and will mention in the following subsections have generalizations to the composite optimization problems:
\begin{equation}
    \min\limits_{x\in\R^n}\left\{F(x) = f(x) + R(x)\right\}, \label{eq:composite_problem}
\end{equation}
where the function $f$ is $L$-smooth, but, possibly, non-convex, while $R(x)$, i.e., composite term/regularizer, is a proper closed convex function which can be non-smooth. Moreover, function $R(x)$ is often chosen in such a way that the proximal operator
\begin{equation}
    \text{prox}_{R}(x) = \argmin\limits_{y\in\R^n}\left\{R(y) + \frac{1}{2}\|y-x\|_2^2\right\} \label{eq:prox_operator}
\end{equation}
can be easily computed, and to make the solution of the problem satisfy certain properties, e.g., sparsity; see \cite{candes2008enhancing, combettes2011proximal, bach2012optimization} for the detailed discussion and examples of regularizers.

In these settings, instead of SGD one can apply prox-SGD defined by the following recurrence:
\begin{equation}
    x^{k+1} = \text{prox}_{h_k R}(x^k - h_kg^k). \label{eq:prox_SGD}
\end{equation}
Moreover, to measure the progress of the method the generalized projected stochastic gradient is used: $\widetilde{g}^k = \nicefrac{(x^{k}-x^{k+1})}{h_k}$. When the regularizer $R(x)$ is a constant $\widetilde{g}^k = g^k$. For proximal stochastic methods we say that the iterate $x^k$ is $\varepsilon$-stationary point if 
\begin{equation}
    \EE\left[\|\widetilde{g}^k\|_2^2\right] \le \varepsilon^2. \label{eq:generalized_stationary_point}
\end{equation}

In \cite{ghadimi2016mini-batch}, it was shown that prox-SGD under uniformly upper-bounded variance assumption \eqref{eq:bounded_var} converges with the rate given in \eqref{eq:SGD_non_cvx_bound}. However, the analysis from \cite{ghadimi2016mini-batch} works only in the large-batch setting, i.e., when \kesha{batch sizes} are of the order $O(\varepsilon^{-2})$. For a long time, there was no analysis establishing the same bound without using $O(\varepsilon^{-2})$ batches, and the problem was recently resolved in \cite{davis2019stochastic}.

\subsubsection{Momentum-SGD}
As we already mentioned, SGD is optimal among stochastic first-order methods for finding $\varepsilon$-stationary points under uniformly bounded variance assumption \cite{arjevani2019lower}. However, it does not imply that there is no sense in using different methods for such problems. In practice, different additional tricks are applied to improve the convergence of SGD, and, perhaps, the most popular one is momentum \cite{polyak1964some}.

Momentum-SGD/Heavy Ball SGD can be written in different forms. Usually it is written as
\begin{eqnarray}
    m^{k+1} &=& \beta_k m^k + g(x^k),\notag\\
    x^{k+1} &=& x^k - h_k g(x^k), \notag
\end{eqnarray}
where parameter $\beta_k \in [0,1)$ is called momentum parameter. In the convex and strongly convex cases this method has some advantages in comparison to SGD like better last-iterate convergence guarantees \cite{tao2018primal, taylor2019stochastic, sebbouh2020convergence}, but does not have an accelerated rate \cite{kidambi2018insufficiency}. In the non-convex case, Momentum-SGD has the same complexity guarantee \eqref{eq:SGD_non_cvx_bound} as SGD under uniformly bounded variance assumption \cite{yan2018unified, defazio2020understanding}. However, in practice, Momentum-SGD often works much better than SGD especially on computer vision problems \cite{sutskever2013importance}, and also navigates ravines and escapes saddle points better than SGD.

Among other works on Momentum-SGD we emphasize the recent paper \cite{defazio2020understanding} establishing the tight convergence rates for Momentum-SGD in Stochastic Primal Averaging \cite{tao2018primal} form via Lyapunov functions analysis. In particular, \cite{defazio2020understanding} justifies (theoretically and/or empirically) the following important insights about the behavior of Momentum-SGD: (i) Momentum-SGD is provably better than SGD during the early stage of the convergence, (ii) it is better to gradually reduce momentum parameter $\beta_k$ rather than the stepsize $h_k$, and (iii) gradual changes of the parameters of Momentum-SGD are preferable than sudden changes.

\subsubsection{Random Reshuffling}
Before this subsection, we always assumed that stochastic gradients are sampled independently from previous iterations. However, in the context of finite sum optimization \eqref{eq:main_pr_non_cvx}+\eqref{eq:f_finite_sum}, the different sampling strategy called Random Reshuffling (or SGD with Without Replacement sampling) is often used: at each epoch (pass through the dataset) random permutation $\{i_1,i_2,\ldots,i_m\}$ of the set $\{1,2,\ldots,m\}$ is generated defining the order of gradients computations (see Algorithm~\ref{alg:SGD_wor}). This strategy implies that stochastic gradient in RR is biased.
\begin{algorithm}[h!]
    \caption{Random Reshuffling (RR)}
    \label{alg:SGD_wor}
    \begin{algorithmic}
        \REQUIRE learning rates $\{h_{s,k}\}_{s,k\ge 0}$, starting point $x^0\in\R^n$, batch size $r \ge 1$, number of epochs $S$
        \STATE Set $x_0^0 = x^0$
        \FOR{ $s=0,1,2,\ldots K-1$ }
        \STATE Generate random permutation $\{i_{s,1},\ldots,i_{s,m}\}$ of the set $\{1,\ldots,m\}$
        \STATE Set $l = \lceil\nicefrac{m}{r}\rceil$
        \FOR{$k = 0,1,\ldots,l-1$}
        \STATE Set $\hat r_s^k = \min\{r,m-kr\}$
        \STATE Compute $g_s^k = \frac{1}{r_s^k}\sum\limits_{j=1}^{r_s^k} \nabla f_{i_{s,kr+j}}(x_s^k)$
        \STATE $x_s^{k+1} = x_s^k - h_{s,k} g_s^k$
        \ENDFOR
        \STATE $x_{s+1}^0 = x_s^l$
        \ENDFOR
        \RETURN $x_{K-1}^{l}$
    \end{algorithmic}
\end{algorithm}

While the superiority of RR to SGD was empirically discovered a long time ago \cite{bottou2009curiously, bottou2012stochastic}, the theoretical justification of this phenomenon was developed only recently \cite{haochen2018random, rajput2020closing, nguyen2020unified, mishchenko2020random}. In particular, authors of \cite{nguyen2020unified} proved that RR under uniformly bounded gradients assumption,
\begin{equation*}
    \|f_i(x)\|_2 \le G \quad \forall i = 1,\ldots,m,\; \forall x\in\R^n,
\end{equation*}
finds $\varepsilon$-stationary point with the rate $O\left(L_{\max}m(f(x^0) - f_*)\left(\varepsilon^{-2} + G\varepsilon^{-3}\right)\right)$, where $L_{\max}$ is the maximal smoothness constant of summands $f_1,\ldots, f_m$. Then, in \cite{mishchenko2020random} this result was generalized and tightened: under the assumption
\begin{equation}
    \frac{1}{m}\sum\limits_{i=1}^m\|\nabla f_i(x) - \nabla f(x)\|_2^2 \le 2A\left(f(x) - f_*\right) + C, \label{eq:rr_assumption}
\end{equation}
which is a special case of \eqref{eq:expected_smoothness} with $B=1$, authors of \cite{mishchenko2020random} derived the following bound:
\begin{equation}
    O\left(L_{\max}\sqrt{m}(f(x^0) - f_*)\left(\frac{\sqrt{m}}{\varepsilon^{2}} + \frac{\sqrt{A(f(x^0)-f_*)} + \sqrt{C}}{\varepsilon^3}\right)\right). \label{eq:rr_sota_bound}
\end{equation}
That is, under uniformly bounded variance assumption \eqref{eq:bounded_var} this bound transforms ($A = 0$, $C = \sigma^2$) into $O\left(L_{\max}\sqrt{m}(f(x^0) - f_*)\left(\sqrt{m}\varepsilon^{-2} + \sigma\varepsilon^{-3}\right)\right)$ which outperforms the corresponding complexity bound for SGD \eqref{eq:SGD_non_cvx_bound} whenever $L_{\max}\sqrt{m}\varepsilon \le L\sigma$. Next, one can show that for $L_{\max}$-smooth $f_i$ uniformly lower bounded by $f_{i,*}$, $i=1,\ldots,m$, \eqref{eq:rr_assumption} holds with $A = L_{\max}$ and $C = 2L_{\max}\Delta_* = \frac{2L_{\max}}{m}\sum_{i=1}^m(f_* - f_{i,*})$, and, as a consequence of \eqref{eq:rr_sota_bound}, RR converges with the rate
\begin{equation*}
    O\left(L_{\max}\sqrt{m}(f(x^0) - f_*)\left(\frac{\sqrt{m}}{\varepsilon^{2}} + \frac{\sqrt{L_{\max}(f(x^0)-f(x^*))} + \sqrt{L_{\max}\Delta_*}}{\varepsilon^3}\right)\right)
\end{equation*}
which is better than corresponding bound for SGD (see Table~\ref{tab:sgd_bounds}) when $L\sqrt{f(x^0)-f_*} \ge \varepsilon\sqrt{L_{\max}m}$ and $L\sqrt{\Delta_*} \ge \varepsilon\sqrt{L_{\max}m}$.

\subsection{Variance-reduced Methods}\label{sec:vr}
In this section, we discuss variance reduction for non-convex optimization -- a special technique aimed at improving the convergence speed of SGD for finite-sum optimization problems \eqref{eq:main_pr_non_cvx}+\eqref{eq:f_finite_sum}. The typical behaviour of SGD with constant stepsize $h$ and \kesha{batch size} $r < m$ is as following: during the first iterations the method converges rapidly to some neighbourhood of the solution or local minimum and then it starts to oscillate in this neighbourhood. Such oscillations of SGD are common even for strongly convex problems meaning that it is not a drawback of the problem. The size of the oscillation region is proportional to $\nicefrac{h\sigma^2}{r}$ and this fact hints two simple and famous remedies: decreasing (gradually or suddenly) or small stepsizes and large enough \kesha{batch sizes}. However, the first option can make the convergence too slow and the second option dramatically increases the iteration cost.

To remove these drawbacks one can apply variance-reduced methods like SAG \cite{schmidt2017minimizing}, SAGA \cite{defazio2014SAGA}, SVRG \cite{johnson2013accelerating}, Finito \cite{defazio2014finito}, MISO \cite{mairal2015incremental}. In particular, all of the mentioned methods have $O\left(\left(m+\nicefrac{L}{\mu}\right)\ln\frac{1}{\varepsilon}\right)$ convergence rate in the $\mu$-strongly convex case. What is more, they use constant stepsize and at each iteration (besides each $m$-th iteration or besides the first one) they require one computation of the stochastic gradient with batch size $r = 1$ in the strongly convex case.

Among variance-reduced methods SAGA and SVRG are the most popular ones (see Algorithm~\ref{alg:SAGA}~and~\ref{alg:SVRG}).
\begin{algorithm}[ht]
    \caption{SAGA \cite{defazio2014SAGA,reddi2016proximal}}
    \label{alg:SAGA}
    \begin{algorithmic}
        \REQUIRE learning rate $h>0$, starting point $x^0\in\R^n$, batch size $r \ge 1$
        \STATE Set $\phi_j^0 = x^0$ for each $j\in[m]$
        \STATE $v^0 = \frac{1}{m}\sum\limits_{i=1}^m\nabla f_i(\phi_j^0)$
        \FOR{ $k=0,1,2,\ldots$ }
        \STATE{Uniformly randomly pick sets $I_k, J_k$ from $\{1,2,\ldots,m\}$ (with replacement) such that $|I_k| = |J_k| = r$}
        \STATE{$g^k = \frac{1}{r}\sum\limits_{i\in I_k} \left(\nabla f_i(x^k) - \nabla f_i(\phi_i^k)\right) + v^k$}
        \STATE{$x^{k+1} = x^k - h g^k$}
        \STATE{$\phi_j^{k+1} = x^{k}$ for $j\in J_k$ and $\phi_j^{k+1} = \phi_j^k$ for $j\not\in J_k$}
        \STATE{$v^{k+1} = v^k - \frac{1}{r}\sum\limits_{j\in J_k}\left(\nabla f_j(\phi_j^k) - \nabla f_j(\phi_j^{k+1})\right)$}
        \ENDFOR
    \end{algorithmic}
\end{algorithm}
\begin{algorithm}[ht]
    \caption{SVRG \cite{johnson2013accelerating,reddi2016proximal}}
    \label{alg:SVRG}
    \begin{algorithmic}
        \REQUIRE learning rate $h>0$, epoch length $T$, starting point $x^0\in\R^n$, batch size $r \ge 1$
        \STATE  $\phi_0 = x_0^0 = x^0$
        \FOR{ $s=0,1,2,\ldots$ }
        \FOR{ $k=0,1,2,\ldots, T-1$ }
        \STATE{Uniformly randomly pick set $I_k$ from $\{1,\ldots, m\}$ (with replacement) such that $|I_k| = r$}
        \STATE{$g^k = \frac{1}{r}\sum\limits_{i\in I_k}\left(\nabla f_i(x_s^k) - \nabla f_i(\phi_s)\right) + \nabla f(\phi_s)$}
        \STATE{$x_s^{k+1} = x_s^k - h g^k$}
        \ENDFOR
        \STATE  $\phi_{s+1} = x_{s+1}^0 = x_s^k$
        \ENDFOR
    \end{algorithmic}
\end{algorithm}
In previous subsections, we already mentioned that to find $\varepsilon$-stationary GD and SGD require\footnote{For simplicity we neglect all parameters except $m$ and $\varepsilon$, see the details in Table~\ref{tab:vr_methods}} $O\left(m\varepsilon^{-2}\right)$ and $O(\varepsilon^{-4})$ calculations of the gradients of the summands respectively. Despite the fact that SAGA and SVRG were initially analysed only in strongly convex cases, now their convergence in non-convex case is also well-known due to \cite{reddi2016proximal,reddi2016stochastic}. Unfortunately, when $r=1$ both SAGA and SVRG guarantee only $O(m\varepsilon^{-2})$ convergence rate as simple GD. However, if $r = m^{\nicefrac{2}{3}}$, then SAGA and SVRG converges with the rate $O(m^{\nicefrac{2}{3}}\varepsilon^{-2})$ which has $m^{\nicefrac{1}{3}}$ times better dependence on $m$ than the complexity bound for GD.

However, the lower bound is $\Omega\left(\sqrt{m}\varepsilon^{-2}\right)$ \cite{fang2018spider, li2020page} and there exist optimal algorithms. Essentially, these methods are variations of SARAH \cite{nguyen2017stochastic}. However, in the original paper on SARAH for non-convex problems authors did not prove complexity bounds for the finite-sum optimization problems. After that, in \cite{fang2018spider} authors proposed the first lower bounds in the small data regime $m = O(L^2(f(x^0)-f^*)\varepsilon^{-4})$ together with the first optimal method called SPIDER. Despite the theoretical optimality of the method, it requires very small stepsize (proportional to $\varepsilon^{-1}$) that leads to the poor behaviour in practice. Moreover, the original proof of the convergence rate for SPIDER is technically tough and, because of it, it is hard to generalize the method for the composite optimization problems. In recent works \cite{wang2018spiderboost,wang2019spiderboost}, much simpler optimal method called SpiderBoost was proposed (see Algorithm~\ref{alg:SpiderBoost}). Moreover, this method works with big constant stepsizes (of order $L^{-1}$), can be easily generalized for the composite optimization problems, and works well with heavy-ball momentum.

\begin{algorithm}[h!]
    \caption{SpiderBoost \cite{wang2018spiderboost,wang2019spiderboost}}
    \label{alg:SpiderBoost}
    \begin{algorithmic}
        \REQUIRE learning rate $h>0$, epoch length $T$, starting point $x^0\in\R^n$, batch size $r \ge 1$, number of iterations $K$
        \FOR{ $k=0,1,2,\ldots$ }
        \IF{$k\mod T = 0$}
        \STATE Compute $g^k = \nabla f(x^k)$
        \ELSE
        \STATE{Uniformly randomly pick set $I_k$ from $\{1,\ldots, m\}$ (with replacement) such that $|I_k| = r$}
        \STATE{Compute $g^k = \frac{1}{r}\sum\limits_{i\in I_k}\left(\nabla f_i(x^k) - \nabla f_i(x^{k-1})\right) + g^{k-1}$}
        \ENDIF
        \STATE{$x^{k+1} = x^k - h g^k$}
        \ENDFOR
        \STATE{Pick $\xi$ uniformly at random from $\{0,\ldots, K-1\}$}
        \RETURN $x^{\xi}$
    \end{algorithmic}
\end{algorithm}

Next, in \cite{li2020page}, the same lower bound $\Omega\left(\sqrt{m}\varepsilon^{-2}\right)$ was derived without any assumptions on $m$. Furthermore, authors of \cite{li2020page} proposed a new optimal method called PAGE (see Algorithm~\ref{alg:PAGE}) which is a variant of SPIDER with random length of the inner loop making the method easier to analyze.

\begin{algorithm}[ht]
    \caption{ProbAbilistic Gradient Estimator (PAGE) Algorithm \cite{li2020page}}
    \label{alg:PAGE}
    \begin{algorithmic}
        \REQUIRE initial point $x^0$, stepsize $h$, minibatch size $r, \; r' < r$, probabilities $\{p_k\}_{k\ge 0} \in (0,1]$ of large-batch stochastic gradient computation, number of iterations $K$
        \STATE  $g^0=\frac{1}{r}\sum\limits_{i \in I_0}\nabla f_i(x^0)$, where $I_0$ denotes indices in the minibatch, $|I_0|=r$
        \FOR{ $k=0,1,2,\ldots, K-1$ }
        \STATE  $x^{k+1} = x^k - h g^k$
        \STATE $g^{k+1} = \begin{cases}
                               \frac{1}{r}\sum\limits_{i \in I_k}\nabla f_i (x^{k+1}) & \text{with probability} \; p_k,\\
                               g^k + \frac{1}{r'}\sum\limits_{i \in I_k'}\left( \nabla f_i(x^{k+1}) - \nabla f_i(x^{k}) \right) & \text{with probability} \; 1-p_k,
                         \end{cases}$ where $|I_k| = r$, $|I_k'| = r'$
        \ENDFOR
        \RETURN $\hat{x}^K$ chosen uniformly from $\{x^k\}_{k = 0}^{K}$
    \end{algorithmic}
\end{algorithm}

However, in deep neural networks training, variance-reduced methods work typically worse than SGD or SGD with momentum \cite{defazio2019ineffectiveness}. This happens often due to the bad behaviour of variance-reduced methods with several widespread in deep learning tricks like batch normalization, data augmentation and dropout (see the details in \cite{defazio2019ineffectiveness}). Moreover, if the model is over-parameterized or, in particular, expected strong growth condition \eqref{eq:expected_strong_growth} or its relaxed version \eqref{eq:relaxed_growth} with small noise level hold, SGD is as fast as GD in terms of iteration complexity, meaning that variance reduction is superfluous. That is, variance reduction trick is often not needed or gives worse rates than the rate of SGD for over-parameterized models from theoretical and practical perspectives. Nevertheless, when the problem is not over-parameterized, it makes sense to use variance-reduced methods.

We summarize the discussed above complexity bounds in Table~\ref{tab:vr_methods}.
\begin{table}[h!]
    \centering
    \normalsize
    \begin{tabular}{|c|c|c|}
        \hline
        Method & Citation & Complexity \\
        \hline
        \makecell{Lower bound} & \makecell{\cite{fang2018spider,li2020page}} & $L\Delta_0\min\{\sigma\varepsilon^{-3},\sqrt{m}\varepsilon^{-2}\}$\\
        \hline
        GD &  & $mL\Delta_0\varepsilon^{-2}$\\
        \hline
        \makecell{SGD,\\ bounded var.} & \cite{ghadimi2013stochastic} & $L\Delta_0\max\{\varepsilon^{-2},\sigma^2\varepsilon^{-4}\}$\\
        \hline
        \makecell{SGD,\\ unbounded var.} & \cite{khaled2020better} & $\frac{L^2\Delta_0}{\varepsilon^4}\max\left\{\Delta_0,\Delta_*\right\}$\\
         \hline
        SVRG, $r=1$ & \cite{reddi2016proximal} & $ mL\Delta_0\varepsilon^{-2}$\\
        \hline
        SVRG, $r=\lceil m^{\nicefrac{2}{3}}\rceil$ & \cite{reddi2016proximal} & $m^{\nicefrac{2}{3}}L\Delta_0\varepsilon^{-2}$\\
        \hline
        SAGA, $r=1$ & \cite{reddi2016proximal} & $ mL\Delta_0\varepsilon^{-2}$\\
        \hline
        SAGA, $r=\lceil m^{\nicefrac{2}{3}}\rceil$ & \cite{reddi2016proximal} & $m^{\nicefrac{2}{3}}L\Delta_0\varepsilon^{-2}$\\
        \hline
        SpiderBoost & \cite{wang2018spiderboost,wang2019spiderboost} & $m^{\nicefrac{1}{2}}L\Delta_0\varepsilon^{-2}$\\
        \hline
        SpiderBoost-M & \cite{wang2019spiderboost} & $m^{\nicefrac{1}{2}}L\Delta_0\varepsilon^{-2}$\\
        \hline
        SPIDER & \cite{fang2018spider} & $L\Delta_0\min\{\sigma\varepsilon^{-3},\sqrt{m}\varepsilon^{-2}\}$\\
        \hline
        PAGE & \cite{li2020page} & $L\Delta_0\min\{\sigma\varepsilon^{-3},\sqrt{m}\varepsilon^{-2}\}$\\
        \hline
    \end{tabular}
    \caption{Overview of the complexity results for different variance-reduced methods applied to solve problem \eqref{eq:main_pr_non_cvx}+\eqref{eq:f_finite_sum} with $L$-smooth summands. The column ``Complexity'' contains an overall number of stochastic first-order oracle calls needed to find $\varepsilon$-stationary point neglecting constant factors. Notation: $\Delta_0 = f(x^0) - \eduard{f_*}$, $\Delta_* = \frac{1}{m}\sum_{i=1}^m(f_* - f_{i,*})$, $\sigma^2 =$ a uniform bound for the variance of the stochastic gradient \eqref{eq:bounded_var} (can be $\infty$ for variance-reduced methods), $r$ = \kesha{batch size}.}
    \label{tab:vr_methods}
\end{table}
We also want to mention some papers not presented in Table~\ref{tab:vr_methods} but being highly relevant. In \cite{li2020unified}, there was developed the generalization of the approach from \cite{khaled2020better} providing a unified analysis of different variants of SGD, non-optimal variance-reduced methods like SAGA or L-SVRG \cite{hofmann2015variance, kovalev2020don}, and some distributed methods with quantization \cite{alistarh2017qsgd} including DIANA-type variance reduction \cite{mishchenko2019distributed,horvath2019stochastic} for non-convex optimization. Next, for the online case \eqref{eq:main_pr_non_cvx}+\eqref{eq:f_expectation} with smooth stochastic trajectories the optimal rate $O(\varepsilon^{-3})$ was shown for STOchastic Recursive Momentum (STORM) method \cite{cutkosky2019momentum}, which does not require periodical large-batch stochastic gradient computations and is more robust to the parameters selection, and for its proximal variant \cite{xu2020momentum}. These results shade a light on the role of momentum in the stochastic first-order methods. Finally, it is optimal to generalize SPIDER and get similar rates for composition optimization problems \cite{zhang2020stochastic, chen2020momentum}.

\subsubsection{Convex and Weakly Convex Sums of Non-Convex Functions}
There are also several results devoted to the case when the objective function $f$ from \eqref{eq:f_finite_sum} is (strongly) convex or almost convex, while the summands $f_i$ are smooth, but can be non-convex. In particular, \cite{zhou2019lower} establishes the lower bounds for the cases when \textbf{(i)} $f$ is $\mu$-strongly convex with $\mu \ge 0$, \textbf{(ii)} $f$ is $\alpha$-weakly convex
\begin{equation*}
    f(x) - f(y) - \langle\nabla f(y),x-y \rangle \ge -\frac{\alpha}{2}\|x-y\|_2^2,
\end{equation*}
and \textbf{(iii)} $f_i$ are $\alpha$-weakly convex. Due to the additional assumptions on the structure of non-convexity in the problem the proposed lower bounds are tighter in these situations than the lower bound from \cite{fang2018spider,li2020page}. \eduard{The lower bounds for the case (i) were further tightened in \cite{xie2019general}.} Moreover, there exist optimal \eduard{and almost optimal} methods for each case\eduard{, see Table~\ref{tab:finite_sums_methods} for the details.} 
\eduard{
\begin{table}[h!]
    \centering
    \normalsize
    \begin{tabular}{|c|c|c|}
        \hline
        Settings & Lower Bound & \makecell{Upper Bound,\\ Methods} \\
        \hline
        \makecell{$f$ is $\mu$-str.\ cvx.\\ and $L$-smooth,\\ $\{f_i\}$ are average $L$-smooth} & 
       \makecell{$(m + m^{\nicefrac{3}{4}}\sqrt{\frac{L}{\mu}})\log\frac{\Delta_0}{\varepsilon}$,\\ \\ \cite{xie2019general}} 
        & \makecell{$(m + m^{\nicefrac{3}{4}}\sqrt{\frac{L}{\mu}})\log\frac{\Delta_0}{\varepsilon}$,\\ Dual-Free SDCA \cite{shalev2016sdca},\\  KatyushaX \cite{allen2018katyusha}}\\
        \hline
        \makecell{$f$ is cvx.\\ and $L$-smooth,\\ $\{f_i\}$ are average $L$-smooth} & \makecell{$m + m^{\nicefrac{3}{4}}\sqrt{\frac{LR_0^2}{\varepsilon}}$,\\ \\ \cite{zhou2019lower}} &  \makecell{$m + m^{\nicefrac{3}{4}}\sqrt{\frac{LR_0^2}{\varepsilon}}$,\\Dual-Free SDCA \cite{shalev2016sdca},\\  KatyushaX \cite{allen2018katyusha}}\\
        \hline
        \makecell{$f$ is $\alpha$-weakly cvx.\\ and $L$-smooth,\\ $\{f_i\}_{i=1}^m$ are average $L$-smooth} & \makecell{$\frac{\Delta_0}{\varepsilon^2}\min\left\{m^{\nicefrac{3}{4}}\sqrt{\alpha L},\sqrt{m}L\right\}$,\\ \\\cite{zhou2019lower}\\ \\} &  \makecell{$\frac{\Delta_0}{\varepsilon^2}\min\left\{m^{\nicefrac{3}{4}}\sqrt{\alpha L},\sqrt{m}L\right\}$,\\RepeatSVRG \cite{carmon2018accelerated,agarwal2017finding},\\ SPIDER \cite{fang2018spider},\\ SNVRG \cite{zhou2018stochastic}}\\
        \hline
        \makecell{$\{f_i\}_{i=1}^m$ are $\alpha$-weakly cvx.\\ and $L$-smooth} & \makecell{$\frac{\Delta_0}{\varepsilon^2}\min\left\{\sqrt{m\alpha L},L\right\}$,\\ \\\cite{zhou2019lower}\\ \\} &  \makecell{$\frac{\Delta_0}{\varepsilon^2}\min\left\{\sqrt{m\alpha L},\sqrt{m}L\right\}$,\\Natasha \cite{allen2017natasha},\\ RapGrad \cite{lan2019accelerated},\\ StagewiseKatyusha \cite{chen2018variance}}\\
        \hline
    \end{tabular}
    \caption{\eduard{Overview of the optimal convergence results for convex and weakly convex sums of non-convex functions. Averaged $L$-smoothness of $\{f_i\}_{i=1}^m$ means that for all $x,y\in\R^n$ the following inequality holds: $\frac{1}{m}\sum_{i=1}^m\|\nabla f_i(x) - \nabla f_i(y)\|_2^2 \le L^2\|x-y\|_2^2$. The column ``Lower Bound'' states for the number of stochastic first-order oracle calls needed to find such $\hat x$ that $\EE[f(\hat x) - f(x^*)] \le \varepsilon$ for the second and the third rows and $\varepsilon$-stationary point for the fourth and the fifth rows. Notation: $R_0$ = distance from $x^0$ to the solutions set (for the third row), $\Delta_0 = f(x^0) - f_*$.}}
    \label{tab:finite_sums_methods}
\end{table}}

\subsection{Adaptive Methods}\label{sec:adaptive_methods}
One of the most significant issues of the methods described above is that they require tuning of the stepsize and other parameters (e.g., \kesha{batch size}) when used in practice. It is often challenging and takes a lot of time, especially for training deep neural networks. That is why, in the recent few years, adaptive methods gained a lot of attention. Below we discuss the most popular ones -- AdaGrad and Adam -- as well as their variants. In fact, all of these methods depend on some parameters, but these algorithms are much more robust than other variants of SGD or variance-reduced methods. Therefore, they are often called adaptive. One can find PyTorch implementation of many popular adaptive first-order methods together with with visualization of their convergence on Rosenbrock and Rastrigin functions in \cite{pytorchOpt}.

\subsubsection{AdaGrad and Adam}
\noindent\textbf{AdaGrad.} As we mentioned above, SGD requires the tuning of the stepsize. The first algorithm aiming to remove this drawback of SGD was AdaGrad \cite{duchi2011adaptive}:
\begin{equation}
    x_i^{k+1} = x_i^k - \frac{h}{\sqrt{G_i^k + \delta}}g_i^k, \label{eq:adagrad}
\end{equation}
where the subscript $i$ denotes the $i$-th component of the vector, $G_i^k = \sum_{t=0}^k(g_i^t)^2$, and $\delta$ is some small positive number preventing from the division by zero and typically taken of the order $10^{-8}$. AdaGrad can be considered as a special case of SGD with different per-coordinate stepsizes.

The main advantage of AdaGrad is in its robustness to the choice of $h$: in practice, it often works well with the default value $h = 10^{-2}$. Moreover, AdaGrad was shown to work well with sparse data \cite{duchi2013estimation}. However, in the dense settings AdaGrad stepsizes rapidly decrease which leads to the slow convergence of the method \cite{wilson2017marginal}.

\noindent\textbf{Adam.} To resolve this issue of AdaGrad one can use exponential moving averages instead of sums $G_i^k$ leading to the method called RMSprop \cite{tieleman2012lecture}. Then, based on RMSprop authors of \cite{kingma2014adam} proposed one the most popular methods in deep learning -- Adam\footnote{To distinguish exponents from superindexes we use braces $(\cdot)$ for exponents.}:
\begin{eqnarray}
    m_i^k &=& \beta_1 m_i^{k-1} + (1-\beta_1) g_i^k,\quad \hat{m}_i^k = \frac{m_i^k}{1-(\beta_1)^k},\notag\\
    v_i^k &=& \beta_2 v_i^{k-1} + (1-\beta_2)(g_i^k)^2,\quad \hat{v}_i^k = \frac{v_i^k}{1-(\beta_2)^2}, \notag\\
    x_i^{k+1} &=& x_i^k - \frac{h}{\sqrt{\hat{v}_i^k}+\delta}\hat{m}_i^k,\quad i=1,\ldots,n, \label{eq:adam}
\end{eqnarray}
$\delta$ is some small positive number preventing from the division by zero and typically taken of the order $10^{-8}$. Default values $\beta_1 = 0.9$ and $\beta_2 = 0.999$ from the original paper \cite{kingma2014adam} often make Adam work well in practice. Adam was initially analyzed in the online convex case, but then authors of \cite{reddi2019convergence} found out the flaw in the proof for Adam and proposed a convergent variant of Adam called AMSGrad.

\noindent\textbf{Convergence Guarantees.} While the superiority of AdaGrad and Adam in comparison to SGD was noticed in many application \cite{duchi2013estimation,lacroix2018canonical,Goodfellow-et-al-2016}, the best-known complexity bounds for AdaGrad, Adam, and their modifications are the same or even worse than ones for SGD \cite{chen2018convergence, zhou2018convergence, zaheer2018adaptive, ward2019adagrad,defossez2020convergence}. Furthermore, these complexity results in non-convex case under more restrictive assumption, e.g., uniformly bounded second moment of the stochastic gradient, than their counterparts for SGD. Among other works providing complexity results for Adam and AdaGrad in the non-convex case we emphasize \cite{defossez2020convergence} because of the generality and the simplicity of the proofs. Moreover, the unified analysis of proximal variants of AdaGrad and Adam was proposed in \cite{yun2020general}. \eduard{Furthermore, we emphasize the recent work \cite{shi2020rmsprop} where authors analyse RMSprop without assuming uniform boundedness of the gradients.}

Next, in \cite{zhang2019adam} the theoretical and empirical study why Adam sometimes behaves significantly better than SGD was conducted. The authors of \cite{zhang2019adam} empirically discovered that Adam performs better than SGD when stochastic gradients are heavy-tailed and the reason is that Adam does an ``adaptive gradient clipping'' \cite{goodfellow2016deep, gorbunov2020stochastic, mikolov2012statistical,pascanu2013difficulty,usmanova2017master,gorbunov2021near-optimal}. In the same work \cite{zhang2019adam} authors showed that in such situations SGD can fail to converge while clipped-SGD (with general and coordinate-wise clipping operators) provably converges to $\epsilon$-stationary point. Moreover, in \cite{zhang2019gradient} it was shown that Gradient Descent with clipping converges even under weaker assumption than $L$-smoothness in the non-convex case with the rate $\sim\epsilon^{-2}$ while Gradient Descent in the same settings can converge arbitrary slower. Then, the bound from \cite{zhang2019gradient} was improved in \cite{zhang2020improved}. Finally, it is known \cite{goodfellow2016deep} that clipped-SGD works better than SGD in the vicinity of extremely steep cliffs. A very similar approach based on the normalization of Gradient Descent was also studied in \cite{hazan2015beyond,levy2016power}.

 \subsubsection{Adaptive SGD}\label{sec:adaptive_sgd}

The approach described in Section \ref{sec:sgd_ubv} for general stochastic optimization problem \eqref{eq:main_pr_non_cvx} with the objective given as \eqref{eq:f_expectation} was recently extended in \cite{dvinskikh2020line-search} to obtain adaptive methods with Armijo-type line-search for stochastic non-convex optimization. 
To do that they consider Algorithm \ref{alg:opt_non_cvx_methods} with the mini-batch stochastic gradient \eqref{eq:mini-batch_stoch_grad} and mini-batch size $r=\max\{1, \nicefrac{8\sigma_0^2}{\eps^{2}}\}$, where $\sigma_0 \geqslant \sigma$. In each iteration $k$ of Algorithm \ref{alg:opt_non_cvx_methods} the stepsize is taken as $h_k=\nicefrac{1}{L_{k}}:=\nicefrac{1}{(2^{i_k-1}L_{k-1})}$ by increasing $i_k \geqslant 0$ until the inequality
\begin{align}\label{eq:choose_L_non_conv3}
            f(x^{k+1}) &\leq f(x^k) + \left\la \frac{1}{r}\sum_{l=1}^r  \nabla f(x,\xi_l), x^{k+1} - x^k \right\ra + L_{k}\|x^{k+1} - x^k \|^2_2 + \frac{\eps^2}{32L_{k}}
\end{align}
is satisfied. This inequality is an inexact upper quadratic bound which follows for sufficiently large $L_k$ from the $L$-smoothness and bounded variance. Thus, $L_k$ plays the role of a guess of the Lipschitz constant $L$ locally between the points $x^k$ and $x^{k+1}$. The authors of \cite{dvinskikh2020line-search} propose also methods for convex problems based on the same idea with the difference that in the convex case the mini-batch size $r$ depends on the iteration counter $k$. Careful choice of this dependence allows to simultaneously adaptively choose both the stepsize $h_k$ and the mini-batch size $r_k$. These methods have the same, up to logarithmic factors, iteration complexity and total number of stochastic oracle calls as their non-adaptive counterparts. In particular, for the non-convex case the iteration complexity to obtain $\varepsilon$-stationary point is $\widetilde{O}\left(\nicefrac{L(f(x^0) - f_*)}{\varepsilon^2}\right)$ and the oracle complexity is $\widetilde{O}\left(L\left(f(x^0)-f_*\right)\max\left\{\nicefrac{1}{\varepsilon^2},\nicefrac{\sigma^2}{\varepsilon^4}\right\}\right)$.
Moreover, empirically, the methods designed for convex problems turned out to be more efficient on non-convex problems than the method designed for non-convex problems.

\section{First-Order Methods under Additional Assumptions}

In the previous parts of the paper, we focused on general non-convex problems. In this section, we consider two subclasses of non-convex objective functions which satisfy assumptions weaker than convexity and, at the same time, strong enough to obtain good global convergence rates of optimization algorithms. For simplicity, we consider an unconstrained optimization problem \eqref{eq:problem_statement} with $Q=\R^n$.

\subsection{Polyak--Łojasiewicz Condition}
A function $f(x)$ is said to satisfy the Polyak--Łojasiewicz (PŁ) condition \cite{polyak1963gradient,lojasiewicz1963topological} (or to be gradient dominated) if for all $x\in\mathbb{R}^n$
\begin{equation}
    f\left( x \right)-f\left( {x^\ast } \right)\leqslant \frac{1}{2\mu }\left\| 
{\nabla f\left( x \right)} \right\|_2^2. \label{eq:PL_cond}
\end{equation}
 This condition implies that any stationary point of $f(x)$ is a global minimum, although it is not necessarily unique. In particular, this property holds for strongly convex functions. It was first shown  %by Polyak in 1963 
in \cite{polyak1963gradient} that if the objective is also $L$-smooth, 
then gradient descent linearly converges to a global minimum, i.e.,
\[
f( {x^k} ){\kern 1pt}-{\kern 1pt}f\left( {x^\ast } 
\right)\;{\kern 1pt}\leqslant {\kern 1pt}\;\exp \left( {-\frac{\mu }{L}k} 
\right){\kern 1pt}\left( {f\left( {x^0} \right){\kern 1pt}-{\kern 
1pt}f\left( {x^\ast } \right)} \right).
\]
The Polyak--Łojasiewicz condition is naturally satisfied for the problems of solving nonlinear systems of equalities $g(x)=0$, 
where $g(x)$ is a vector-valued function. This problem can be equivalently reformulated as
\[
\mathop  {\min }\limits_{x\in {\mathbb{R}}^n} \left\{f\left( x \right)=\frac{1}{2}\left\| {g\left( x \right)} \right\|_2^2\right\}
 .
\]
Assuming that, for all $x \in \R^n$ 
$$\lambda _{\min } \left( J_g(x)J^T_g(x)\right)\ge    \mu >0     ,$$
    where $J_g(x)$ is the Jacobian matrix of $g(x)$, one can show that
\[
\left\|{\nabla f\left( x \right)}\right\|^2=\|J^T_g(x)g(x)\|^2\geqslant\mu\|g(x)\|^2=2\mu f(x),
\] 
which is exactly the Polyak--Łojasiewicz condition since $g(x^*)=0$. An extensive \eduard{survey} of first-order optimization methods under this condition, as well as its relationship with other classes of functions, can be found in \cite{karimi2016linear}. An interesting example of the emergence of PŁ condition in Linear Feedback Control theory was recently described in \cite{fatkhullin2020optimizing} and in over-parameterized deep learning in \cite{belkin2021fit}.

Next, consider the convergence of gradient descent under the PŁ condition in terms of relative accuracy $\widetilde{\nabla} f(x)$
\[
\|\widetilde{\nabla} f(x) - \nabla f(x)\|_2 \le \alpha \|\nabla f(x)\|_2,
\]
where $\alpha \in [0,1)$. Let the stepsize $h$ in gradient descent 
\[
x^{k+1}=x^k - h \widetilde{\nabla} f(x^k)
\]
be computed using the following formula:
\[
h = \frac{1}{L}\frac{1-\alpha}{(1+\alpha)^2}.
\]
Combining this with the Lipschitz condition, we obtain
\[
f(x^{k+1}) \le f(x^k) - \frac{1}{2L}\frac{(1-\alpha)^2}{(1+\alpha)^2}\|\nabla f(x^k)\|^2_2,
\]
leading to
\[
f(x^N) - f(x^\ast) \le \left( 1 - \frac{\mu}{L}\frac{(1-\alpha)^2}{(1+\alpha)^2} \right)^N \left( f(x^0) - f(x^\ast)\right).
\]
As a result, we achieve a linear convergence rate for the gradient descent under the PŁ condition.

In general case the main ingredient that guaranties global linear convergence under PŁ condition is an estimate like $$\|\nabla f(x^N)\|_2^2 \le \theta(N)\cdot\left(f(x^0) - f(x^*)\right),$$ where $\theta(N)$ -- some decreasing function, i.e. \eqref{eq:nonconv_compl_grad}. We assume that there exists such $N(\mu)$, that $\theta\left(N(\mu)\right)\le \mu$, i.e. for \eqref{eq:nonconv_compl_grad} $N(\mu) = 2L/\mu$. In this case from PŁ condition $$\|\nabla f(x^N(\mu))\|_2^2 \le \frac{1}{2}\|\nabla f(x^0)\|_2^2.$$
By applying restarts we obtain oracle complexity $\tilde{O}\left(N(\mu)\right)$.

\subsubsection{Stochastic First-Order Methods under Polyak--Łojasiewicz Condition}
The majority of the methods described in Section~\ref{sec:stoch_methods} are analyzed under PŁ condition as well. That is, one can find the state-of-the-art results for different variants of SGD and non-accelerated variance reduced methods like SVRG and SAGA in \cite{li2020unified}, accelerated variance reduced methods like PAGE in \cite{li2020page}, the tightest known analysis of Random Reshuffling under PŁ condition in \cite{ahn2020sgd}, and the convergence results for SGD in the over-parameterized case with constant, Armijo-type, and stochastic Polyak's stepsizes in \cite{vaswani2018fast}, \cite{vaswani2019painless}, and \cite{loizou2020stochastic} respectively. \eduard{The summary of known complexity results for the stochastic methods under PŁ condition is given in Table~\ref{tab:PL_condition_summary}. We emphasize that the analysis from \cite{gower2020sgd} is derived under so-called \textbf{expected residual (ER)} assumption on the stochastic gradient $g(x)$: there exists such constant $\rho >0$ that
\begin{equation}
    \EE\left[\left\|g(x)-g(x^*) - \left(\nabla f(x) - f(x^*)\right)\right\|_2^2\right] \le 2\rho\left(f(x) - f(x^*)\right).\label{eq:expected_residual}
\end{equation}
Moreover, in the analysis of Random Reshuffling from \cite{ahn2020sgd} it is used that the norms of the gradients of individual functions from the sum \eqref{eq:f_finite_sum} are uniformly upper by some constant $G$ on the sublevel set:
\begin{equation}
    \|\nabla f_i(x)\|_2 \le G,\quad \forall i=1,\ldots,m,~~~ \forall x\in\R^n: f(x) \le f(x^0). \label{eq:uniform_upper_bounded_grad}
\end{equation}}

Rather simple introduction (close to the state of the art results) for SGD with bias under PŁ condition can be find in \cite{ajalloeian2020analysis}.

\begin{table}[h!]
    \centering
    % \normalsize
    \begin{threeparttable}
    \begin{tabular}{|c|c|l c r|}
         \hline
         Problem & Method & Citation & Complexity & Assumptions\\ 
\hline\hline
    \makecell{\eqref{eq:main_pr_non_cvx}} & GD & \cite{polyak1963gradient} & $\frac{L}{\mu}\log\frac{\Delta_0}{\varepsilon}$ & \\
    \hline\hline
    \multirow{4}{1cm}{\eqref{eq:main_pr_non_cvx}+\eqref{eq:f_expectation}} &\multirow{2.5}{1cm}{\centering SGD} &\cite{khaled2020better,karimi2016linear} & $\frac{L}{\mu}\log\left(\frac{\Delta_0}{\varepsilon}\right) + \frac{L\sigma^2}{\mu^2\varepsilon}$ & UV \eqref{eq:bounded_var}\\
    & & \cite{vaswani2018fast,khaled2020better} & $\frac{\alpha L}{\mu}\log\left(\frac{\Delta_0}{\varepsilon}\right) + \frac{L\beta}{\mu^2\varepsilon}$ & RG \eqref{eq:relaxed_growth} \\
    \cline{2-5}
    & PAGE & \cite{li2020page} & $\left(\frac{\sigma^2}{\mu\varepsilon} + \sqrt{\frac{\sigma^2}{\mu\varepsilon}}\frac{L_{\text{avg}}}{\mu}\right)\log\left(\frac{\Delta_0}{\varepsilon}\right)$ & \begin{tabular}{r}
        UV \eqref{eq:bounded_var},\\ Avg.\ $L_{\text{avg}}$-smth. 
    \end{tabular}\\
    \hline\hline
    \multirow{20}{1cm}{\eqref{eq:main_pr_non_cvx}+\eqref{eq:f_finite_sum}}& GD & \cite{polyak1963gradient} & $m\frac{L}{\mu}\log\frac{\Delta_0}{\varepsilon}$ &\\
    \cline{2-5}
    &\multirow{5.5}{1cm}{\centering SGD} & \cite{khaled2020better} & $\frac{L}{\mu}\left(\left(\frac{A}{\mu}+B\right)\log\left(\frac{\Delta_0}{\varepsilon}\right)+\frac{C}{\mu\varepsilon}\right)$ & ES \eqref{eq:expected_smoothness} \\
    & &\cite{vaswani2018fast,khaled2020better} & $\frac{\alpha L}{\mu}\log\left(\frac{\Delta_0}{\varepsilon}\right) + \frac{L\beta}{\mu^2\varepsilon}$ & RG \eqref{eq:relaxed_growth} \\
    & & \cite{khaled2020better} & $\frac{L}{\mu}\left(\frac{\max_i L_i}{\mu}\log\left(\frac{\Delta_0}{\varepsilon}\right)+\frac{\max_i L_i\Delta_*}{\mu\varepsilon}\right)$ & Unif.\ sampl. \\
    & & \cite{khaled2020better} & $\frac{L}{\mu}\left(\frac{\overline{L}}{\mu}\log\left(\frac{\Delta_0}{\varepsilon}\right)+\frac{\overline{L}\Delta_*}{\mu\varepsilon}\right)$ & Imp.\ sampl. \\
    & &\cite{gower2020sgd} & $\frac{L}{\mu}\left(\left(\frac{\rho}{\mu}+1\right)\log\left(\frac{\Delta_0}{\varepsilon}\right)+\frac{\sigma_*^2}{\mu\varepsilon}\right)$ & ER \eqref{eq:expected_residual} \\
    &\begin{tabular}{c}
        + Armijo\\
        line-search 
    \end{tabular} & \cite{vaswani2019painless} & $\left(\frac{\alpha L}{\mu} + \frac{\max_i L_i}{\mu}\right)\log\left(\frac{\Delta_0}{\varepsilon}\right)$ & E-SG \eqref{eq:expected_strong_growth} \\
    &\begin{tabular}{c}
        + Polyak\\
        stepsizes 
    \end{tabular} & \cite{loizou2020stochastic} & $\frac{\max_i L_i^2}{\mu^2}\log\left(\frac{\Delta_0}{\varepsilon}\right)$ & Interpolation \eqref{eq:interpolation_cond} \\
    \cline{2-5}
    & RR & \cite{ahn2020sgd} & $\left(\frac{m\Delta_0}{\varepsilon} + \frac{mL^2G^2\log^3(\varepsilon^{-1})}{\mu^3\varepsilon}\right)^{\nicefrac{1}{2}}$ & Bounded gradients \eqref{eq:uniform_upper_bounded_grad}\\
    \cline{2-5}
    & SVRG &
        \cite{reddi2016stochastic,reddi2016proximal}
    & $\left(m+\frac{m^{\nicefrac{2}{3}}\max_i L_i}{\mu}\right)\log\left(\frac{\Delta_0}{\varepsilon}\right)$
    & \\
    \cline{2-5}
    & \begin{tabular}{c}
        L-SVRG\\
        SAGA
    \end{tabular} & \cite{li2020unified,reddi2016proximal} & $\left(m+\frac{m^{\nicefrac{2}{3}}L}{\mu}\right)\log\left(\frac{\Delta_0}{\varepsilon}\right)$ & Avg.\ $L_{\text{avg}}$-smth.\\
    \cline{2-5}
    & PAGE & \cite{li2020page} & \begin{tabular}{c}
        $\left(b + \sqrt{b}\frac{L_{\text{avg}}}{\mu}\right)\log\left(\frac{\Delta_0}{\varepsilon}\right)$,\\
        where $b = \min\{\frac{\sigma^2}{\mu\varepsilon},m\}$ 
    \end{tabular}& \begin{tabular}{r}
        UV \eqref{eq:bounded_var}\\ with $\sigma^2 \le +\infty$, \\ Avg.\ $L_{\text{avg}}$-smth. 
    \end{tabular}\\
    \hline
    \end{tabular}
    \end{threeparttable}
    \caption{\eduard{Summary of the state-of-the-art complexity results for different stochastic first-order methods under assumption that $f$ is $L$-smooth and satisfies PŁ condition \eqref{eq:PL_cond}. Columns: ``Complexity'' -- an overall number of stochastic first-order oracle calls needed to find such $\hat x$ that $\EE[f(\hat x) - f(x^*)] \le$ neglecting constant factors; ``Assumptions'' -- the assumptions used to derive the corresponding complexity bound in addition to $L$-smoothness of $f$ and PŁ condition \eqref{eq:PL_cond} for $f$. For finite-sum case \eqref{eq:f_finite_sum} it is additionally assumed that each $f_i$ is $L_i$-smooth, $i=1,\ldots,m$. Abbreviations: UV -- uniform variance bound assumption \eqref{eq:bounded_var}; RG -- relaxed growth condition \eqref{eq:relaxed_growth}; Avg. $L_{\text{avg}}$-smth.\ -- averaged $L$-smoothness assumption meaning that there exist such $L$ that $\EE\left[\|\nabla f(x,\xi) - \nabla f(y,\xi)\|_2^2\right] \le L^2\|x-y\|_2^2$ in the online case \eqref{eq:f_expectation}, and $\EE\left[\|\nabla f_j(x) - \nabla f_j(y)\|_2^2\right] \le L^2\|x-y\|_2^2$ in the finite-sum case, where $j$ is sampled uniformly at random from $\{1,\ldots,m\}$; Unif.\ Sampl.\ and Imp.\ Sampl.\ denote the sampling strategies described in Section~\ref{sec:SGD_assumptions}. Notation: $\Delta_0 = f(x^0) - \eduard{f_*}$; $\sigma^2 =$ a uniform bound for the variance of the stochastic gradient \eqref{eq:bounded_var}; $\alpha,\beta$ = relaxed growth condition parameters; $\Delta_* = \frac{1}{m}\sum_{i=1}^m(f_* - f_{i,*})$; $\max_i L_i$ = maximal smoothness constant of $f_i$ in \eqref{eq:f_finite_sum}; $\overline{L}$ = averaged smoothness constant of $f_i$ in \eqref{eq:f_finite_sum}; $\sigma_*^2 = \EE[\|g(x^*)\|_2^2]$ -- the variance of the stochastic gradient at the solution.}}
    \label{tab:PL_condition_summary}
\end{table}

\subsection{Star-convexity and $\alpha$-weak-quasi-convexity}
A function $f(x)$ is called star-convex if for some global minimizer $x^\ast$ and for all $\lambda \in [0,1]$ and $x\in\mathbb{R}^n$
\[f(\lambda x+(1-\lambda)x^\ast)\leqslant \lambda f(x)+(1-\lambda) f(x^\ast).\] While any interval connecting two points on the graph of a convex function lies not lower than the graph, for a star-convex functions this is assumed only for intervals connecting some fixed global minimizer and any other point on the graph. This condition is considerably weaker than convexity, even for functions of one variable. For example, the function $|x|(1-e^{-|x|})$ is a non-convex star-convex function. 
The authors of \cite{lee2016optimizing} analyze a cutting plane method for minimization of this class of functions and obtain a polylogarithmic in $\eps$ and polynomial in $n$ complexity bound using only function evaluations. The authors of \cite{guminov2019accelerated,nesterov2020primal-dual} prove that the same Algorithm \ref{AGMsDR} possesses the following convergence rate for star-convex $L$-smooth functions
\begin{equation}\label{grad}
\min_{k=[N/2],...,N}\| \nabla f(y^k)\|_*^2 \leqslant \frac{64L^2V[x^0](x^*)}{N^{3}},
\end{equation}
$$
f(x^N) - f(x^{\pd{*}}) \leqslant \frac{4L\pd{V[x^0](x^*)}}{N^2}.
$$

A more general class of functions is the class of $\alpha$-weakly-quasi-convex functions satisfying
\begin{equation}
    f\left( x \right)-f\left( {x^\ast } \right)\leqslant \frac{1}{\alpha}\left \langle{\nabla f\left( x \right),x-x^\ast } \right\rangle \label{eq:alpha_weak_cvx}
\end{equation}
for some $\alpha\in(0,1]$ and some global minimizer $x^\ast$. Continuously differentiable 1-weakly-quasi-convex functions are exactly the star-convex functions. 
The authors of \cite{guminov2017accelerated} propose an algorithm with iteration complexity $O(\alpha^{-1}L^{1/2}R\eps^{-1/2})$, where $R$ is an upper bound on the initial distance to the point $x^*$. A slightly worse bound $O(\alpha^{-3/2}L^{1/2}R\eps^{-1/2})$ is obtained in \cite{nesterov2020primal-dual} by restarting Algorithm \ref{AGMsDR}. Both approaches require a line search for which the complexity also needs to be estimated. The authors of \cite{hinder2020near} analyze this complexity and propose an algorithm with $O(\alpha^{-1}L^{1/2}R\eps^{-1/2})$ iteration complexity and the same up to a logarithmic factor in $\alpha^{-1}\eps^{-1}$ number of function and gradient evaluations. Moreover, they provide a similar lower complexity bound, thus proving that their method is optimal. Further, they also consider a class of $(\alpha,\mu)$-strongly quasi-convex functions satisfying
\[
f\left( x \right)-f\left( {x^\ast } \right)\leqslant \frac{1}{\alpha}\left \langle{\nabla f\left( 
x \right),x-x^\ast } \right\rangle - \frac{\mu}{2}\|x-x^*\|^2 
\] 
and provide an algorithm which has iteration complexity $$O(\alpha^{-1}L^{1/2}\mu^{-1/2}\log(\alpha^{-1}\eps^{-1}))$$ and requires up to a logarithmic factor the same number function and gradient evaluations.
Similar optimal complexity bounds for accelerated gradient method for $\alpha$-weakly-quasi-convex functions  and $(\alpha,\mu)$-strongly quasi-convex functions were obtained in \cite{bu2020note} by extending the estimating sequence technique.

\eduard{
\subsubsection{Stochastic Methods and $\alpha$-weak-quasi-convexity}
The most general analysis of SGD under $\alpha$-weak-quasi-convexity is provided in \cite{gower2020sgd}. As it was mentioned earlier, authors of \cite{gower2020sgd} consider finite-sum optimization problems\footnote{In fact, most of the results from \cite{gower2020sgd} do not rely on the finite-sum structure of $f$. } \eqref{eq:main_pr_non_cvx}+\eqref{eq:f_finite_sum} and derive complexity bounds for SGD under expected residual \eqref{eq:expected_residual} assumption on the stochastic gradient for the $\alpha$-weak-quasi-convex function and functions satisfying PŁ condition. In particular, for SGD in these settings the following bound was established:
\begin{equation}
    O\left(\frac{(\rho+L)R_0^2}{\alpha^2\varepsilon} + \frac{\sigma_*^2 R_0^2}{\alpha^2\varepsilon^2}\right),
\end{equation}
where $\sigma_*^2 = \EE[\|g(x^*)\|_2^2]$ is the variance of the stochastic gradient at the solution. Note, that when interpolation condition \eqref{eq:interpolation_cond} holds this bounds reduces to $ O\left(\nicefrac{(\rho+L)R_0^2}{(\alpha^2\varepsilon)}\right)$. Moreover, under interpolation condition the authors of \cite{gower2020sgd} also derived that the generalized version of stochastic Polyak stepsize \eqref{eq:stochastic_polyak_stepsize} for stochastically reformulated problem \eqref{eq:main_pr_non_cvx}+\eqref{eq:f_finite_sum} converges with the rate
\begin{equation}
    O\left(\frac{\mathcal{L} R_0^2}{\alpha^2\varepsilon} \right),
\end{equation}
where $\mathcal{L}$ is the expected smoothness constant of stochastic reformulation (see the details in \cite{gower2019sgd,gower2020sgd}). In the full-batch case, i.e., when $g(x) = \nabla f(x)$, we have $\mathcal{L} = L$, and in the importance sampling case, i.e., when $g(x) = \nabla f_j(x)$ where $j = i$ with probability $\nicefrac{L_i}{\sum_{t=1}^mL_t}$, we have $\mathcal{L} = \overline{L} = \frac{1}{m}\sum_{i=1}^mL_i$.
}

\subsubsection{\eduard{Further Generalizations}}
A more wide class of functions that covers the class of $\alpha$-weakly-quasi-convex functions referred to as approximately homogeneous functions satisfying the condition 
\[
N (f\left( x \right)-f\left( {x^\ast } \right))\leqslant \left \langle{\partial f\left( 
x \right),x-x^\ast } \right\rangle \leqslant M (f\left( x \right)-f\left( {x^\ast } \right)),
\]
where $\partial f(x)$ is a subgradient of $f(x)$ and $N,M$ are some constants. This class of functions was first defined in \cite{shor1967generalized} and discussed in \cite{polyak1987introduction}.

\alexander{In general, if there exist good lower and upper convex models for non-convex target function, one can derive that complexity of such problem is similar to convex ones rather than non-convex (see  \cite{bazarova2020linearly} and references therein).}

\section{Higher-Order Methods} 
\subsection{Second-Order Methods}
Another branch of optimization incremental methods for solving \eqref{eq:main_pr_non_cvx} are methods that use the second-order information about the function. This information is very helpful to escape saddle-points by using a negative curvature. 
Next we define an $(\varepsilon,\delta)$-second-order stationary point $x^{\ast}$ if
\begin{align*}
    \|\nabla f(x^{\ast}) \|_2 \leq \varepsilon,\quad \lambda_{\min} \left(\nabla^2 f(x^{\ast})\right) \geq - \delta.
\end{align*}

Next in this section we suppose that $f(x)$ has $L_2$-Lipschitz second-order derivative. The basic method for this class of problems is a \dk{Cubic Regularization method (CR)} \cite{nesterov2006cubic}.
\begin{equation}
    \label{eq_cubic_newton}
x^{k+1}_{Cubic} = x^{k} + \argmin_{s\in \R^n} \left[ \nabla f(x^{k})^\top s + \frac{1}{2} s^\top \nabla^2 f(x^{k}) s + \frac{\dk{H}}{6} \|s\|_2^3 \right],
\end{equation}
where $\dk{H} \ge 0$. It globally converges to the minimum for convex functions and converges to a $(\varepsilon,\delta)$-second-order stationary point for non-convex function within $O(\varepsilon^{-3/2})$ number of iterations. Note, that the subproblem \eqref{eq_cubic_newton} is also non-convex but in \cite{nesterov2006cubic} authors proposed a method to solve this problem as a convex problem via special choose of $\dk{H}$ and line-search for a dual problem. 
\pd{A related line of work considers trust region methods \cite{conn2000trust,cartis2011adaptive,cartis2011adaptive2,cartis2017improved,cartis2019universal}, where a classical Newton step is calculated on a Euclidean ball of a carefully chosen radius. 
Both cubic regularized Newton methods and trust region methods can be extended to work for constrained problems  with linear and conic constraints \cite{haeser2019optimality,dvurechensky2019generalized,dvurechensky2021hessian}.
In general, all these algorithms work well for the problems in moderate dimensions.
%It works very good for small problems, but 
Unfortunately, for many large-scale Machine Learning problems it is hard to calculate the full Hessian and the inverse such a large matrix.} Recent work has therefore explored the use of Hessian-vector products $\nabla^2 f(x) \cdot \dk{s}$, which can be computed as efficiently as gradients in many cases including neural networks by using autogradient technique. By this Hessian-vector product we can efficiently find \dk{$x^{k+1}_{\text{Cubic}}$} by variants of gradient descent \cite{carmon2016gradient}.  
Several algorithms incorporating Hessian-vector products \cite{allen2018make,allen2018natasha} have been shown to achieve faster convergence rates than gradient descent in the non-stochastic setting. However, in the stochastic setting where we only have access to stochastic Hessian-vector products, significantly less progress has been made.

One of the improvement of this method was done in \cite{wang2020cubic}. \pd{The authors introduce a momentum step and obtain faster convergence rate}. This technique \pd{is} widely used to speed up the first order methods and also can speed up the second order method.

\begin{algorithm}  
	\caption{CRm}
	\label{Adaptive_version}
	\begin{algorithmic}[1] 
		\STATE {\bfseries Input:} Initialization $x^0 =  {y}^0\in \mathbb{R}^\dk{n},\rho < 1, \dk{H} > L_2$.
		\FOR{$k=0,1,\dots $}
		\STATE {\bfseries Cubic step:}   
		{ 
		\begin{align*}  {s}^{k+1} &=\argmin_{s} \left[\nabla f(x^{k})^\top s  + \frac{1}{2} s^\top \nabla^2 f(x^{k}) s  + \frac{\dk{H}}{6} \|s\|_2^3  \right], \\
		y^{k+1} &= x^{k} + s^{k+1}. 
			\end{align*} 
		}
		\vspace{-.3cm}
		\STATE {\bfseries Momentum step:}   
		\begin{align} 
			\dk{\beta_{k+1}} &=    \min\{\rho , \|\nabla f( {y}^{k+1})\|_2, \| {y}^{k+1} - x^{k}\|_2  \} ,  \\
			 \dk{z}^{k+1} &=  {y}^{k+1} + \dk{\beta_{k+1}} ( {y}^{k+1} -  {y}^{k}).
		\end{align}
		\vspace{-.3cm} 
		\STATE \textbf{Monotone Step:}
		\vspace{-.3cm}
		\begin{align}
		x^{k+1}= \argmin_{x \in \{ {y}^{k+1},  \dk{z}^{k+1}\}} f(x). 
		\end{align}  
		\ENDFOR
	\end{algorithmic} 	 
\end{algorithm}

Also, second-order methods that have access to the Hessian of $f$ can exploit negative curvature to more effectively escape saddles and arrive at local minima. To show this concept we introduce one of such methods \cite{wright2018optimization}.
There are two types of steps: gradient steps and a step in a negative curvature for the Hessian. So

\begin{itemize}
    \item If $\|\nabla f(x^k)\|_2>\varepsilon$, we do gradient step.
    \item Otherwise, if $\lambda_{\min}\left(\nabla^2 f(x^k)\right) < - \delta$, choose $\dk{s}^k$ to be the eigenvector corresponding to $\lambda_{\min}\left(\nabla^2 f(x^k)\right))$ and do step $x^{k+1}=x^{k}+\alpha_k \dk{s}^k$.
\end{itemize}
There are different policies to $\alpha_k$ and gradient steps. The main idea here is to use the first-order methods as a cheap main method and switch to expensive second-order methods when we reach local stationary point and want to escape it to find a better local minimum. Methods with this idea are still developing. In \cite{ge2015escaping, jin2017how} \pd{it} was proved that gradient methods with additive noise are able to escape from nondegenerate saddle points and find approximate local minima. These ideas lead to the state of art first-order methods to find local minima with Hessian-vector product \cite{carmon2018accelerated, royer2018complexity, allen2018natasha, xu2018first, allen2018neon2, jin2018accelerated, fang2018spider, nguyen2017sarah}. In recent works \cite{fang2019sharp,jin2019nonconvex,roy2020escaping} \pd{it} was proved that stochastic gradient descent can escape from saddle point and converges to approximate local minima.

\subsection{Stochastic Second-Order Methods}
Now we move to stochastic version of problem \eqref{alg:opt_non_cvx_methods}. Firstly, we speak about online version \eqref{eq:f_expectation}, where we minimize expectation of some stochastic function. 
In the work \cite{tripuraneni2018stochastic} authors propose a stochastic optimization method that utilizes stochastic gradients and Hessian-vector products to find an $(\varepsilon,\delta)$-second-order stationary point using only $O(\varepsilon^{-3.5})$ oracle evaluations. This rate improves upon the $O(\varepsilon^{-4})$ rate of stochastic gradient descent, and matches the best-known result for finding local minima without the need for any delicate acceleration or variance reduction techniques.

\begin{algorithm}[!h]
\caption{Stochastic Cubic Regularization}\label{algo:SCRN}
\begin{algorithmic}[1]
\REQUIRE mini-batch sizes $\dk{r}_1, \dk{r}_2$, initialization $x_{0}$, number of iterations $\dk{N}$, and final tolerance $\varepsilon$.
 \FOR{$\dk{k} = 0, \hdots, \dk{N}$}
 \STATE Sample $S_1 \leftarrow \{\mathbf{\xi}_{i} \}_{i=1}^{\dk{r}_1}$, $S_2 \leftarrow \{\mathbf{\xi}_{i} \}_{i=1}^{\dk{r}_2}$.
 \STATE $g^\dk{k} \dk{=} \frac{1}{\dk{r_1}} \sum_{\mathbf{\xi}_i \in S_1} \nabla f(x^{k}; \ \mathbf{\xi}_i)$
 \STATE $B^\dk{k} [\cdot] \dk{=} \frac{1}{\dk{r_2}} \sum_{\xi_i\in S_2} \nabla^2 f(x^{k}, \xi_i)(\cdot)$
 \STATE $\dk{s^k = \argmin\limits_{s} \left\{ \psi_k(s) =  s^\top g^k  + \frac{1}{2} s^\top B^k  s + \frac{\dk{L_2}}{6} \|s\|_2^3 \right\}}  $
 \STATE $x^{\dk{k} +1} \leftarrow x^{\dk{k} } + \dk{s^k}$
 %\IF{ $\dk{\psi_k(s^k)} \geq -\frac{1}{100} \sqrt{\frac{\varepsilon^3}{\dk{L_2}}}$}  \label{algo1:line:small_descent}
  %\STATE $\Delta \leftarrow \text{Cubic-Finalsolver}(g^\dk{k} , \ B^\dk{k} [\cdot], \varepsilon)$
  %\STATE $x^* \leftarrow x^{\dk{k}}  + \Delta$
  %\STATE \textbf{brake}
 %\ENDIF
 \ENDFOR
\ENSURE %$x^*$ if the early termination condition was reached, otherwise 
\dk{The final iterate $x_{N + 1}$}.
\end{algorithmic}
\end{algorithm}
\dk{
This is a stochastic cubic regularization algorithm in Algorithm \ref{algo:SCRN}.  To obtain stochastic gradients and Hessians, we can sample independent batches of $S_1 $and $S_2$ in each iteration, but they can also be connected so that $S_2 \subseteq S_1$. The average gradient is denoted by} 
\begin{align}
  g^\dk{k}  = \frac{1}{\dk{r_1}} \sum_{\xi_i \in S_1} \nabla f(x^{k}, \xi_i)
\end{align}
and the average Hessian by
\begin{align}
  B^\dk{k}  = \frac{1}{\dk{r_2}} \sum_{\xi_i \in S_2} \nabla^2 f (x^{k}, \xi_i),
\end{align}
this implies a \textit{stochastic cubic submodel}:
\begin{equation*}
\dk{\psi_{k}(s)=  s^\top g^k  + \frac{1}{2} s^\top B^k  s + \frac{\dk{L_2}}{6} \|s\|_2^3  }
\label{eq:approx_cubic} .
\end{equation*}

%After sampling minibatches for the gradient and the Hessian, Algorithm \ref{algo:SCRN} makes a call to a black-box cubic subsolver to optimize the stochastic submodel $\dk{\psi_{k}(s)}$. The subsolver returns a parameter change $\dk{s}$, i.e., an approximate minimizer of the submodel, along with the corresponding change in submodel value, $\Delta_m = m^{k}(x^{k} + \Delta) - m^{k}(x^{k})$. The algorithm then updates the parameters by adding $\Delta$ to the current iterate, and checks whether $\Delta_m$ satisfies a stopping condition.

\dk{This subproblem should be solved by special gradient-based subroutine. It is written in details in \cite{tripuraneni2018stochastic}.
Since only the gradient is used to solve the subproblem, we need to compute only a Hessian-vector product $B^k [s]$ but not a full Hessian $B^k$. If our function can be represented by a computational tree, then we can use autogradient techniques and compute Hessian-vector products as fast as we compute gradients up to a small constant.}

How many Hessians should we take? By concentration inequalities it is possible to show that we need
\begin{equation}
      |S_2|\dk{=r_2}= O\left(\varepsilon^{-1}\right).
\end{equation}

So in total, the method converges with $O\left(\varepsilon^{-3/2}\right)$ iterations and $O\left(\varepsilon^{-5/2}\right)$ Hessian calculations of the function. 

In paper \cite{arjevani2020second} this approach is improved by using special variance reduction technique. Authors get method that needs only $O(\varepsilon^{-3})$ gradients and Hessian-vector products for finding second-order stationary point. Also, in this article authors prove lower bounds for higher-order stochastic problems.

\dk{What is the main advantage of such methods? We calculate fewer Hessians than in the full CR version and also do it in parallel if we have many cores for computing. The simplicity of the algorithms, both at fast rates and when escaping from saddle-points, leads us to very good optimization methods for non-convex stochastic problems.}

Next we go to offline version that works with sum of functions \eqref{eq:f_finite_sum}.
\begin{equation}
\label{eq:f_finite_sum_tensor}
    f(x) = \frac{1}{m}\sum\limits_{i=1}^m f_i(x), 
\end{equation}
\dk{where $f_i(x)$ has Lipschitz continuous Hessian.}
In this regime we have $m$ functions and hence classic CR needs to compute $O(m\varepsilon^{-3/2})$ Hessians. To reduce it in papers \cite{kohler2017sub,xu2020newton} authors used subsampled gradient and subsampled Hessian, which achieve $\tilde{O}(m \varepsilon^{-3/2} \wedge \varepsilon^{-7/2})$ gradient complexity and $\tilde{O}(m \varepsilon^{-3/2} \wedge \varepsilon^{-5/2})$ Hessian complexity similarly to the previous section. Next appears many articles with different stochastic variance-reduced cubic(SVRC) methods. To collect this results in one place we add a table \eduard{(see Table~\ref{tab:second_vr_methods})} with the convergence rates, where $a \wedge b = \eduard{\min\{a, b\}}$.

\begin{table}[h!]
    \centering
    \begin{tabular}{|c|c|c|c|}
        \hline
        Method & Gradient & Hessian \\
        \hline
        CR \cite{nesterov2006cubic}  &$O\left(m \cdot\varepsilon^{-3/2} \right)$ & $O\left( m \cdot\varepsilon^{-3/2}  \right)$\\
        \hline
        SCR \cite{kohler2017sub,xu2020newton}  &$\tilde{O}\left(m \cdot\varepsilon^{-3/2} \wedge \varepsilon^{-7/2} \right)$ & $\tilde{O}\left( m \cdot\varepsilon^{-3/2} \wedge \varepsilon^{-5/2} \right)$\\
        \hline
        SVRC1 \cite{zhou2019stochastic}  &$\tilde{O}\left(m^{4/5}\cdot\varepsilon^{-3/2}\right)$ & $\tilde{O}\left(m^{4/5}\cdot\varepsilon^{-3/2}\right)$\\
        \hline
        SVRC2 \cite{wang2019stochastic, zhou2019stochastic_jmlr}  &$\tilde{O}\left(m \cdot\varepsilon^{-3/2}\right)$ & $\tilde{O}\left(m^{2/3}\cdot\varepsilon^{-3/2}\right)$\\
        \hline
        SVRC3 \cite{zhang2018adaptive}  &$\tilde{O}\left(m \cdot \varepsilon^{-3/2} \wedge m^{2/3}\cdot\varepsilon^{-5/2}\right)$ & $\tilde{O}\left(m^{2/3}\cdot\varepsilon^{-3/2}\right)$\\
        \hline
        STR \cite{shen2019stochastic}  &$\tilde{O}\left(m \cdot \varepsilon^{-3/2} \wedge m^{1/2}\cdot\varepsilon^{-2}\right)$ & $\tilde{O}\left(m^{1/2}\cdot\varepsilon^{-3/2}\wedge \varepsilon^{-2}\right)$\\
        \hline
        SRVRC \cite{zhou2020stochastic}  &$\tilde{O}\left(m \cdot \varepsilon^{-3/2} \wedge m^{1/2}\cdot\varepsilon^{-2}\wedge \varepsilon^{-3}\right)$ & $\tilde{O}\left(m^{1/2}\cdot\varepsilon^{-3/2}\wedge \varepsilon^{-2}\right)$\\
        \hline
         Lower bound \cite{emmenegger2021oracle}  &$\Omega\left(m^{1/4} \cdot \varepsilon^{-3/2}\right)$ & $\Omega\left(m^{1/4} \cdot \varepsilon^{-3/2}\right)$\\
        \hline
    \end{tabular}
    \caption{An Overview of the number of computations of gradients and Hessians of functions in \eqref{eq:f_finite_sum_tensor}.}
    \label{tab:second_vr_methods}
\end{table}

As a result, we have a method that not only works efficiently with the big sum by utilizing stochastic nature, but also employs Hessian information to escape saddles more effectively and arrive at to better local minimum. This statement is supported by the experiments described in \cite{xu2020second, osawa2018second,martens2010deep,park2020combining}. The authors of these papers experiment with various second-order methods and show how they compete with first-order methods without any second-order information in practice. These papers' main conclusions are that second-order methods find deeper local minima and avoid saddle-points. They are more robust when hyperparameters are used. Subsampling speeds up computations and allows for the parallelization of such methods. As a result, second-order methods may be competitive with first-order methods in practice.

%-----------------------------

\subsection{Tensor Methods}

Next, we present high-order or tensor methods for finding local minima of a highly smooth and non-convex objective function. High-order derivatives better describe functions and enable you to use curvature to improve convergence. 

First, we lay out some standard assumptions about the smoothness of the function $f$. 
In the following, we will denote the directional derivative of the function $f$ at $x$ along the directions $h^j \in \R^n,\, j = 1, \dots, p$ as
\begin{equation*}
\nabla^p f(x)[h^1, \dots, h^p].
\end{equation*}
For instance, $\nabla f(x)[h] = \nabla f(x)^\top h$ and $\nabla^2 f(x)[h]^2 = h^\top \nabla^2 f(x) h$.

The functions $f_i$ for each $p = 0, \dots, 3$ has $L_p$-Lipschitz-continuous derivatives,
\begin{equation*}
\| \nabla^p f_i (x) - \nabla^p f_i (y) \|_2 \leq L_p \| x - y \|_2
\end{equation*}
for all $x, y \in \R^\dk{n}$.

\dk{From this inequality we get next tensor method for $p=3$,
\begin{equation}
    \label{eq_tensor_newton}
x^{k+1}_{Tensor} = x^{k} + \argmin_{s\in \R^n} \left[ \nabla f(x^{k})[s] + \frac{1}{2} \nabla^2 f(x^{k}) [s]^2 + \frac{1}{6} \nabla^3 f(x^{k}) [s]^3 + \frac{H}{4!} \|s\|_2^4 \right],
\end{equation}
}
\dk{In papers \cite{birgin2017worst,carmon2019lowerI,carmon2019lowerII} it was proved that tensor $p$-order method with Taylor approximation is optimal, match lower bounds, and  converges with the rate $O(\varepsilon^{-(p+1)/p})$ for non-convex problems, hence for the third-order methods we get the rate $O(\varepsilon^{-4/3})$ instead of $O(\varepsilon^{-3/2})$ for the second-order methods. So, we get that third-order methods are faster than second-order methods in terms of iterations.}

\dk{Another crucial motivation is that the second-order method could get stuck at the so-called degenerate saddle point, where the Hessian matrix has nonnegative eigenvalues with some eigenvalues equal to 0 \cite{anandkumar2016efficient}.}

In paper \cite{zhu2020adaptive} it is shown how gradient descent and cubic regularization method stuck in such points for even small problems, like $f(x,y) = x^3 - 3x y^2$ in degenerate saddle \dk{point} $(0,0)$. So, we should use third-order information to escape \dk{them}. 

This lead us to the third-order critical point. We define next critically measures
\begin{align*}
    \chi_{f,1}(x_k) &= \|\nabla f(x_k)\|_2,\\
    \chi_{f,2}(x_k) &= \max\left\lbrace0, - \lambda_{\min}\left(\nabla^2 f(x_k)\right)\right\rbrace,\\
    \chi_{f,3}(x_k) &= \max\limits_{y\in Z_{k+1}} \left|\nabla^3 f(x_k)[y]^3 \right|,
\end{align*}
where $Z_{k+1}$ is the kernel of $\nabla^2 f(x_k)$. 
Then, we define $x^{\ast}$ a $(\eps_1,\eps_2,\eps_3)$-third-order critical point if 
\begin{align*}
    \chi_{f,1}(x_k) \leq \eps_1, \quad
    \chi_{f,2}(x_k) \leq \eps_2, \quad
    \chi_{f,3}(x_k) \leq \eps_3,
\end{align*}

\dk{Third-order method converges to a $(\eps_1,\eps_2,\eps_3)$-third-order critical point with the rate $O\left(\max\left(\eps_1^{-4/3},\eps_2^{-2},\eps_3^{-4}\right)\right)$.}

\dk{But the calculation of the third-order derivative would be very computationally expensive. This problem leads us to stochastic tensor methods.}   
The main idea of the stochastic method that by different concentration inequalities we can compute much fewer Hessians and third-order derivatives for sum type problems, than gradients. Correct proportions is written in \eqref{eq_tensor_samplesize}. For example, if we have $200000$ functions in sum, we may compute full gradient, only $10000$ Hessians and $100$ third-order derivatives and get the same speed as for full Hessian and full third-order derivatives.

\dk{ In paper by \cite{lucchi2019stochastic} introduce such method that work with batch tensors and converges as fast as for full-batch methods. The optimization algorithm we consider is detailed in Algorithm~\ref{alg:stm}.
This algorithm uses sub-sampled derivatives instead of exact quantities and its implementation relies on tensor-vector products only. The proposed approach is shown to find an $(\eps_1,\eps_2,\eps_3)$-third-order critical point in at most $O\left(\max\left(\eps_1^{-4/3}, \eps_2^{-2}, \eps_3^{-4}\right)\right)$ iterations, thereby matching the rate of deterministic approaches. }

We construct an inexact Taylor approximation model and add a fourth-order regularization defined as:
\begin{align}
\phi_k(s) &= f(x^k) + \dk{g^k [s]+ \frac12 B^k [s]^2+ \frac{1}{6} T^k [s]^3},\\
\dk{\psi}_k(s) &= \phi_k(s) + \frac{\dk{H}_k}{4} \|s\|^4_2,
\label{eq:model_tensor}
\end{align}
where $\dk{g^k}, \dk{B^k}$ and $\dk{T^k}$ approximate the derivatives $\nabla f(x^k), \nabla^2 f(x^k)$ and $\nabla^3 f(x^k)$ through sampling as follows. Three sample sets $S^g, S^b$ and $S^t$ are drawn and the derivatives are then estimated as
\begin{align*}
\dk{g^k} &= \frac{1}{|S^g|} \sum_{i \in S^g} \nabla f_i(x^k),
\dk{B^k} = \frac{1}{|S^b|} \sum_{i \in S^b} \nabla^2 f_i(x^k), \nonumber \\
\dk{T^k} &= \frac{1}{|S^t|} \sum_{i \in S^t} \nabla^3 f_i(x^k).
\end{align*}

\dk{It is worth mentioning that the implementation of the algorithm does not require the computation of the Hessian or the third-order tensor, both of which would demand significant computational resources, but rather directly computes Tensor-vector products with a complexity of order $O(n)$.}

We will make use of the following condition in order to reach an $\varepsilon$-critical point \dk{(where $\varepsilon=\eps_1$)}. For a given $\varepsilon$ accuracy, one can choose the size of the sample sets $S^g, S^b, S^t$ for sufficiently small $\kappa_g, \kappa_b, \kappa_t > 0$ such that:
\begin{align}
\label{eq:sampling_g_tensor}
\| \dk{g^k} - \nabla f(x^k) \|_2 &\leq \kappa_g \varepsilon, \\ 
\label{eq:sampling_b_tensor}
\| (\dk{B^k} - \nabla^2 f(x^k)) s\|_2 &\leq \kappa_b \varepsilon^{2/3} \| s\|_2, \;\; \forall s\in \R^\dk{n}, \\
\label{eq:sampling_t_tensor}
\| \dk{T^k}[s]^2 - \nabla^3 f(x^k) [s]^2 \|_2 &\leq \kappa_t \varepsilon^{1/3} \| s\|_2^2, \;\; \forall s\in \R^\dk{n}.
\end{align}

\dk{In practice, we can choose the size of the sample sets $S^g, S^b$ and $S^t$ as follows}
\begin{align}
\dk{r}_g = \tilde{O} \left( \frac{L_0^2}{\kappa_g^2 \varepsilon^{2}} \right), \quad
\dk{r}_b = \tilde{O} \left( \frac{L_1^2}{\kappa_b^2 \varepsilon^{4/3}} \right),  \quad
\dk{r}_t = \tilde{O} \left( \frac{L_2^2}{\kappa_t^2 \varepsilon^{2/3}} \right), 
\label{eq_tensor_samplesize}
\end{align}
where $\tilde{O}$ hides poly-logarithmic factors and a polynomial dependency to $\dk{n}$. \dk{We can see that due to the stochastic nature of the data and tensor concentration inequalities, we can use far fewer computations while still achieving the same convergence speed as a full-batch method.}

%%%%%%%%%%%%%%%%%%%%%%%%% ALGORITHM %%%%%%%%%%%%%%%%%%%%%%%%%%%%%%%%%%%%%%%%%%%%
 \begin{algorithm*}[tb]
   \caption{Stochastic Tensor Method (STM)}
   \label{alg:stm}
\begin{algorithmic}[1]
   \STATE {\bfseries Input:} \\ 
   %$\quad$ training examples $s= (x_1,x_2,\ldots,x_n)$, $x_i \sim P$\\
   $\quad$ Starting point $x^0 \in \R^\dk{n}$ (e.g~$x^0 = {\bf 0}$) \\
   $\quad 0 < \gamma_1 < 1 < \gamma_2 < \gamma_3, 1>\eta_2>\eta_1>0$, and $\dk{H}_0>0, \dk{H}_{min}>0$
   \FOR{$k=0,1,\dots,\text{until convergence}$}
   \STATE Sample gradient $\dk{g^k}$, Hessian $\dk{B^k}$ and $\dk{T^k}$ such that Eq.~\eqref{eq:sampling_g_tensor}, Eq.~\eqref{eq:sampling_b_tensor} \& Eq.~\eqref{eq:sampling_t_tensor} hold.
   \STATE Obtain $s^k$ by solving $\psi_k(s^k)$ (Eq.~\eqref{eq:model_tensor}).
   \STATE Compute $f(x^k + s^k)$ and 
	\begin{equation}
	\rho_k=\dfrac{f(x^k)-f(x^k + s^k)}{f(x^k) - \phi_k(s^k)}.
	\end{equation}
	\STATE Set
	\begin{equation}
	x^{k+1} = \begin{cases}
	x^k + s^k & \text{ if } \rho_k \geq \eta_1\\
	x^k & \text{ otherwise.}
	\end{cases}
	\end{equation}
	\STATE Set
	\begin{equation*} \label{eq:sigma_update}
	\dk{H}_{k+1}= \begin{cases}
	[\max\{\dk{H}_{min},\gamma_1 \dk{H}_k \}, \dk{H}_k] & \text{ if } \rho_k>\eta_2 \text{ (very successful iteration)}\\
	[\dk{H}_k, \gamma_2 \dk{H}_k] & \text{ if } \eta_2\geq\rho_k\geq \eta_1 \text{ (successful iteration)}\\
	[\gamma_2 \dk{H}_k, \gamma_3 \dk{H}_k] & \text{ otherwise}\text{ (unsuccessful iteration)}.
	\end{cases}
	\end{equation*}
   \ENDFOR
\end{algorithmic}
\end{algorithm*}
%%%%%%%%%%%%%%%%%%%%%%%%% ALGORITHM %%%%%%%%%%%%%%%%%%%%%%%%%%%%%%%%%%%%%%%%%%%%

\dk{As shown in \cite{emmenegger2021oracle}, the lower bounds for sum type problem are still rather far from upper bound even for the second-order methods. Hence, further research in this area may lead to new methods for sum-type problems by using variance reduction techniques. Another branch of possible research is a combination of tensor methods with first or second-order methods.}

\section{Zeroth-Order Methods}

Gradient free or zeroth-order optimization methods, which use only function values, are becoming increasingly important in machine learning problems, especially in reinforcement learning \cite{malik2020derivativefree}, black-box adversarial
attacks on deep neural networks \cite{papernot2017practical} and other
problems with structure making gradients difficult or infeasible to obtain. 

While there is a class of methods that does not have any connection to the gradient, for example, random search algorithms \cite{Schumer1968} (which are one of the first methods of zeroth-order optimization, beside grid search), the Nelder--Mead algorithm \cite{nelder1965simplex}, the model-based methods (see Chapters 2-6 and 10-11 in \cite{conn2009introduction}) or the recent stochastic three points (STP) method \cite{bergou2019stochastic} and its momentum variant STMP \cite{gorbunov2020smtp} most zeroth-order optimization methods use gradient estimations, such as $g(x) = \sum_{i=1}^{n}\kesha{\ls\nicefrac{\lp f(x+\mu e_i) - f(x)\rp}{\mu}\rs} e_i$ (where $e_i$ are columns of $n\times n$ identity matrix $I_n$, $i\in\{1,\ldots, n\}$), then for good enough functions ($f\in C^{1,1}_{L}$ i.e. continuously differentiable with Lipschitz-continuous gradient) it can be shown, for example, that $\|g(x) - \nabla f(x)\|_2 \leqslant \kesha{\mu L}\sqrt{n}$. One then can consider some first-order optimization scheme, replace actual gradients with their estimations, and use bounds like this to return to gradients from estimations in proofs, obtaining the results for the zeroth-order case relatively easy.

While such deterministic zeroth-order schemes (like the $GD$ with gradient estimation of the same form as above) often suffer from the problem dimensionality because of the number of oracle calls needed to reconstruct the gradient ($n$ for the estimation mentioned above, see also \cite{berahas2020theoretical} for other examples), in a randomized approach one can use two- or one- point schemes of gradient approximation which makes every iteration simpler, sometimes leading to better results in terms of oracle calls \cite{liu2018zerothorder}. Another benefit of the stochastic approach is that such methods often have good theoretical properties, for example, the Gaussian smoothing approach \cite{nesterov2017random} that gives a smoothed version of the initial function, for which the convergence of stochastic zeroth-order algorithm can be easily proved, which can be later used to show the convergence of the algorithm for the initial function. And there are setups (for example online learning \cite{bubeck2011introduction}) where one is limited to use only several (or even one) oracle queries thus being unable to construct the full gradient approximation, so the stochastic approach becomes the only option. 

We begin with the formalization of these zeroth-order randomized schemes - we have a problem with the form
\begin{align*}
    \min\limits_{x\in Q \subseteq \mathbb{R}^n}f(x)
\end{align*}
then stochastic zeroth-order methods generate $\{x^k\}$ s.t.
\begin{align}\label{stochastic_zero_order_base_procedure}
    x^{k+1} = A\left( \hat{f}, X, P, \{x^i\}_{i=0}^{k}, \{u^i\}_{i=0}^{k}\right)
\end{align}
so the procedure $A$ gives us $x^{k+1}$ based on function values (obtained via oracle $\hat{f}$), history of $\{x^k\}$, random vectors  $\{u^k\}$, and parameters $P$ such as dimension $n$ of $X$, $L_\nu$ and $\nu$ -- H\"older parameters, etc. Function $\hat{f}$ is not necessarily equal to $f$, we can, for example, use $\hat{f}(x) = f(x) + \eps(x)$ where $|\eps(x)|\ll |f(x)|$, or $\hat{f}(x, u) = f(x) + \eps(x, u)$ s.t.  $\Exp_u[\hat{f}(x, u)] = f(x)$.

In the subsections, we will discuss the characteristics of several zeroth-order gradient estimations and then the zeroth-order methods for sum-minimization type problems in a non-convex setup. Other information on gradient-free optimization (such as structured objectives) can be found in the recent survey \cite{Larson_2019}.

\subsection{Random Directions Gradient Estimations}

Let us start with the methods following the standard zeroth-order scheme of using gradient approximation to benefit from the analysis of first-order methods. In this section all methods have a form similar to the classic gradient descent 
\begin{align*}
    x^{k+1} = x^k - h_k g(x^k, u^k)
\end{align*}
with only difference that instead of the true gradient we use the gradient approximation $g(x, u)$. One way to build such gradient approximations is to use random directions to compute finite differences in the form
\begin{align*}
    g(x^k, u^k) :=\frac{\hat{f}(x^k+\mu u^k) - \hat{f}(x^k)}{\mu}\cdot u^k
\end{align*}

It makes sense to use centrally symmetric distributions for $u^k$, for example uniformly distributed over the unit Euclidean sphere $S^{n-1}=\{x\in\mathbb{R}^n:\|x\|_2=1\}$ (see \cite{flaxman2005online,gorbunov2018accelerated,dvurechensky2021accelerated}), or $u^k\sim\mathcal{N}(0,I_n)$ --- so-called Gaussian smoothing introduced in \cite{nesterov2017random}. In this article, the authors proved Gaussian approximation 
\begin{align*}
    f_{\mu}(x) = \frac{1}{\kappa}\int\limits_{\mathbb{R}^n}f(x+\mu u) e^{-\tfrac{1}{2}\|u\|_2^2}du
\end{align*}
(there \kesha{$\kappa = \int_{E}e^{-\nicefrac{\|u\|_2^2}{2}}du = (2\pi)^{\nicefrac{n}{2}}$}) to have several good properties, such as convexity preservation (if $f$ is convex then $f_\mu$ is convex too), differentiability, and if $f\in C^{0,0}_{L_0}$ or $f\in C^{1,1}_{L_1}$ (i.e. Lipschitz-continuous function with constant $L_0$ or function with Lipschitz-continuous gradient with $L_1$ respectively) then the same holds for $f_{\mu}$ with $L_0(f_\mu)\leqslant L_0(f)$ and  $L_1(f_\mu)\leqslant L_1(f)$ respectively. It can be also shown that $|f_{\mu}(x)-f(x)|\leqslant \mu L_{0}\kesha{\sqrt{n}}$ for the case of $f\in C^{0,0}_{L_0}$.

While in that paper the authors mostly discuss the convex case, there are some results (\cite{nesterov2017random}[Section 7]) for a non-convex objective $f$ too. They consider a process $x^{k+1} = x^k - h_k g(x^k, u^k)$, with $g$ defined above, $\hat{f} = f$ and $u^k\sim\mathcal{N}(0,I_n)$, and show that for the case of $f\in C^{1,1}_{L_1}$ this process converges in the sense of $\Exp_{U}\|\nabla f_{\mu}(x)\|_{2}$ (where $U = \{u^k\}_{k=0}^{N-1}$):
\begin{align*}
    \frac{1}{N}\sum\limits_{k=0}^{N - 1}\Exp_{U}\ls\|\nabla f_{\mu}(x^k)\|^2_{2}\rs\leqslant 8(n+4)L_1\ls \frac{f_{\mu}(x^0) - f^{\ast}}{N}+\frac{3\mu^2(n+4)}{32}L_1\rs
\end{align*}
then using the fact that (\cite{nesterov2017random}[Lemma 3]) $\|\nabla f_\mu(x)  - \nabla f(x)\|_{2} \leqslant \kesha{\ls\nicefrac{\mu L_1}{2}\rs}(n+3)^{\kesha{\nicefrac{3}{2}}}$ we obtain (from $\|\nabla f(x)\|_{2}^2 \leqslant 2\|\nabla f_\mu(x)  - \nabla f(x)\|_{2}^2 + 2\|\nabla f_\mu(x)\|_{2}^2 $)
\begin{align*}
    \frac{1}{N}\sum\limits_{k=0}^{N - 1}\Exp_{U}\ls\|\nabla f(x^k)\|^2_{2}\rs\leqslant & 2\frac{\mu^2 L_1^2}{4}(n+3)^{3} \\ 
    & + 16(n+4)L_1\ls \frac{f_{\mu}(x^0) - f^{\ast}}{N}+\frac{3\mu^2(n+4)}{32}L_1\rs
\end{align*}
and choosing $\mu = O\lp\kesha{\nicefrac{\eps}{\ls n^{3/2}L_1 \rs}}\rp$ we ensure $\frac{1}{N}\kesha{\sum_{k=0}^{N - 1}}\Exp_{U}\ls\|\nabla f(x^k)\|_{2}^2\rs\leqslant \eps^2$ with the upper bound for the expected number of steps $N=O\lp\kesha{\nicefrac{n}{\varepsilon^2}}\rp$.

For the case of $f\in C^{0,0}_{L_0}$
\begin{align*}
    \frac{1}{S_N}\sum\limits_{k=0}^{N - 1}h_k\Exp_{U}\ls\|\nabla f_{\mu}(x^k)\|^2_{2}\rs\leqslant \frac{1}{S_N}\ls (f_{\mu}(x^0) - f^{\ast})+\frac{1}{\mu}n^{\kesha{\nicefrac{1}{2}}}(n+4)^2 L_0^3\sum\limits_{k=0}^{N - 1}h_k^2\rs
\end{align*}
they show only that this process converges to the stationary point of $f_{\mu}(x)$ -- consider $Q$ with $\text{diam}(Q)\leqslant R$, then it can be shown that we need to make 
\begin{align*}
    N = O\lp\frac{n(n+4)^2L_0^5R}{\varepsilon^4 \delta}\rp
\end{align*}
steps to ensure that $\frac{1}{N}\kesha{\sum_{k=0}^{N - 1}}\Exp_{U}\ls\|\nabla f_{\mu}(x^k)\|^2_{2}\rs\leqslant \varepsilon^2$ keeping functional gap $|f_{\mu}(x)-f(x)|\leqslant \delta$ small. Authors also mention that with the $h_k\to 0$ and $\mu\to 0$ the convergence in the sense of $\Exp_{U}\|\nabla f(x)\|_{2}$ can be proved too.

These results can be extended \cite{shibaev2021zeroth-order} to the case of noisy $\hat{f}$ i.e. $|\hat{f}(x)-f(x)|\leqslant\delta$ for $f$ with H\"older continuous gradient ($\|\nabla f(x) - \nabla f(y) \|_2 \leqslant L_{\nu} \|x-y\|_2^{\nu}$) -- it can be shown that for a small enough noise $\delta$ these convergence rates can be preserved. More specifically, to ensure $\frac{1}{N}\kesha{\sum_{k=0}^{N - 1}}\Exp_{U}\ls\|\nabla f(x^k)\|_{2}^2\rs\leqslant \eps^2$ one need to make 
\begin{align*}
    N = O\lp\frac{n^{2 + \frac{1-\nu}{\kesha{2}\nu}}}{\eps^{\frac{2}{\nu}}}\rp\text{ steps under the assumption that noise }\delta = O\lp\frac{\eps^{\frac{3+\nu}{2\nu}}}{n^{\kesha{\frac{3+7\nu}{4\nu}}}}\rp
\end{align*}
where $\nu$ is a H\"older parameter. For the convergence in the sense of smoothed function gradient norm $\frac{1}{N}\kesha{\sum_{k=0}^{N - 1}}\Exp_{U}\ls\|\nabla f_{\mu}(x^k)\|^2_{2}\rs\leqslant \varepsilon^2$ it can be shown
\begin{align*}
    N = O\lp\frac{n^{\frac{7-3\nu}{2}}}{\eps^{\frac{\kesha{2(3-\nu)}}{1+\nu}}}\rp\text{ with }\delta = O\lp\frac{\eps^{\frac{5-\nu}{1+\nu}}}{n^{\frac{13-3\nu}{4}}}\rp
\end{align*}
with functional gap $|f_{\mu}(x)-f(x)|=O\lp \kesha{\nicefrac{\eps}{\ls n^{\nicefrac{(1+\nu)}{2}}\rs}}\rp$. For the case of $\nu = 1$ (i.e. $f\in C^{1,1}_{L_1}$) these results can be improved to $N=O\lp\kesha{\nicefrac{n}{\eps^{2}}}\rp$ ($n$ times better) achieving the same rate of convergence as in previous paper \cite{nesterov2017random}.

Such noisy setup is also interesting because it can be shown \cite{Risteski2016AlgorithmsAM}, that for a non-convex function $\hat{f}(x)$ s.t. $|\hat{f}(x)-f(x)|\leqslant\eps_f$, where initial $f$ is convex and $1$-Lipschitz and $\eps_f \sim \max\lf \nicefrac{\eps^2}{\sqrt{n}}, \nicefrac{\eps}{n}\rf$ there exists an algorithm which finds a point $\tilde{x}$ s.t. $\hat{f}(\tilde{x}) \leqslant \hat{f}_{\ast} + \eps$ with complexity $Poly\left(n,\frac{1}{\eps}\right)$. The dependence $\eps_f(\eps)$ is optimal in this class of algorithms.

This Gaussian smoothing technique was later used in works \cite{ghadimi2013stochastic} (RSGF) and \cite{ghadimi2016mini-batch} (RSPGF) to obtain complexity guarantees for stochastic zeroth-order optimization. In the first one (\cite{ghadimi2013stochastic}), the unconstrained problem $Q=\mathbb{R}^{n}$ is considered, where $\hat{f}=F(x,\xi)$ s.t. $\Exp_{\xi}[F(x, \xi)] = f(x)$ and $F(\cdot,\xi)$ has a Lipschitz-continuous gradient with constant $L_1$, $\xi$ is a  random variable whose distribution $P$ is supported on $\Xi_k\subseteq R^n$. The procedure (\ref{stochastic_zero_order_base_procedure}) has a form similar to the one proposed in \cite{nesterov2017random}
\begin{align*}
    x^{k+1} = x^k - h_k G(x^k, \xi^k, u^k),~G(x^k, \xi^k, u^k):=\frac{\hat{f}(x^k+\mu u^k, \xi^k) - \hat{f}(x^k, \xi^k)}{\mu}\cdot u^k,
\end{align*}
and from $\Exp_{\xi}[F(x, \xi)] = f(x)$ it follows that
\begin{align*}
    \Exp_{\xi,u}\ls G(x,\xi,u)\rs = \nabla f_{\mu}(x).
\end{align*}

The method then chooses the $x^k$ from generated $\{x^k\}_{k=1}^{N}$ as $k=R$ where $R$ is some random variable with a probability mass function $P_{R}$ supported on $\{1,\ldots, N\}$. The main goal to introduce this random iteration count $R$ is to derive new complexity results for non-convex stochastic optimization case.

For the case of $f\in C^{1,1}_{L_1}$, smoothing parameter $\mu$, $D_f=\kesha{\sqrt{ \nicefrac{2(f(x^1)-f^{\ast})}{L}}}$,  variance $\sigma^2$  ($\Exp_{\xi}\ls\|\nabla\hat{f}(x, \xi) - \nabla f(x)\|^2_2\rs\leqslant\sigma^2$) and \kesha{the} probability mass function 
\kesha{
\begin{align*}
P_R(k) = \frac{h_k-2L(n+4)h_k^2}{\sum\limits_{i=1}^{N}(h_i-2L(n+4)h_i^2)}
\end{align*}
}
they obtain (\cite{ghadimi2013stochastic}[Theorem 3.2])
\begin{align*}
    & \frac{1}{L}\Exp\ls\|\nabla f(x^R)\|^2_2\rs\leqslant  \\
    & \leqslant \frac{ D_f^2 + 2\mu^2(n+4) \lp 1 + L(n+4)^2\sum\limits_{k=1}^{N}\lp\frac{h_k}{4} + L h_k^2\rp\rp + 2(n+4)\sigma^2 \sum\limits_{k=1}^N h_k^2}{\sum\limits_{k=1}^N\ls h_k - 2L(n+4)h_k^2\rs} 
\end{align*}
where the expectation is taken with respect to $R$, $\{\xi^k\}$. After choosing specific constant stepsizes $h_k = \kesha{\nicefrac{1}{\ls\sqrt{n+4}\rs}\cdot\min\lf\nicefrac{1}{\ls 4L\sqrt{n+4}\rs},\nicefrac{\tilde{D}}{\ls\sigma \sqrt{N}\rs}\rf}$ (note that this makes $P_R$ uniform on $\{1,\ldots, N\}$) they get (\cite{ghadimi2013stochastic}[Corollary 3.3])
\begin{align*}
    \frac{1}{L}\Exp\ls\|\nabla f(x^R)\|^2_{2}\rs\leqslant \frac{12(n+4)LD_f^2}{N} + \frac{2\sigma\sqrt{n+4}}{\sqrt{N}}\lp\tilde{D}+\frac{D_f^2}{\tilde{D}}\rp
\end{align*}
where $\tilde{D}>0$ is our estimation of $D_f$ (for example some upper bound). It can be shown that to ensure $\mathbb{P}\{\|\nabla f(x^R)\|^2_{2}\leqslant \varepsilon\}\geqslant 1-\Lambda$ (so-called $(\varepsilon,\Lambda)$-solution) the total number of calls to the oracle $\hat{f}$ can be bounded as 
\begin{align*}
    O\lp \frac{nL^2D_f^2}{\Lambda\varepsilon}+\frac{nL^2}{\Lambda^2}\lp\tilde{D}+\frac{D_f^2}{\tilde{D}}\rp^2\frac{\sigma^2}{\varepsilon^2}\rp
\end{align*}

Another method that is considered in \cite{ghadimi2013stochastic} is a two-phase method (2-RSGF), which uses the first one (RSGF) $S = \log\lp \kesha{\nicefrac{2}{\Lambda}}\rp$ times as a subroutine producing a list of candidates $\{\bar{x}^k\}_{k=1}^S$ and then the output point $\bar{x}^{\ast}$ is chosen in such a way that
\begin{align*}
    \|g(\bar{x}^{\ast})\|_{2} = \min\limits_{k=1,\ldots,S}\|g(\bar{x}^k)\|_{2},~g(\bar{x}^k)  := \frac{1}{T}\sum\limits_{i=1}^{T}G(\bar{x}^k,\xi^k,u^k)
\end{align*}
then it can be shown (\cite{ghadimi2013stochastic}[Theorem 3.4]) that $(\varepsilon,\Lambda)$-solution will be achieved after taking
\begin{align*}
    O\lp \frac{nL^2D_f^2\log\kesha{\lp\nicefrac{1}{\Lambda}\rp}}{\varepsilon}+nL^2\lp\tilde{D}+\frac{D_f^2}{\tilde{D}}\rp^2\frac{\sigma^2}{\varepsilon^2}\log\kesha{\lp\nicefrac{1}{\Lambda}\rp}+\frac{n\log^2\kesha{\lp\nicefrac{1}{\Lambda}\rp}}{\Lambda}\lp1+\frac{\sigma^2}{\varepsilon}\rp\rp
\end{align*}
calls to the $\hat{f}$ which is better than the previous one in terms of $\Lambda$.

A more general problem $\min_{x\in Q \subseteq \mathbb{R}^n}\Psi(x)=f(x)+h(x)$, where $f\in C^{1,1}_L$  and $h(x)$ is a simple convex and possibly non-smooth function is considered in \cite{ghadimi2016mini-batch}. They use a mini-batched version of gradient estimation from the previous paper \cite{ghadimi2013stochastic} and generalized projection obtaining (\cite{ghadimi2016mini-batch}[Theorem 4, Corollaries 6-7]) similar bounds for the gradient norm.

In \cite{Sener2020Learning}, the authors use symmetric gradient estimations based on uniform distribution over the sphere to build a less dimension depending method. They consider the minimization problem  $\min_{x\in\mathbb{R}^n} f(x) = \Exp_{\xi}[F(x,\xi)] = \Exp_{\xi}[\hat{f}(x,\xi)]$ (note that in this paper authors consider both $\mathbb{R}^d$ and $\mathbb{R}^n$ with $d\ll n$) where $f(x)$ is $L$-Lipschitz, and $\mu$-smooth, $|F(x,\xi)|\leqslant{\Omega}$ and $F$ variance is bounded by $V_f$. It was shown that using
\begin{align*}
    g(x^k,\xi^k,u^k) := n\frac{\hat{f}(x^k+\mu u^k, \xi^k) - \hat{f}(x^k-\mu u^k, \xi^k)}{2\mu}\cdot u^k
\end{align*}
where $u^k\sim\mathcal{U}\lp S^{n-1}\rp$ (uniform distribution on the unit sphere $S^{n-1}$) and the process $x^{k+1}=x^k-\alpha g(x^k,\xi^k,u^k)$ after $N$ steps
\begin{align*}
    \frac{1}{N}\sum\limits_{i = 1}^{N}\Exp\ls\|\nabla f(x^i)\|^2_{2}\rs =  O\lp\frac{n}{N^{\kesha{\nicefrac{1}{2}}}}+\frac{n^{\kesha{\nicefrac{2}{3}}}}{N^{\kesha{\nicefrac{1}{3}}}}\rp
\end{align*}

Now consider the case when for a given $\xi$, $F(x,\xi) = g(r(x,\theta^{\ast}),\psi^{\ast})$ (there $g(\cdot,\psi)$ and $r(\cdot,\theta)$ are parameterized function classes), where $r(\cdot, \theta^{\ast}):\mathbb{R}^n\to\mathbb{R}^d$ where $d\ll n$. To put it simply, the authors consider the case when $F(\cdot,\xi):\mathbb{R}^n\to\mathbb{R}$ while it is actually defined on an $d$-dimensional manifold $\mathcal{M}$ for all $\xi$. That means that if one knows the manifold (i.e. $\theta^{\ast}$), and $g$ and $r$ are smooth the chain rule can be applied giving $\nabla f(x) = J(x,\theta^{\ast})\nabla_{r}g(r,\psi)$ (where $J(x,\theta^{\ast}) = \kesha{\nicefrac{\partial r(x,\theta^{\ast})}{\partial x}}$) leading to
\begin{align*}
    g(x^k,\xi^k,u^k)) := d\frac{\hat{f}(x^k+\mu J_q u^k,\xi^k)) - \hat{f}(x^k-\mu J_q u^k,\xi^k))}{2\mu}\cdot u^k
\end{align*}
where $J_q$ is the orthonormalized $J(x^k,\kesha{\theta^{\ast}})$ and $u^k\sim\mathcal{U}\lp S^{d-1}\rp$, and this gives 
\begin{align*}
    \frac{1}{N}\sum\limits_{i = 1}^{N}\Exp\ls\|\nabla f(x^i)\|^2_{2}\rs =  O\lp\frac{d}{N^{\kesha{\nicefrac{1}{2}}}}+\frac{d^{\kesha{\nicefrac{2}{3}}}}{N^{\kesha{\nicefrac{1}{3}}}}\rp
\end{align*}
which is much better than the previous one (because $d\ll n$). However, this is impractical due to the fact that it requires the knowledge of $\theta^{\ast}$. Authors mix two previous estimations and estimate $\theta$ and $\psi$ on every step, obtaining the method that (\cite{Sener2020Learning}[Theorem 1]) after $N$ steps ensures
\begin{align*}
    \frac{1}{N}\sum\limits_{i = 1}^{N}\Exp\ls\|\nabla f(x^i)\|^2_{2}\rs = O\lp\frac{n^{\kesha{\nicefrac{1}{2}}}}{N}+\frac{n^{\kesha{\nicefrac{1}{2}}} + d + dn^{\kesha{\nicefrac{1}{2}}}}{N^{\kesha{\nicefrac{1}{2}}}} + \frac{d^{\kesha{\nicefrac{2}{3}}}+ n^{\kesha{\nicefrac{1}{2}}}d^{\kesha{\nicefrac{2}{3}}}}{N^{\kesha{\nicefrac{1}{3}}}}\rp
\end{align*}
which is better than the initial bound for $d\leqslant n^{\kesha{\nicefrac{1}{2}}}$.

While such gradient estimates based on random directions are common it can be shown that in terms of the number of samples required to the approximate gradient to ensure norm condition (or at least ensure it with some probability) random directions based methods lose to standard finite differences \cite{berahas2020theoretical, berahas2019global, berahas2019linear}. In these papers, authors consider an unconstrained optimization problem $\min_{x\in\mathbb{R}^n}f(x)$ where $\hat{f}(x) = f(x) + \eps(x)$ is computable, the noise $\eps$ is bounded uniformly: $|\eps(x)|\leqslant \eps_f$ and $f(x)\in C^{1,1}_L$ or $f(x)\in C^{2,2}_M$ (i.e. twice continuously differentiable function with $M$-Lipschitz continuous Hessian) . 

The main idea in \cite{berahas2020theoretical} is to compare the number of calls \kesha{r (essentially a batch size)} to the oracle $\hat{f}(x)$ that will be enough to ensure norm condition
\begin{align}\label{norm_condition}
    \|g(x)-\nabla f(x)\|_2\leqslant \theta \|\nabla f(x)\|_2,~\theta\in[0,1)
\end{align}
for zeroth-order gradient estimation $g(x)$. This condition simplifies the transition from gradient estimations to gradient when proving the convergence of algorithms. One of its implications is that $g(x)$ is a descent direction for the function $\phi$. In \cite{berahas2019global} the line-search method that uses such gradient approximations, ensuring the norm condition, is shown to converge. 

They consider several methods of gradient estimation, deterministic (Forward and Central Finite Differences ($FFD$ and $CFD$) and Linear Interpolation ($LI$) as generalization) and stochastic (Gaussian Smoothed Gradients ($GSG$ and its centered version $cGSG$) and Sphere Smoothed Gradients ($BSG$ and $cBSG$)), for the latter authors obtain the number of calls needed to ensure the norm condition with probability $1-\delta$.

\begin{table}[h]
    \centering
    \begin{tabular}{|c|c|c|c|}
        \hline
        Name & Gradient estimation $g(x)$ form & Number of calls \kesha{$r$} & $\|\nabla f(x)\|_2$\\
        \hline
        $FFD$ & $\sum\limits_{i = 1}^{n}\frac{\hat{f}(x+\mu e_i)-\hat{f}(x)}{\mu}e_i$ & $n$ & $\frac{2\sqrt{nL\eps_f}}{\theta}$\\ 
        \hline
        $CFD$ & $\sum\limits_{i = 1}^{n}\frac{\hat{f}(x+\mu e_i)-\hat{f}(x-\mu e_i)}{2\mu}e_i$  & $n$ & $\frac{2\sqrt{n}\sqrt[3]{M\eps_f^2}}{\sqrt[3]{6}\theta}$\\ 
        \hline
        $LI$ & $\sum\limits_{i = 1}^{n}\frac{\hat{f}(x+\mu u^i)-\hat{f}(x)}{\mu}u^i$, $u^i=[Q]_{i}$,   & $n$ & $\frac{2\|Q^{-1}\|\sqrt{nL\eps_f}}{\theta}$ \\ 
        \hline
        $GSG$ & $\frac{1}{\kesha{r}}\sum\limits_{i = 1}^{\kesha{r}}\frac{\hat{f}(x+\mu u^i)-\hat{f}(x)}{\mu}u^i$, $u^i\sim\mathcal{N}(0,I_n)$  & $\frac{12n}{\delta\theta^2} + \frac{n+20}{16\delta}$ & $\frac{6n\sqrt{L\eps_f}}{\theta}$\\ 
        \hline
        $cGSG$ & $\frac{1}{\kesha{r}}\sum\limits_{i = 1}^{\kesha{r}}\frac{\hat{f}(x+\mu u^i)-\hat{f}(x - \mu u^i)}{2 \mu}u^i$, $u^i\sim\mathcal{N}(0,I_n)$  & $\frac{12n}{\delta\theta^2} + \frac{n+30}{144\delta}$& $\frac{12\sqrt[3]{n^{\kesha{\nicefrac{7}{2}}}M\eps_f^2}}{\theta}$\\
        \hline
        $BSG$ & $\frac{n}{\kesha{r}}\sum\limits_{i = 1}^{\kesha{r}}\frac{\hat{f}(x+\mu u^i)-\hat{f}(x)}{\mu}u^i$, $u^i\sim\mathcal{U}\lp S^{n-1}\rp$  & $\ls\frac{8n}{\theta^2} + \frac{8n}{3\theta} + \frac{11n+104}{24}\rs\log\frac{n+1}{\delta}$& $\frac{4n\sqrt{L\eps_f}}{\theta}$ \\ 
        \hline
        $cBSG$ & $\frac{n}{\kesha{r}}\sum\limits_{i = 1}^{\kesha{r}}\frac{\hat{f}(x+\mu u^i)-\hat{f}(x - \mu u^i)}{2 \mu}u^i$, $u^i\sim\mathcal{U}\lp S^{n-1}\rp$  & $\ls\frac{8n}{\theta^2} + \frac{8n}{3\theta} + \frac{9n+192}{27}\rs\log\frac{n+1}{\delta}$& $\frac{4\sqrt[3]{n^{\kesha{\nicefrac{7}{2}}}M\eps_f^2}}{\theta}$\\ 
        \hline
    \end{tabular}
    \caption{Bounds on number of $\hat{f}$ calls \kesha{$r$}, and $\|\nabla f(x)\|_2$ that ensure the norm  condition $\|g(x)-\nabla f(x)\|_2\leqslant \theta \|\nabla f(x)\|_2$. For the $GSG$, $cGSG$, $BSG$ and $cBSG$ these are the results with probability $1-\delta$. The gradient norm bound (last column) essentially means that for a noisy oracle $\hat{f}$ we can ensure norm condition only for big enough gradients. The $LI$ method is basically $FFD$ with directions given as columns of the nonsingular  matrix $Q$. When $Q$ is orthonormal the $g(x)$ takes a form from the table. }
    \label{tab:zo_gradient_estimations_comparison}
\end{table}

Let us take a look at two of these methods: $FFD$ and $GSG$. For the first one, the gradient estimation takes the form
\begin{align*}
    g(x) := \sum\limits_{i = 1}^{n}\frac{\hat{f}(x+\mu e_i)-\hat{f}(x)}{\mu}e_i
\end{align*}
where $e_i$ are the columns of $I_n$. It can be shown that for such $g(x)$ the following holds
\begin{align*}
    \|g(x)-\nabla f(x)\|_{2}\leqslant\frac{\kesha{\mu L\sqrt{n}} }{2} + \frac{2\kesha{\eps_f\sqrt{n}}}{\mu}.
\end{align*}
If there was no noise ($\eps_f=0$) we could make this approximation as close to the gradient as we want, so we would be able to ensure the norm condition in $n$ calls to the $\hat{f}$. This is also true for a small enough noise (for example even from this inequality we can take $\eps_f = \kesha{\nicefrac{L\mu^2}{4}}$ obtaining $\|g(x)-\nabla f(x)\|_{2}\leqslant\kesha{\mu L\sqrt{n}}$). Authors provide such noise bound in form of lower bound on $\|\nabla f(x)\|_{2}$ for which the norm condition can still be ensured
\begin{align*}
     2\sqrt{\frac{\eps_f}{L}}\leqslant\mu\leqslant\frac{\theta\|\nabla f (x)\|_{2}}{\sqrt{n} L}\Rightarrow\frac{2\sqrt{n L\eps_f}}{\theta}\leqslant\|\nabla f (x)\|_{2}
\end{align*}
In other words, that means that we can converge to the neighborhood where $\|\nabla f(x)\|_2\approx \kesha{\nicefrac{2\sqrt{n L\eps_f}}{\theta}}$. 

For the $GSG$ they consider the mini-batched version of Gaussian smoothing from \cite{nesterov2017random}
\begin{align*}
    g(x,\{u^i\}):=\frac{1}{\kesha{r}}\sum\limits_{i = 1}^{\kesha{r}}\frac{\hat{f}(x+\mu u^i)-\hat{f}(x)}{\mu}u^i,~u^i\sim\mathcal{N}(0,I_n)
\end{align*}
and prove that the norm condition will be ensured with probability $1-\delta$ after 
\begin{align*}
    \kesha{r} \geqslant \frac{3n}{\delta \theta^2}\frac{n}{(\sqrt{n} - 1)^2} + \frac{(n+4)}{16\delta} + \frac{1}{\delta} = \Omega\lp\frac{3n}{\theta^2 \delta }\rp
\end{align*}
calls, which is while linear on $n$ is still worse than the plain $n$ in $FFD$, because of $\delta$, and additional constants. However, this is a sufficient number of calls, not a necessary, so authors derive the lower bound for \kesha{$r$} (\cite{berahas2020theoretical}[Section 2.3.1]) 
\begin{align*}
    \kesha{r} \geqslant \frac{1-\sqrt{\delta}}{\theta^2}(n+1)
\end{align*}
necessary to have probability $\mathbb{P}(\|g(x)-\nabla f(x)\|_2\leqslant \theta \|\nabla f(x)\|_2) > 1 -\delta $. In their numerical experiments they show that to ensure the norm condition with $\theta < \kesha{\nicefrac{1}{2}}$ with probability of at least $\kesha{\nicefrac{1}{2}}$ more than $n$ oracle calls are needed, so this lower bound is weak.

The sufficient lower bound can be improved using smoothing on a sphere for which they obtain $\Omega\lp\kesha{\nicefrac{n}{\theta^2}\cdot\log\ls\nicefrac{(n+1)}{\delta}\rs}\rp$, yet it is still worse than deterministic variants, and in practice its behavior is very similar to the Gaussian directions based approach.

There are also results for the case of $f(x)\in C^{2,2}_M$ (centered versions of the estimations), they can be found in Table \ref{tab:zo_gradient_estimations_comparison}.

\subsection{Variance-Reduced Zeroth-Order Methods}

One special case of the $\min f(x)$ problem is the finite sum minimization which was considered in previous sections for the first-order methods. These problems in zeroth-order setup arise in reinforcement learning \cite{fazel2019global} (there as a minimization of a long-term cost which is essentially a sum of functions) and non-stationary online optimization problems \cite{zhang2020boosting}. 

Let us start with the ZO-SVRG from \cite{liu2018zerothorder} -- a zeroth-order version of SVRG from \cite{johnson2013accelerating}. 

There a non-convex finite-sum problem of the form
\begin{align*}
    \min\limits_{x\in\mathbb{R}^n} f(x) = \frac{1}{m}\sum\limits_{i=1}^m f_i(x)
\end{align*}
where $f_i\in C^{1,1}_{L}$ i.e. $\|\nabla f_i(x) - \nabla f_i(y)\|_{2} \leqslant L\|x - y\|_{2}$ for any $x,y\in\mathbb{R}^n$ and $i\in\{1,\ldots,m\}$ is considered. Authors use the standard assumption that the variance of stochastic gradients is bounded 
\begin{align*}
    \frac{1}{m}\sum\limits_{i=1}^m\|\nabla f_i(x) - \nabla f_i(y)\|_{2}^2 \leqslant \sigma^2
\end{align*}
and consider several different gradient estimates: two based on random directions on a unit sphere (in notation of \cite{berahas2019global} these are $BSG$ with $N=1$ and $N=q$ (see Table \ref{tab:zo_gradient_estimations_comparison}), called RandGradEst and Avg-RandGradEst respectively), and one deterministic coordinate estimation (variant of $CFD$ from Table \ref{tab:zo_gradient_estimations_comparison} with possibly different $\mu_j$ for each direction $e_j$ called CoordGradEst)
\begin{align*}
    \text{RandGradEst}:&~\hat{\nabla} f_{i}(x) = \frac{n}{\mu}[f_{i}(x+\mu u^i) - f_i(x)]u^i,\\
    \text{Avg-RandGradEst}:&~\hat{\nabla} f_{i}(x) = \frac{n}{\mu q}\sum\limits_{j=1}^{q}[f_{i}(x+\mu u^{i,j}) - f_i(x)]u^{i,j},\\
    \text{CoordGradEst}:&~\hat{\nabla} f_{i}(x) = \frac{1}{2\mu}\sum\limits_{j=1}^{n}[f_{i}(x+\mu_j e_j) - f_i(x-\mu_j e_j)]e_j
\end{align*}
there $i\in\{1,\ldots,m\}$, $\mu > 0$, and $\{e_j\}_{j=1}^{n}$ are standard basis vectors (columns of $I_n$). 

\begin{algorithm}[ht]
    \caption{ZO-SVRG \cite{liu2018zerothorder}}
    \label{alg:ZO-SVRG}
    \begin{algorithmic}
        \REQUIRE stepsizes $\{h^k_s\}$, epoch length $T$, starting point $x^0\in\R^n$, batch size $r \ge 1$, smoothing parameter $\mu>0$, number of iterations $N = S\cdot T$
        \STATE  $\phi_0 = x_0^0 = x^0$
        \FOR{ $s=0,1,2,\ldots, S-1$}
        \FOR{ $k=0,1,2,\ldots, T-1$ }
        \STATE{Uniformly randomly pick set $I_k$ from $\{1,\ldots, m\}$ such that $|I_k| = r$}
        \STATE{$g^k = \frac{1}{r}\sum\limits_{i\in I_k}\left(\hat{\nabla} f_i(x_s^k) - \hat{\nabla} f_i(\phi_s)\right) + \hat{\nabla} f(\phi_s)$}
        \STATE{$x_s^{k+1} = x_s^k - h^k_s g^k$}
        \ENDFOR
        \STATE  $\phi_{s+1} = x_{s+1}^0 = x_s^k$
        \ENDFOR
        \STATE{Pick $\xi$ uniformly at random from $\{0,\ldots, N-1\}$}
        \RETURN $x^{\xi}$
    \end{algorithmic}
\end{algorithm}

For a mini-batch $I\subseteq\{1,\ldots,m\}$ of size $r$, authors denote 
\begin{align*}
    \hat{\nabla} f_{I}(x) = \frac{1}{r}\sum\limits_{i\in I}\hat{\nabla} f_i(x)
\end{align*}
and the algorithm is the same as for SVRG (Algorithm \ref{alg:SVRG}), with the only difference that instead of true gradients update
\begin{align*}
    x^{k+1}_s = x^k_s - h^k_s v^k_s,~v^k_s=\nabla f_{I_k}(x^k_s) - \nabla f_{I_k}(x^0_s) + \nabla f(x^0_s)
\end{align*}
they use gradient estimations
\begin{align*}
    x^{k+1}_s = x^k_s - h^k_s \hat{v}^k_s,~\hat{v}^k_s=\hat{\nabla} f_{I_k}(x^k_s) - \hat{\nabla} f_{I_k}(x^0_s) + \hat{\nabla} f(x^0_s)
\end{align*}
This estimation $\hat{\nabla} f(x^0_s)$ is no longer unbiased for zeroth-order gradient estimations, and that is the main problem for the convergence analysis of this method. They show that under assumptions mentioned above ZO-SVRG algorithm after $N=S\cdot T$  (there $S$ is a number of epochs) steps ensures that 
\begin{align*}
    \text{RandGradEst}:&~\Exp\ls\|\nabla f(\bar{x})\|^2_{2}\rs = O\lp\frac{n}{N}+\frac{\delta_n}{r}\rp\\
    \text{Avg-RandGradEst}:&~\Exp\ls\|\nabla f(\bar{x})\|^2_{2}\rs =  O\lp\frac{n}{N}+\frac{\delta_n}{r\cdot \min\{n,q\}}\rp\\
    \text{CoordGradEst}:&~\Exp\ls\|\nabla f(\bar{x})\|^2_{2}\rs = 
    O\lp\frac{n}{N}\rp
\end{align*}
there $n$ is a dimension, $r = |I|$ -- batch size, $q$ is the number of directions used to estimate gradient via Avg-RandGradEst, $\bar{x}$ is uniformly chosen from $\{x^k_s\}_{s,k=0}^{S-1,T-1}$, $N=S\cdot T$ is a total number of steps and
\begin{align*}
    \delta_n = 
    \begin{cases}
        1,&\text{if $I_k$ draws samples from \{1,\ldots, m\} with replacement}\\
        j(b < n),&\ldots\text{ without replacement}
    \end{cases}
\end{align*}
where $j(b < n) = 1$ if $b<n$ and $j(b < n) = 0$ otherwise.

Basically, that means that CoordGradEst, the deterministic policy of gradient estimations, achieves the convergence rates of the original SVRG. In their tests, however, in terms of training loss versus function queries ZO-SVRG (the variant without mini-batching and with random directions on the sphere) beats ZO-SVRG-Ave (based on Avg-RandGradEst) and ZO-SVRG-Coord (based on CoordGradEst). 

\begin{algorithm}[ht]
    \caption{SpiderSZO \cite{fang2018spider}}
    \label{alg:SpiderSZO}
    \begin{algorithmic}
        \REQUIRE $n_0\in [1,\kesha{\nicefrac{n^{1/2}}{6}}]$, Lipschitz constant $L$, epoch length $T$, starting point $x^0\in\R^n$, outer batch size $r_1 \ge 1$, inner batch size $r_2 \ge 1$, number of iterations $N = S\cdot T$
        \FOR{ $k=0,1,2,\ldots, N - 1$ }
        \IF{$k\mod T = 0$}
        \STATE{Uniformly randomly pick set $I_k$ from $\{1,\ldots, m\}$ (with replacement) such that $|I_k| = r_1$}
        \STATE Compute $g^k = \sum\limits_{j = 1}^{n}\lp\frac{1}{r_1}\sum\limits_{i\in I_k}\frac{[f_{i}(x^k+\mu e_j) - f_i(x^k)]}{\mu}\rp e_j$
        \ELSE
        \STATE{Create set of pairs $I_k = \{(i, u^i)\}$ where $i$ uniformly randomly picked from $\{1,\ldots, m\}$ (with replacement) and independent $u^i\sim\mathcal{N}(0,I_n)$ such that $|I_k| = r_2$}
        \STATE{Compute $g^k = \frac{1}{r_2}\sum\limits_{(i,u^i)\in I_k}\lp\frac{f_i(x^k+\mu u^i)-f_i(x^k)}{\mu}u^i - \frac{f_i(x^{k-1}+\mu u^i)-f_i(x^{k-1})}{\mu}u^i\rp + g^{k-1}$}
        \ENDIF
        \STATE{$x^{k+1} = x^k - h_k g^k$ where $h_k = \min\lp\frac{\eps}{L n_0\|v^k\|_{2}},~\frac{1}{2Ln_0}\rp$}
        \ENDFOR
        \STATE{Pick $\xi$ uniformly at random from $\{0,\ldots, N-1\}$}
        \RETURN $x^{\xi}$
    \end{algorithmic}
\end{algorithm}

Another discussed above algorithm that can be used in the zeroth-order finite-sum minimization setting is SPIDER \cite{fang2018spider}. The zeroth-order variant (Algorithm \ref{alg:SpiderSZO}) of the algorithm blends stochastic and deterministic gradient estimations, using mini-batched $FFD$ (Table \ref{tab:zo_gradient_estimations_comparison}) every $p$ steps to reconstruct $v^k$, which is later updated by mini-batched $GSG$. 

The $h_k = \min\lp\kesha{\nicefrac{\eps}{\ls L n_0\|v^k\|_{2}\rs},~\nicefrac{1}{\ls 2Ln_0\rs}}\rp$ is a stepsize policy from Normalized Gradient Descent (NGD, \cite{nesterov2004introduction}), where the stepsize is inverse-proportional to the norm of the gradient. 

Authors show, that after $N = O\lp\kesha{\nicefrac{1}{\eps^2}}\rp$ iterations and $O\lp n\min\lp \kesha{\nicefrac{m^{1/2}}{\eps^2},\nicefrac{1}{\eps^3}}\rp\rp$ (there $n$ is a dimension and $m$ is a number of functions) IZO calls (i.e. calls of the oracle that returns the value of $f_{i}(x)$ given $x$ and $i$) this algorithm ensures
\begin{align*}
    \Exp[\|\nabla f(\bar{x})\|_{2}]\leqslant 6\eps
\end{align*}
where $\bar{x}$ is uniformly chosen from $\{x^k\}_{k=0}^{N-1}$. This result is better than what follows directly from \cite{nesterov2017random}, at least by the factor of $m^{\kesha{\nicefrac{1}{2}}}$ (the direct application of the results from \cite{nesterov2017random} requires $m$ calls on every step, and gives $\Exp[\|\nabla f(\bar{x})\|_{2}]\leqslant \eps$ in $O\lp\kesha{\nicefrac{n}{\eps^2}}\rp$ steps so the number of IZO calls would be $O\lp\kesha{\nicefrac{nm}{\eps^2}}\rp$).

The results of two previously discussed papers \cite{liu2018zerothorder, fang2018spider} were improved in the recent work \cite{ji2019improved}. Authors show that ZO-SVRG-Coord actually has a better convergence rate (\cite{ji2019improved}[Theorem 2]) of $\Exp\ls\|\nabla f(\bar{x})\|^2_{2}\rs = O\lp\kesha{\nicefrac{1}{N}}\rp$ ($n$ times better than the previous analysis). At first they consider an intermediate variant of ZO-SVRG-Coord and ZO-SVRG-Ave called ZO-SVRG-Coord-Rand, that uses $CFD$ and $BSG$ (Table \ref{tab:zo_gradient_estimations_comparison}) for $\hat{\nabla} f(\phi_s)$ and $\hat{\nabla} f_i(x_s^k) - \hat{\nabla} f_i(\phi_s)$ parts of 
\begin{align*}
    g^k = \frac{1}{r}\sum\limits_{i\in I_k}\left(\hat{\nabla} f_i(x_s^k) - \hat{\nabla} f_i(\phi_s)\right) + \hat{\nabla} f(\phi_s)
\end{align*}
(from Algorithm \ref{alg:ZO-SVRG}) respectively, while variants in \cite{liu2018zerothorder} used only one type of gradient estimation at once. Then authors  proof (\cite{ji2019improved}[Corollary 1]) the convergence rate $\Exp\ls\|\nabla f(\bar{x})\|^2_{2}\rs = O\lp\kesha{\nicefrac{1}{N}}\rp$ and show (\cite{ji2019improved}[Lemmas 1-2]) that although the replacement of $BSG$ with $CFD$ requires $n$ more oracle calls it achieves more accurate gradient estimation so the convergence rate stays the same for the ZO-SVRG-Coord.

Another part of this work is devoted to SPIDER. Authors construct a new algorithm (called ZO-SPIDER-Coord) in a way similar to the previous one -- they use $CFD$ instead of $GSG$ in Algorithm \ref{alg:SpiderSZO} and show that it has the same rate of convergence, but with bigger stepsize $h_k = \kesha{\nicefrac{1}{\ls 4L \rs}}$ (that doesn't depend on $\eps$), which is better in practice. 

One particular case of finite-sum minimization is considered in \cite{zhang2020boosting}. In this paper, authors consider non-stationary online optimization problems, when the objective function being queried is
time-varying, so one is limited to the use of one-point estimators.

Such estimators can be constructed easily in the stochastic zeroth-order case. For example we can consider $GSG$ (Table \ref{tab:zo_gradient_estimations_comparison}) with $N=1$ then 
\begin{align*}
    \Exp_{u}(g(x)) = \Exp_{u}\ls\frac{f(x+\mu u)-f(x)}{\mu}u\rs = \Exp_{u}\ls\frac{f(x+\mu u)}{\mu}u\rs = \nabla f_{\mu}(x)
\end{align*}
so we can chose $g(x) := \kesha{\ls\nicefrac{f(x+\mu u)}{\mu}\rs}u$ and obtain a reasonable one-point estimation. The problem is that the variance of such estimations explodes as $\mu\to 0$ (see \cite{berahas2019global}).

In this work, authors consider the residual feedback estimator
\begin{align*}
    \tilde{g}_k(x^k) := \frac{u^k}{\mu}\lp f_k(x^k+\mu u^k) - f_{k-1}(x^{k-1}+\mu u^{k-1})\rp
\end{align*}
where $u^k,u^{k-1}\sim\mathcal{N}(0,I_n)$. They show that (Lemma 2.4)
\begin{align*}
    \mathbb{E}[\tilde{g}_k(x^k)] = \nabla f_{\mu,k}(x^k),~\forall x^k\in X\text{ and }k
\end{align*}
(there $\nabla f_{\mu,k}$ is a gradient of smoothed $f_k$). They consider the online bandit problem with regret function 
\begin{align*}
    R^T_{g,\mu} = \sum\limits_{k=0}^{T-1}\Exp\ls\|\nabla f_{\mu,k}(x^k)\|^2_{2}\rs
\end{align*}
and show (\cite{zhang2020boosting}[Theorem 4.2]) that for $x^{k+1} = \Pi_{X}\lp x^k-\eta\tilde{g}_k(x^k)\rp$ (where $\Pi_{X}$ is the projection operator onto set $X$) if $f\in C^{0,0}_{L_0}$
\begin{align*}
    R^T_{g,\mu} = O\lp \frac{n^{\kesha{\nicefrac{3}{2}}}L_0^2}{\eps_f^{\kesha{\nicefrac{3}{2}}}}\lp W_T+\tilde{W}_T T^{-1}\rp T^{\kesha{\nicefrac{1}{2}}} + n^{\kesha{\nicefrac{3}{2}}}L_0 \eps_f^{\kesha{\nicefrac{1}{2}}}T^{\kesha{\nicefrac{1}{2}}}\rp
\end{align*}
and if additionally $f\in C^{1,1}_{L_1}$ (\cite{zhang2020boosting}[Theorem 4.3])
\begin{align*}
    R^T_{g} = \sum\limits_{k=0}^{T-1}\Exp\ls\|\nabla f_{k}(x^k)\|^2_{2}\rs = O\lp n^{\kesha{\nicefrac{4}{3}}}L_0W_T T^{\kesha{\nicefrac{1}{2}}} + n^{\kesha{\nicefrac{4}{3}}}L_1L_0^{-1} \tilde{W}_T\rp
\end{align*}
where $W_T$ and $\tilde{W}_T$ are constants s.t. 
\begin{align*}
    \sum\limits_{k=1}^{T}\Exp\ls f_k(x) - f_{k-1}(x)\rs \leqslant W_T,~\forall T,x \\
    \sum\limits_{k=1}^{T}\Exp\ls |f_k(x) - f_{k-1}(x)|^2\rs \leqslant \tilde{W}_T,~\forall T,x.
\end{align*}

That bound implies that $\nicefrac{R^T_{g}}{T}\to 0$ if $W_T=o\lp T^{\kesha{\nicefrac{1}{2}}}\rp$ and $\tilde{W}_T = o\lp T\rp$. Authors also consider (\cite{zhang2020boosting}[Section 5]) the stochastic online optimization case where $\hat{f}_{t} = F_{t}(x,\xi_t)$ s.t. $\Exp\ls F_{t}(x,\xi_t)\rs = f_t(x)$ and show that under the assumptions of the same form as above (with $W_{T,\xi}$ and $\tilde{W}_{T,\xi}$) similar regret bounds can be achieved. 

In their numerical experiments, authors compare conventional one-point and two-point approaches with one-point residual feedback. Even though the latter works worse than the two-point variant, it has lower variance and achieves better results than the conventional one-point feedback, and can be used in practice, in contrast to two-point feedback. 

\section{Globalization Techniques}
In the previous sections we mainly considered guarantees for the methods to converge to a stationary point or local extremum. Global performance guarantees are available only for some subclasses of non-convex minimization problems. Despite that there are several practical techniques for convergence globalization for the local methods, which we briefly describe next, following \cite{zhigljavsky2007stochastic}.

\subsection{Multistart Technique}

The first approach involves using an algorithm which converges to a local minimum and running it multiple times from different starting points. This may result in the algorithm for finding multiple local minima of the objective, some of which might in fact be global solutions.

To be more concrete, we consider the problem
\[
\mathop {\min }\limits_{x\in \left[ {0,1} \right]^n}  f\left( x \right) .
\]
Let the initial points be sampled from the uniform distribution on $[0,1]^n$. If the Lebesgue measure of the attraction basin (the set of points, initialized at which the local algorithm converges to the global minimum) of the global minimum is $\mu>0$, then the expected number of points required to find the global minimum is $m=\tilde {{\rm O}}\left( 1/\mu \right)$.
If the attraction basin is a ball of radius $r$, then $\mu\sim r^n$. Hence, it is reasonable to expect that the number of initial points required depends on $n$ exponentially. For that reason, this approach to global optimization becomes impractical as $n$ grows.

The effectiveness of this approach also depends on the chosen initial points. The quality of a family of initial points $\left\{ {x^{0,i}} \right\}_{i=1}^m $ can be characterized by the quantity \[
d_n \left( {\left\{ {x^{0,i}} \right\}_{i=1}^m } \right)=\mathop {\max 
}\limits_{x\in \left[ {0,1} \right]^n} \mathop {\min }\limits_{i=1,...,m} 
\left\| {x-x^{0,i}} \right\|_2.
\]

One of the ways to iteratively generate the starting points $\left\{ {x^{0,k}} 
\right\}_{k=1}^m $ is called the quasi Monte Carlo scheme using low-discrepancy sequences, for example, the Van der Corput sequence. Let $\left\{ {p_i } \right\}_{i=1}^n $ be a sequence of distinct prime numbers, and let $\phi _{i } \left( k \right)$ be the $k$-th element of the Van der Corput sequence in base $p_i$. Explicitly, $\phi_{i}(k)=\sum\limits_{j=0}^{l_{k,i } } {a_j 
p_i^{-j-1} } $, where $l_{k,i}$ is the length of the representation of $k$ in base $p_i$ $k=\sum\limits_{j=0}^{l_{k,i } } {a_j p_i^j }$.
Finally, set $x^{0,k}=\left( {\phi _{1 } \left( k \right),..1.,\phi_{n} \left( k 
\right)} \right)$, $k=1,...,m$.
In this case $d_n \left( {\left\{ {x^{0,i}} \right\}_{i=1}^m } \right)={\rm
O}\left( {\sqrt n m^{-1 \mathord{\left/ {\vphantom {1 n}} \right. 
\kern-\nulldelimiterspace} n}\ln m} \right)$, while the optimal value, which is achieved at the uniform grid, is ${\rm O}\left( {\sqrt n m^{-1 \mathord{\left/ {\vphantom {1 {\left( {2n} \right)}}} \right. \kern-\nulldelimiterspace} {\left( {2n} \right)}}} \right).$

\subsection{Multidimensional Bisection}

The main shortcoming of the approach described above is that the family $\left\{ 
{x^{0,k}} \right\}_{k=1}^m $ is constructed without taking into account any properties of $f\left( x \right)$. 
Assume now that, for all $ x,y\in[0,1]^n$ $,
\left| {f\left( y \right)-f\left( x \right)} \right|\leqslant M\left\| {y-x} 
\right\|$.
Then, for any $y$, the function $f(y)-M\|x-y\|$ is a minorant of $f(x)$. Consequently, for any $\{y^k\}_{k=1}^m$ the function $\max\limits_{k=1,\ldots,m} f(y^k)-M\|x-y^k\|$ is also a minorant of $f(x)$. Then one may choose the next initial point to be the minimizer of the minorant constructed using the previous initial points \alexander{\cite{evtushenko1971numerical}}:
\[
x^{0,m+1}=\arg\min\limits_x \mathop {\max }\limits_{k=1,...,m} 
\left\{ {f\left( {x^{0,k}} \right)-M\left\| {x-x^{0,k}} \right\|} \right\}.
\]
In the one-dimensional case, each minorant is just a piecewise linear function, and its minimum is easy to compute explicitly. In higher-dimensions, this idea is more difficult to implement, and the resulting algorithms also tend to become slower as $n$ increases. This method also requires an estimate of the Lipschitz constant and is sensitive to the accuracy of this estimate.

\subsection{Langevin Dynamics}
The last but not least approach which we consider in this section is inspired by the Langevin dynamics, which is defined by the stochastic differential equation
\[
dx(t)=-\nabla f(x(t))dt+\sqrt {2T} dW\left( t \right),
\]
where $W(t)$ is a Wiener process (also known as Brownian motion) and $T$ is the temperature parameter. It has been shown that the distribution of $x\left( t \right)$ converges to a distribution with density\[
\frac{\exp \left( {{-f\left( x \right)} \mathord{\left/ {\vphantom 
{{-f\left( x \right)} T}} \right. \kern-\nulldelimiterspace} T} 
\right)}{\int {\exp \left( {{-f\left( y \right)} \mathord{\left/ {\vphantom 
{{-f\left( y \right)} T}} \right. \kern-\nulldelimiterspace} T} \right)} 
dy}
\] as $t\to\infty$, and as $T\to 0+$ this distribution concentrates around the global minima.
To apply this in practice, the continuous dynamics has to be discretized. One of the ways to do that is as follows:
\[x_{k+1} = x_k-h\nabla f(x_k) + \sqrt {2hT}\epsilon_k,
\] where $h>0$ is the stepsize and $\epsilon_k$ is standard gaussian random variable. 
%with different policy of choosing $T\left( t \right)$, such that. 
Non-asymptotic results demonstrating the convergence of this method to an approximate global minimum were presented in the work \cite{xu2018global}. In this paper, the temperature parameter $T$ was assumed to be constant. However, other strategies are sometimes used in practice, for example,\[
T_k =\frac{c}{\ln \left( {2+k} \right)},
\]
which ensures $T_k\to 0+$ as $k\to\infty$. 

\section*{Acknowledgements}
The authors are grateful to A.~Gornov, A.~Nazin, Yu.~Nesterov, B.~Polyak and K. Scheinberg for fruitful discussions and their suggestions which helped to improve the quality of the text.  

The research was partially supported by the Ministry of Science and Higher Education of the Russian Federation (Goszadaniye) No.075-00337-20-03, project No. 0714-2020-0005.

\bibliographystyle{abbrv}
\bibliography{published, nonpublished}

\begin{thebibliography}{100}

\bibitem{pytorchOpt}
Collection of optimizers for pytorch.
\newblock \url{https://github.com/jettify/pytorch-optimizer}.

\bibitem{agarwal2017finding}
N.~Agarwal, Z.~Allen-Zhu, B.~Bullins, E.~Hazan, and T.~Ma.
\newblock Finding approximate local minima faster than gradient descent.
\newblock In {\em Proceedings of the 49th Annual ACM SIGACT Symposium on Theory
  of Computing}, pages 1195--1199, 2017.

\bibitem{ahn2020sgd}
K.~Ahn, C.~Yun, and S.~Sra.
\newblock Sgd with shuffling: optimal rates without component convexity and
  large epoch requirements.
\newblock {\em Advances in Neural Information Processing Systems}, 33, 2020.

\bibitem{ajalloeian2020analysis}
A.~Ajalloeian and S.~U. Stich.
\newblock Analysis of sgd with biased gradient estimators.
\newblock {\em arXiv preprint arXiv:2008.00051}, 2020.

\bibitem{alistarh2017qsgd}
D.~Alistarh, D.~Grubic, J.~Li, R.~Tomioka, and M.~Vojnovic.
\newblock Qsgd: Communication-efficient sgd via gradient quantization and
  encoding.
\newblock In {\em Advances in Neural Information Processing Systems}, pages
  1709--1720, 2017.

\bibitem{allen2017natasha}
Z.~Allen-Zhu.
\newblock Natasha: Faster non-convex stochastic optimization via strongly
  non-convex parameter.
\newblock In {\em International Conference on Machine Learning}, pages 89--97,
  2017.

\bibitem{allen2018make}
Z.~Allen-Zhu.
\newblock How to make the gradients small stochastically: Even faster convex
  and nonconvex sgd.
\newblock In {\em Advances in Neural Information Processing Systems}, pages
  1157--1167, 2018.

\bibitem{allen2018katyusha}
Z.~Allen-Zhu.
\newblock Katyusha x: Simple momentum method for stochastic sum-of-nonconvex
  optimization.
\newblock In {\em International Conference on Machine Learning}, pages
  179--185, 2018.

\bibitem{allen2018natasha}
Z.~Allen-Zhu.
\newblock Natasha 2: Faster non-convex optimization than sgd.
\newblock In {\em Advances in Neural Information Processing Systems}, pages
  2675--2686, 2018.

\bibitem{allen2018neon2}
Z.~Allen-Zhu and Y.~Li.
\newblock Neon2: Finding local minima via first-order oracles.
\newblock In {\em Advances in Neural Information Processing Systems}, pages
  3716--3726, 2018.

\bibitem{allen2019can}
Z.~Allen-Zhu and Y.~Li.
\newblock Can sgd learn recurrent neural networks with provable generalization?
\newblock In {\em Advances in Neural Information Processing Systems}, pages
  10331--10341, 2019.

\bibitem{allen2019learning}
Z.~Allen-Zhu, Y.~Li, and Y.~Liang.
\newblock Learning and generalization in overparameterized neural networks,
  going beyond two layers.
\newblock In {\em Advances in neural information processing systems}, pages
  6158--6169, 2019.

\bibitem{allen2019convergence}
Z.~Allen-Zhu, Y.~Li, and Z.~Song.
\newblock A convergence theory for deep learning via over-parameterization.
\newblock In {\em International Conference on Machine Learning}, pages
  242--252. PMLR, 2019.

\bibitem{allen2019on}
Z.~Allen-Zhu, Y.~Li, and Z.~Song.
\newblock On the convergence rate of training recurrent neural networks.
\newblock In {\em Advances in neural information processing systems}, pages
  6676--6688, 2019.

\bibitem{anandkumar2016efficient}
A.~Anandkumar and R.~Ge.
\newblock Efficient approaches for escaping higher order saddle points in
  non-convex optimization.
\newblock In {\em Conference on learning theory}, pages 81--102. PMLR, 2016.

\bibitem{arjevani2020second}
Y.~Arjevani, Y.~Carmon, J.~C. Duchi, D.~J. Foster, A.~Sekhari, and
  K.~Sridharan.
\newblock Second-order information in non-convex stochastic optimization: Power
  and limitations.
\newblock In {\em Conference on Learning Theory}, pages 242--299, 2020.

\bibitem{arjevani2019lower}
Y.~Arjevani, Y.~Carmon, J.~C. Duchi, D.~J. Foster, N.~Srebro, and B.~Woodworth.
\newblock Lower bounds for non-convex stochastic optimization.
\newblock {\em arXiv preprint arXiv:1912.02365}, 2019.

\bibitem{arora2018convergence}
S.~Arora, N.~Cohen, N.~Golowich, and W.~Hu.
\newblock A convergence analysis of gradient descent for deep linear neural
  networks.
\newblock {\em arXiv preprint arXiv:1810.02281}, 2018.

\bibitem{bach2012optimization}
F.~Bach, R.~Jenatton, J.~Mairal, G.~Obozinski, et~al.
\newblock Optimization with sparsity-inducing penalties.
\newblock {\em Foundations and Trends{\textregistered} in Machine Learning},
  4(1):1--106, 2012.

\bibitem{baraniuk2008simple}
R.~Baraniuk, M.~Davenport, R.~DeVore, and M.~Wakin.
\newblock A simple proof of the restricted isometry property for random
  matrices.
\newblock {\em Constructive Approximation}, 28(3):253--263, 2008.

\bibitem{bazarova2020linearly}
A.~Bazarova, A.~Beznosikov, and A.~Gasnikov.
\newblock Linearly convergent gradient-free methods for minimization of
  symmetric parabolic approximation.
\newblock {\em arXiv preprint arXiv:2009.04906}, 2020.

\bibitem{ben-tal2001lectures}
A.~Ben-Tal and A.~Nemirovski.
\newblock {\em Lectures on Modern Convex Optimization.}
\newblock Society for Industrial and Applied Mathematics, 2001.

\bibitem{berahas2019linear}
A.~S. Berahas, L.~Cao, K.~Choromanski, and K.~Scheinberg.
\newblock Linear interpolation gives better gradients than gaussian smoothing
  in derivative-free optimization, 2019.

\bibitem{berahas2020theoretical}
A.~S. Berahas, L.~Cao, K.~Choromanski, and K.~Scheinberg.
\newblock A theoretical and empirical comparison of gradient approximations in
  derivative-free optimization, 2020.

\bibitem{berahas2019global}
A.~S. Berahas, L.~Cao, and K.~Scheinberg.
\newblock Global convergence rate analysis of a generic line search algorithm
  with noise, 2019.

\bibitem{bergou2019stochastic}
E.~H. Bergou, E.~Gorbunov, and P.~Richtárik.
\newblock Stochastic three points method for unconstrained smooth minimization,
  2019.

\bibitem{beznosikov2020biased}
A.~Beznosikov, S.~Horv{\'a}th, P.~Richt{\'a}rik, and M.~Safaryan.
\newblock On biased compression for distributed learning.
\newblock {\em arXiv preprint arXiv:2002.12410}, 2020.

\bibitem{bhojanapalli2016dropping}
S.~Bhojanapalli, A.~Kyrillidis, and S.~Sanghavi.
\newblock Dropping convexity for faster semi-definite optimization.
\newblock In {\em Conference on Learning Theory}, pages 530--582, 2016.

\bibitem{birgin2017worst}
E.~G. Birgin, J.~Gardenghi, J.~M. Mart{\'\i}nez, S.~A. Santos, and P.~L. Toint.
\newblock Worst-case evaluation complexity for unconstrained nonlinear
  optimization using high-order regularized models.
\newblock {\em Mathematical Programming}, 163(1-2):359--368, 2017.

\bibitem{blum2016foundations}
A.~Blum, J.~Hopcroft, and R.~Kannan.
\newblock {\em Foundations of data science}.
\newblock Cambridge University Press, 2016.

\bibitem{blum1989training}
A.~Blum and R.~L. Rivest.
\newblock Training a 3-node neural network is np-complete.
\newblock In {\em Advances in neural information processing systems}, pages
  494--501, 1989.

\bibitem{blumensath2009iterative}
T.~Blumensath and M.~E. Davies.
\newblock Iterative hard thresholding for compressed sensing.
\newblock {\em Applied and computational harmonic analysis}, 27(3):265--274,
  2009.

\bibitem{bogolubsky2016learning}
L.~Bogolubsky, P.~Dvurechensky, A.~Gasnikov, G.~Gusev, Y.~Nesterov, A.~M.
  Raigorodskii, A.~Tikhonov, and M.~Zhukovskii.
\newblock Learning supervised pagerank with gradient-based and gradient-free
  optimization methods.
\newblock In D.~D. Lee, M.~Sugiyama, U.~V. Luxburg, I.~Guyon, and R.~Garnett,
  editors, {\em Advances in Neural Information Processing Systems 29}, pages
  4914--4922. Curran Associates, Inc., 2016.
\newblock arXiv:1603.00717.

\bibitem{bottou2009curiously}
L.~Bottou.
\newblock Curiously fast convergence of some stochastic gradient descent
  algorithms.
\newblock In {\em Proceedings of the symposium on learning and data science,
  Paris}, 2009.

\bibitem{bottou2010large}
L.~Bottou.
\newblock Large-scale machine learning with stochastic gradient descent.
\newblock In {\em Proceedings of COMPSTAT'2010}, pages 177--186. Springer,
  2010.

\bibitem{bottou2012stochastic}
L.~Bottou.
\newblock Stochastic gradient descent tricks.
\newblock In {\em Neural networks: Tricks of the trade}, pages 421--436.
  Springer, 2012.

\bibitem{bottou2018optimization}
L.~Bottou, F.~E. Curtis, and J.~Nocedal.
\newblock Optimization methods for large-scale machine learning.
\newblock {\em Siam Review}, 60(2):223--311, 2018.

\bibitem{boyd2004convex}
S.~Boyd and L.~Vandenberghe.
\newblock {\em Convex Optimization}.
\newblock NY Cambridge University Press, 2004.

\bibitem{bu2020note}
J.~Bu and M.~Mesbahi.
\newblock A note on {N}esterov's accelerated method in nonconvex optimization:
  a weak estimate sequence approach.
\newblock {\em arXiv preprint arXiv:2006.08548}, 2020.

\bibitem{bubeck2011introduction}
S.~Bubeck.
\newblock Introduction to online optimization.
\newblock 2011.

\bibitem{bubeck2015convex}
S.~Bubeck.
\newblock Convex optimization: Algorithms and complexity.
\newblock {\em Found. Trends Mach. Learn.}, 8(3–4):231–357, nov 2015.

\bibitem{candes2015phase}
E.~J. Candes, X.~Li, and M.~Soltanolkotabi.
\newblock Phase retrieval via wirtinger flow: Theory and algorithms.
\newblock {\em IEEE Transactions on Information Theory}, 61(4):1985--2007,
  2015.

\bibitem{candes2009exact}
E.~J. Cand{\`e}s and B.~Recht.
\newblock Exact matrix completion via convex optimization.
\newblock {\em Foundations of Computational mathematics}, 9(6):717, 2009.

\bibitem{candes2005decoding}
E.~J. Candes and T.~Tao.
\newblock Decoding by linear programming.
\newblock {\em IEEE transactions on information theory}, 51(12):4203--4215,
  2005.

\bibitem{candes2010power}
E.~J. Cand{\`e}s and T.~Tao.
\newblock The power of convex relaxation: Near-optimal matrix completion.
\newblock {\em IEEE Transactions on Information Theory}, 56(5):2053--2080,
  2010.

\bibitem{candes2008enhancing}
E.~J. Candes, M.~B. Wakin, and S.~P. Boyd.
\newblock Enhancing sparsity by reweighted $\ell_1$ minimization.
\newblock {\em Journal of Fourier analysis and applications}, 14(5-6):877--905,
  2008.

\bibitem{carmon2016gradient}
Y.~Carmon and J.~C. Duchi.
\newblock Gradient descent efficiently finds the cubic-regularized non-convex
  newton step.
\newblock {\em arXiv preprint arXiv:1612.00547}, 2016.

\bibitem{carmon2017convex}
Y.~Carmon, J.~C. Duchi, O.~Hinder, and A.~Sidford.
\newblock ``{C}onvex until proven guilty'': Dimension-free acceleration of
  gradient descent on non-convex functions.
\newblock volume~70 of {\em Proceedings of Machine Learning Research}, pages
  654--663, International Convention Centre, Sydney, Australia, 06--11 Aug
  2017. PMLR.

\bibitem{carmon2018accelerated}
Y.~Carmon, J.~C. Duchi, O.~Hinder, and A.~Sidford.
\newblock Accelerated methods for nonconvex optimization.
\newblock {\em SIAM Journal on Optimization}, 28(2):1751--1772, 2018.

\bibitem{carmon2019lowerII}
Y.~Carmon, J.~C. Duchi, O.~Hinder, and A.~Sidford.
\newblock Lower bounds for finding stationary points {II}: first-order methods.
\newblock {\em Mathematical Programming}, Sep 2019.

\bibitem{carmon2019lowerI}
Y.~Carmon, J.~C. Duchi, O.~Hinder, and A.~Sidford.
\newblock Lower bounds for finding stationary points i.
\newblock {\em Mathematical Programming}, 184(1):71--120, Nov 2020.

\bibitem{cartis2011adaptive}
C.~Cartis, N.~I. Gould, and P.~L. Toint.
\newblock Adaptive cubic regularisation methods for unconstrained optimization.
  part i: motivation, convergence and numerical results.
\newblock {\em Mathematical Programming}, 127(2):245--295, 2011.

\bibitem{cartis2019universal}
C.~Cartis, N.~I. Gould, and P.~L. Toint.
\newblock Universal regularization methods: Varying the power, the smoothness
  and the accuracy.
\newblock {\em SIAM Journal on Optimization}, 29(1):595--615, 2019.

\bibitem{cartis2011adaptive2}
C.~Cartis, N.~I.~M. Gould, and P.~L. Toint.
\newblock Adaptive cubic regularisation methods for unconstrained optimization.
  part ii: worst-case function- and derivative-evaluation complexity.
\newblock {\em Mathematical Programming}, 130(2):295--319, Dec 2011.

\bibitem{cartis2017improved}
C.~Cartis, N.~I.~M. Gould, and P.~L. Toint.
\newblock Improved second-order evaluation complexity for unconstrained
  nonlinear optimization using high-order regularized models.
\newblock {\em arXiv:1708.04044}, 2018.

\bibitem{charisopoulos2020entrywise}
V.~Charisopoulos, A.~R. Benson, and A.~Damle.
\newblock Entrywise convergence of iterative methods for eigenproblems.
\newblock {\em arXiv preprint arXiv:2002.08491}, 2020.

\bibitem{chen2018convergence}
X.~Chen, S.~Liu, R.~Sun, and M.~Hong.
\newblock On the convergence of a class of adam-type algorithms for non-convex
  optimization.
\newblock {\em arXiv preprint arXiv:1808.02941}, 2018.

\bibitem{chen2018harnessing}
Y.~Chen and Y.~Chi.
\newblock Harnessing structures in big data via guaranteed low-rank matrix
  estimation.
\newblock {\em arXiv preprint arXiv:1802.08397}, 2018.

\bibitem{chen2019gradient}
Y.~Chen, Y.~Chi, J.~Fan, and C.~Ma.
\newblock Gradient descent with random initialization: Fast global convergence
  for nonconvex phase retrieval.
\newblock {\em Mathematical Programming}, 176(1-2):5--37, 2019.

\bibitem{chen2018variance}
Z.~Chen and T.~Yang.
\newblock A variance reduction method for non-convex optimization with improved
  convergence under large condition number.
\newblock {\em arXiv preprint arXiv:1809.06754}, 2018.

\bibitem{chen2020momentum}
Z.~Chen and Y.~Zhou.
\newblock Momentum with variance reduction for nonconvex composition
  optimization.
\newblock {\em arXiv preprint arXiv:2005.07755}, 2020.

\bibitem{chi2018nonconvex}
Y.~Chi, Y.~M. Lu, and Y.~Chen.
\newblock Nonconvex optimization meets low-rank matrix factorization: An
  overview.
\newblock {\em arXiv preprint arXiv:1809.09573}, 2018.

\bibitem{combettes2011proximal}
P.~L. Combettes and J.-C. Pesquet.
\newblock Proximal splitting methods in signal processing.
\newblock In {\em Fixed-point algorithms for inverse problems in science and
  engineering}, pages 185--212. Springer, 2011.

\bibitem{conn2000trust}
A.~Conn, N.~Gould, and P.~Toint.
\newblock {\em Trust Region Methods}.
\newblock Society for Industrial and Applied Mathematics, 2000.

\bibitem{conn2009introduction}
A.~Conn, K.~Scheinberg, and L.~Vicente.
\newblock {\em Introduction to Derivative-Free Optimization}.
\newblock Society for Industrial and Applied Mathematics, 2009.

\bibitem{curtis2017optimization}
F.~E. Curtis and K.~Scheinberg.
\newblock Optimization methods for supervised machine learning: From linear
  models to deep learning.
\newblock {\em arXiv preprint arXiv:1706.10207}, 2017.

\bibitem{cutkosky2019momentum}
A.~Cutkosky and F.~Orabona.
\newblock Momentum-based variance reduction in non-convex sgd.
\newblock In {\em Advances in Neural Information Processing Systems}, pages
  15236--15245, 2019.

\bibitem{dang2015stochastic}
C.~D. Dang and G.~Lan.
\newblock Stochastic block mirror descent methods for nonsmooth and stochastic
  optimization.
\newblock {\em SIAM J. on Optimization}, 25(2):856--881, Apr. 2015.

\bibitem{davis2019stochastic}
D.~Davis and D.~Drusvyatskiy.
\newblock Stochastic model-based minimization of weakly convex functions.
\newblock {\em SIAM Journal on Optimization}, 29(1):207--239, 2019.

\bibitem{defazio2020understanding}
A.~Defazio.
\newblock Understanding the role of momentum in non-convex optimization:
  Practical insights from a lyapunov analysis.
\newblock {\em arXiv preprint arXiv:2010.00406}, 2020.

\bibitem{defazio2014SAGA}
A.~Defazio, F.~Bach, and S.~Lacoste-Julien.
\newblock Saga: A fast incremental gradient method with support for
  non-strongly convex composite objectives.
\newblock In {\em Proceedings of the 27th International Conference on Neural
  Information Processing Systems}, NIPS'14, pages 1646--1654, Cambridge, MA,
  USA, 2014. MIT Press.

\bibitem{defazio2019ineffectiveness}
A.~Defazio and L.~Bottou.
\newblock On the ineffectiveness of variance reduced optimization for deep
  learning.
\newblock In {\em Advances in Neural Information Processing Systems}, pages
  1753--1763, 2019.

\bibitem{defazio2014finito}
A.~Defazio, J.~Domke, et~al.
\newblock Finito: A faster, permutable incremental gradient method for big data
  problems.
\newblock In {\em International Conference on Machine Learning}, pages
  1125--1133, 2014.

\bibitem{defossez2020convergence}
A.~D{\'e}fossez, L.~Bottou, F.~Bach, and N.~Usunier.
\newblock On the convergence of adam and adagrad.
\newblock {\em arXiv preprint arXiv:2003.02395}, 2020.

\bibitem{demin2021necessary}
V.~Demin, D.~Nekhaev, I.~Surazhevsky, K.~Nikiruy, A.~Emelyanov, S.~Nikolaev,
  V.~Rylkov, and M.~Kovalchuk.
\newblock Necessary conditions for stdp-based pattern recognition learning in a
  memristive spiking neural network.
\newblock {\em Neural Networks}, 134:64--75, 2021.

\bibitem{devlin2018bert}
J.~Devlin, M.-W. Chang, K.~Lee, and K.~Toutanova.
\newblock Bert: Pre-training of deep bidirectional transformers for language
  understanding.
\newblock {\em arXiv preprint arXiv:1810.04805}, 2018.

\bibitem{diakonikolas2019generalized}
J.~Diakonikolas and M.~I. Jordan.
\newblock Generalized momentum-based methods: A {H}amiltonian perspective.
\newblock {\em arXiv preprint arXiv:1906.00436}, 2019.

\bibitem{ding2019spurious}
T.~Ding, D.~Li, and R.~Sun.
\newblock Spurious local minima exist for almost all over-parameterized neural
  networks.
\newblock 2019.

\bibitem{duchi2011adaptive}
J.~Duchi, E.~Hazan, and Y.~Singer.
\newblock Adaptive subgradient methods for online learning and stochastic
  optimization.
\newblock {\em Journal of Machine Learning Research}, 12(Jul.):2121--2159,
  2011.

\bibitem{duchi2013estimation}
J.~Duchi, M.~I. Jordan, and B.~McMahan.
\newblock Estimation, optimization, and parallelism when data is sparse.
\newblock In {\em Advances in Neural Information Processing Systems}, pages
  2832--2840, 2013.

\bibitem{dvinskikh2020line-search}
D.~Dvinskikh, A.~Ogaltsov, A.~Gasnikov, P.~Dvurechensky, and V.~Spokoiny.
\newblock On the line-search gradient methods for stochastic optimization.
\newblock {\em IFAC-PapersOnLine}, 53(2):1715--1720, 2020.
\newblock 21th IFAC World Congress, arXiv:1911.08380.

\bibitem{dvurechensky2017gradient}
P.~Dvurechensky.
\newblock Gradient method with inexact oracle for composite non-convex
  optimization.
\newblock {\em arXiv:1703.09180}, 2017.

\bibitem{dvurechensky2021accelerated}
P.~Dvurechensky, E.~Gorbunov, and A.~Gasnikov.
\newblock An accelerated directional derivative method for smooth stochastic
  convex optimization.
\newblock {\em European Journal of Operational Research}, 290(2):601 -- 621,
  2021.

\bibitem{dvurechensky2021first-order}
P.~Dvurechensky, S.~Shtern, and M.~Staudigl.
\newblock First-order methods for convex optimization.
\newblock {\em EURO Journal on Computational Optimization}, 9:100015, 2021.
\newblock arXiv:2101.00935.

\bibitem{dvurechensky2021hessian}
P.~Dvurechensky and M.~Staudigl.
\newblock Hessian barrier algorithms for non-convex conic optimization.
\newblock {\em arXiv:2111.00100}, 2021.

\bibitem{dvurechensky2019generalized}
P.~Dvurechensky, M.~Staudigl, and C.~A. Uribe.
\newblock Generalized self-concordant hessian-barrier algorithms.
\newblock {\em arXiv:1911.01522}, 2019.
\newblock WIAS Preprint No. 2693.

\bibitem{dvurechensky2020advances}
P.~E. Dvurechensky, A.~V. Gasnikov, E.~A. Nurminski, and F.~S. Stonyakin.
\newblock {\em Advances in Low-Memory Subgradient Optimization}, pages 19--59.
\newblock Springer International Publishing, Cham, 2020.
\newblock arXiv:1902.01572.

\bibitem{emmenegger2021oracle}
N.~Emmenegger, R.~Kyng, and A.~N. Zehmakan.
\newblock On the oracle complexity of higher-order smooth non-convex finite-sum
  optimization.
\newblock {\em arXiv preprint arXiv:2103.05138}, 2021.

\bibitem{evtushenko1971numerical}
Y.~G. Evtushenko.
\newblock Numerical methods for finding global extrema (case of a non-uniform
  mesh).
\newblock {\em USSR Computational Mathematics and Mathematical Physics},
  11(6):38--54, 1971.

\bibitem{fang2018spider}
C.~Fang, C.~J. Li, Z.~Lin, and T.~Zhang.
\newblock Spider: Near-optimal non-convex optimization via stochastic
  path-integrated differential estimator.
\newblock In {\em Advances in Neural Information Processing Systems}, pages
  689--699, 2018.

\bibitem{fang2019sharp}
C.~Fang, Z.~Lin, and T.~Zhang.
\newblock Sharp analysis for nonconvex sgd escaping from saddle points.
\newblock In {\em Conference on Learning Theory}, pages 1192--1234, 2019.

\bibitem{fatkhullin2020optimizing}
I.~Fatkhullin and B.~Polyak.
\newblock Optimizing static linear feedback: Gradient method.
\newblock {\em arXiv preprint arXiv:2004.09875}, 2020.

\bibitem{fazel2019global}
M.~Fazel, R.~Ge, S.~M. Kakade, and M.~Mesbahi.
\newblock Global convergence of policy gradient methods for the linear
  quadratic regulator, 2019.

\bibitem{feizi2017porcupine}
S.~Feizi, H.~Javadi, J.~Zhang, and D.~Tse.
\newblock Porcupine neural networks:(almost) all local optima are global.
\newblock {\em arXiv preprint arXiv:1710.02196}, 2017.

\bibitem{flaxman2005online}
A.~D. Flaxman, A.~T. Kalai, and H.~B. McMahan.
\newblock Online convex optimization in the bandit setting: Gradient descent
  without a gradient.
\newblock In {\em Proceedings of the Sixteenth Annual ACM-SIAM Symposium on
  Discrete Algorithms}, SODA '05, pages 385--394, Philadelphia, PA, USA, 2005.
  Society for Industrial and Applied Mathematics.

\bibitem{floudas2008encyclopedia}
C.~A. Floudas and P.~M. Pardalos.
\newblock {\em Encyclopedia of optimization}.
\newblock Springer Science \& Business Media, 2008.

\bibitem{GasnikovBook}
A.~Gasnikov.
\newblock {\em Universal gradient descent}.
\newblock MCCME, Moscow, 2021.

\bibitem{gasnikov2018power}
A.~Gasnikov, P.~Dvurechensky, M.~Zhukovskii, S.~Kim, S.~Plaunov, D.~Smirnov,
  and F.~Noskov.
\newblock About the power law of the pagerank vector component distribution.
  part 2. the buckley--osthus model, verification of the power law for this
  model, and setup of real search engines.
\newblock {\em Numerical Analysis and Applications}, 11(1):16--32, 2018.

\bibitem{ge2015escaping}
R.~Ge, F.~Huang, C.~Jin, and Y.~Yuan.
\newblock Escaping from saddle points—online stochastic gradient for tensor
  decomposition.
\newblock In {\em Conference on Learning Theory}, pages 797--842, 2015.

\bibitem{ge2015intersecting}
R.~Ge and J.~Zou.
\newblock Intersecting faces: Non-negative matrix factorization with new
  guarantees.
\newblock In {\em International Conference on Machine Learning}, pages
  2295--2303. PMLR, 2015.

\bibitem{ghadimi2013stochastic}
S.~Ghadimi and G.~Lan.
\newblock Stochastic first- and zeroth-order methods for nonconvex stochastic
  programming.
\newblock {\em SIAM Journal on Optimization}, 23(4):2341--2368, 2013.
\newblock arXiv:1309.5549.

\bibitem{ghadimi2016accelerated}
S.~Ghadimi and G.~Lan.
\newblock Accelerated gradient methods for nonconvex nonlinear and stochastic
  programming.
\newblock {\em Mathematical Programming}, 156(1):59--99, 2016.

\bibitem{ghadimi2016mini-batch}
S.~Ghadimi, G.~Lan, and H.~Zhang.
\newblock Mini-batch stochastic approximation methods for nonconvex stochastic
  composite optimization.
\newblock {\em Mathematical Programming}, 155(1):267--305, 2016.
\newblock arXiv:1308.6594.

\bibitem{ghadimi2019generalized}
S.~Ghadimi, G.~Lan, and H.~Zhang.
\newblock Generalized uniformly optimal methods for nonlinear programming.
\newblock {\em Journal of Scientific Computing}, 79(3):1854--1881, Jun 2019.

\bibitem{goemans1995improved}
M.~X. Goemans and D.~P. Williamson.
\newblock Improved approximation algorithms for maximum cut and satisfiability
  problems using semidefinite programming.
\newblock {\em Journal of the ACM (JACM)}, 42(6):1115--1145, 1995.

\bibitem{Goodfellow-et-al-2016}
I.~Goodfellow, Y.~Bengio, and A.~Courville.
\newblock {\em Deep Learning}.
\newblock MIT Press, 2016.
\newblock \url{http://www.deeplearningbook.org}.

\bibitem{goodfellow2016deep}
I.~Goodfellow, Y.~Bengio, A.~Courville, and Y.~Bengio.
\newblock {\em Deep learning}, volume~1.
\newblock MIT press Cambridge, 2016.

\bibitem{pmlr-v139-gorbunov21a}
E.~Gorbunov, K.~P. Burlachenko, Z.~Li, and P.~Richtarik.
\newblock Marina: Faster non-convex distributed learning with compression.
\newblock In M.~Meila and T.~Zhang, editors, {\em Proceedings of the 38th
  International Conference on Machine Learning}, volume 139 of {\em Proceedings
  of Machine Learning Research}, pages 3788--3798. PMLR, 18--24 Jul 2021.

\bibitem{gorbunov2020stochastic}
E.~Gorbunov, M.~Danilova, and A.~Gasnikov.
\newblock Stochastic optimization with heavy-tailed noise via accelerated
  gradient clipping.
\newblock In H.~Larochelle, M.~Ranzato, R.~Hadsell, M.~F. Balcan, and H.~Lin,
  editors, {\em Advances in Neural Information Processing Systems}, volume~33,
  pages 15042--15053. Curran Associates, Inc., 2020.

\bibitem{gorbunov2021near-optimal}
E.~Gorbunov, M.~Danilova, I.~Shibaev, P.~Dvurechensky, and A.~Gasnikov.
\newblock Near-optimal high probability complexity bounds for non-smooth
  stochastic optimization with heavy-tailed noise.
\newblock {\em arXiv:2106.05958}, 2021.

\bibitem{gorbunov2018accelerated}
E.~Gorbunov, P.~Dvurechensky, and A.~Gasnikov.
\newblock An accelerated method for derivative-free smooth stochastic convex
  optimization.
\newblock {\em arXiv preprint arXiv:1802.09022 (accepted to SIOPT)}, 2018.

\bibitem{gorbunov2020unified}
E.~Gorbunov, F.~Hanzely, and P.~Richt{\'a}rik.
\newblock A unified theory of sgd: Variance reduction, sampling, quantization
  and coordinate descent.
\newblock In {\em International Conference on Artificial Intelligence and
  Statistics}, pages 680--690, 2020.

\bibitem{gorbunov2020smtp}
E.~A. Gorbunov, A.~Bibi, O.~Sener, E.~H. Bergou, and P.~Richt{\'a}rik.
\newblock A stochastic derivative free optimization method with momentum.
\newblock In {\em ICLR}, 2020.

\bibitem{gorodetskiy2020delta}
A.~Gorodetskiy, A.~Shlychkova, and A.~I. Panov.
\newblock Delta schema network in model-based reinforcement learning.
\newblock In B.~Goertzel, A.~I. Panov, A.~Potapov, and R.~Yampolskiy, editors,
  {\em Artificial General Intelligence}, pages 172--182, Cham, 2020. Springer
  International Publishing.

\bibitem{gotmare2018closer}
A.~Gotmare, N.~S. Keskar, C.~Xiong, and R.~Socher.
\newblock A closer look at deep learning heuristics: Learning rate restarts,
  warmup and distillation.
\newblock {\em arXiv preprint arXiv:1810.13243}, 2018.

\bibitem{gower2020sgd}
R.~Gower, O.~Sebbouh, and N.~Loizou.
\newblock Sgd for structured nonconvex functions: Learning rates, minibatching
  and interpolation.
\newblock In {\em International Conference on Artificial Intelligence and
  Statistics}, pages 1315--1323. PMLR, 2021.

\bibitem{gower2019sgd}
R.~M. Gower, N.~Loizou, X.~Qian, A.~Sailanbayev, E.~Shulgin, and
  P.~Richt{\'a}rik.
\newblock Sgd: General analysis and improved rates.
\newblock In {\em International Conference on Machine Learning}, pages
  5200--5209, 2019.

\bibitem{goyal2017accurate}
P.~Goyal, P.~Doll{\'a}r, R.~Girshick, P.~Noordhuis, L.~Wesolowski, A.~Kyrola,
  A.~Tulloch, Y.~Jia, and K.~He.
\newblock Accurate, large minibatch sgd: Training imagenet in 1 hour.
\newblock {\em arXiv preprint arXiv:1706.02677}, 2017.

\bibitem{griewank1981generalized}
A.~O. Griewank.
\newblock Generalized descent for global optimization.
\newblock {\em Journal of optimization theory and applications}, 34(1):11--39,
  1981.

\bibitem{guminov2021combination}
S.~Guminov, P.~Dvurechensky, N.~Tupitsa, and A.~Gasnikov.
\newblock On a combination of alternating minimization and {N}esterov's
  momentum.
\newblock In {\em Proceedings of the 38th International Conference on Machine
  Learning}, volume 145 of {\em Proceedings of Machine Learning Research},
  Virtual, 18--24 Jul 2021. PMLR.
\newblock arXiv:1906.03622, WIAS Preprint No. 2695.

\bibitem{guminov2017accelerated}
S.~Guminov and A.~Gasnikov.
\newblock Accelerated methods for alpha-weakly-quasi-convex problems.
\newblock {\em arXiv preprint arXiv:1710.00797}, 2017.

\bibitem{guminov2019accelerated}
S.~V. Guminov, Y.~E. Nesterov, P.~E. Dvurechensky, and A.~V. Gasnikov.
\newblock Accelerated primal-dual gradient descent with linesearch for convex,
  nonconvex, and nonsmooth optimization problems.
\newblock {\em Doklady Mathematics}, 99(2):125--128, Mar 2019.

\bibitem{haeffele2017global}
B.~D. Haeffele and R.~Vidal.
\newblock Global optimality in neural network training.
\newblock In {\em Proceedings of the IEEE Conference on Computer Vision and
  Pattern Recognition}, pages 7331--7339, 2017.

\bibitem{haeser2019optimality}
G.~Haeser, H.~Liu, and Y.~Ye.
\newblock Optimality condition and complexity analysis for linearly-constrained
  optimization without differentiability on the boundary.
\newblock {\em Mathematical Programming}, 178(1):263--299, Nov 2019.

\bibitem{haochen2018random}
J.~Z. HaoChen and S.~Sra.
\newblock Random shuffling beats sgd after finite epochs.
\newblock {\em arXiv preprint arXiv:1806.10077}, 2018.

\bibitem{hazan2015beyond}
E.~Hazan, K.~Levy, and S.~Shalev-Shwartz.
\newblock Beyond convexity: Stochastic quasi-convex optimization.
\newblock In {\em Advances in Neural Information Processing Systems}, pages
  1594--1602, 2015.

\bibitem{he2016deep}
K.~He, X.~Zhang, S.~Ren, and J.~Sun.
\newblock Deep residual learning for image recognition.
\newblock In {\em Proceedings of the IEEE conference on computer vision and
  pattern recognition}, pages 770--778, 2016.

\bibitem{hinder2020near}
O.~Hinder, A.~Sidford, and N.~Sohoni.
\newblock Near-optimal methods for minimizing star-convex functions and beyond.
\newblock In {\em Conference on Learning Theory}, pages 1894--1938. PMLR, 2020.

\bibitem{hofmann2015variance}
T.~Hofmann, A.~Lucchi, S.~Lacoste-Julien, and B.~McWilliams.
\newblock Variance reduced stochastic gradient descent with neighbors.
\newblock In {\em Advances in Neural Information Processing Systems}, pages
  2305--2313, 2015.

\bibitem{horvath2019stochastic}
S.~Horv{\'a}th, D.~Kovalev, K.~Mishchenko, S.~Stich, and P.~Richt{\'a}rik.
\newblock Stochastic distributed learning with gradient quantization and
  variance reduction.
\newblock {\em arXiv preprint arXiv:1904.05115}, 2019.

\bibitem{ilyuhin2020recognition}
S.~A. Ilyuhin, A.~V. Sheshkus, and V.~L. Arlazarov.
\newblock {Recognition of images of Korean characters using embedded networks}.
\newblock In W.~Osten and D.~P. Nikolaev, editors, {\em Twelfth International
  Conference on Machine Vision (ICMV 2019)}, volume 11433, pages 273 -- 279.
  International Society for Optics and Photonics, SPIE, 2020.

\bibitem{jain2017non-convex}
P.~Jain and P.~Kar.
\newblock Non-convex optimization for machine learning.
\newblock {\em Found. Trends Mach. Learn.}, 10(3–4):142–336, Dec. 2017.

\bibitem{jain2013low}
P.~Jain, P.~Netrapalli, and S.~Sanghavi.
\newblock Low-rank matrix completion using alternating minimization.
\newblock In {\em Proceedings of the forty-fifth annual ACM symposium on Theory
  of computing}, pages 665--674, 2013.

\bibitem{ji2019improved}
K.~Ji, Z.~Wang, Y.~Zhou, and Y.~Liang.
\newblock Improved zeroth-order variance reduced algorithms and analysis for
  nonconvex optimization, 2019.

\bibitem{ji2019gradient}
Z.~Ji and M.~J. Telgarsky.
\newblock Gradient descent aligns the layers of deep linear networks.
\newblock In {\em 7th International Conference on Learning Representations,
  ICLR 2019}, 2019.

\bibitem{jin2017how}
C.~Jin, R.~Ge, P.~Netrapalli, S.~M. Kakade, and M.~I. Jordan.
\newblock How to escape saddle points efficiently.
\newblock volume~70 of {\em Proceedings of Machine Learning Research}, pages
  1724--1732, International Convention Centre, Sydney, Australia, 06--11 Aug
  2017. PMLR.

\bibitem{jin2019nonconvex}
C.~Jin, P.~Netrapalli, R.~Ge, S.~M. Kakade, and M.~I. Jordan.
\newblock On nonconvex optimization for machine learning: Gradients,
  stochasticity, and saddle points.
\newblock {\em Journal of the ACM (JACM)}, 68(2):1--29, 2021.

\bibitem{jin2018accelerated}
C.~Jin, P.~Netrapalli, and M.~I. Jordan.
\newblock Accelerated gradient descent escapes saddle points faster than
  gradient descent.
\newblock In {\em Conference On Learning Theory}, pages 1042--1085. PMLR, 2018.

\bibitem{johnson2013accelerating}
R.~Johnson and T.~Zhang.
\newblock Accelerating stochastic gradient descent using predictive variance
  reduction.
\newblock In {\em Advances in neural information processing systems}, pages
  315--323, 2013.

\bibitem{karimi2016linear}
H.~Karimi, J.~Nutini, and M.~Schmidt.
\newblock Linear convergence of gradient and proximal-gradient methods under
  the polyak-{\l}ojasiewicz condition.
\newblock In {\em Joint European Conference on Machine Learning and Knowledge
  Discovery in Databases}, pages 795--811. Springer, 2016.

\bibitem{khaled2020better}
A.~Khaled and P.~Richt{\'a}rik.
\newblock Better theory for sgd in the nonconvex world.
\newblock {\em arXiv preprint arXiv:2002.03329}, 2020.

\bibitem{khot2007optimal}
S.~Khot, G.~Kindler, E.~Mossel, and R.~O’Donnell.
\newblock Optimal inapproximability results for max-cut and other 2-variable
  csps?
\newblock {\em SIAM Journal on Computing}, 37(1):319--357, 2007.

\bibitem{khritankov2021hidden}
A.~Khritankov.
\newblock Hidden feedback loops in machine learning systems: A simulation model
  and preliminary results.
\newblock In D.~Winkler, S.~Biffl, D.~Mendez, M.~Wimmer, and J.~Bergsmann,
  editors, {\em Software Quality: Future Perspectives on Software Engineering
  Quality}, pages 54--65, Cham, 2021. Springer International Publishing.

\bibitem{kidambi2018insufficiency}
R.~Kidambi, P.~Netrapalli, P.~Jain, and S.~Kakade.
\newblock On the insufficiency of existing momentum schemes for stochastic
  optimization.
\newblock In {\em 2018 Information Theory and Applications Workshop (ITA)},
  pages 1--9. IEEE, 2018.

\bibitem{kiefer2020iterative}
L.~Kiefer, M.~Storath, and A.~Weinmann.
\newblock Iterative potts minimization for the recovery of signals with
  discontinuities from indirect measurements: The multivariate case.
\newblock {\em Foundations of Computational Mathematics}, pages 1--46, 2020.

\bibitem{kingma2014adam}
D.~P. Kingma and J.~Ba.
\newblock Adam: A method for stochastic optimization.
\newblock {\em arXiv preprint arXiv:1412.6980}, 2014.

\bibitem{kniaz2021adversarial}
V.~V. Kniaz, V.~A. Knyaz, V.~Mizginov, A.~Papazyan, N.~Fomin, and
  L.~Grodzitsky.
\newblock Adversarial dataset augmentation using reinforcement learning and 3d
  modeling.
\newblock In B.~Kryzhanovsky, W.~Dunin-Barkowski, V.~Redko, and Y.~Tiumentsev,
  editors, {\em Advances in Neural Computation, Machine Learning, and Cognitive
  Research IV}, pages 316--329, Cham, 2021. Springer International Publishing.

\bibitem{paparoditis2020wire}
V.~V. Kniaz, S.~Y. Zheltov, F.~Remondino, V.~A. Knyaz, and A.~Gruen.
\newblock Wire structure image-based 3d reconstruction aided by deep learning.
\newblock volume XLIII-B2-2020, pages 435 -- 441, Göttingen, 2020. Copernicus.
\newblock XXIV ISPRS Congress 2020 (virtual); Conference Location: Online;
  Conference Date: August 31 - September 2, 2020; Due to the Corona virus
  (COVID-19) the conference was conducted virtually.

\bibitem{kohler2017sub}
J.~M. Kohler and A.~Lucchi.
\newblock Sub-sampled cubic regularization for non-convex optimization.
\newblock In {\em International Conference on Machine Learning}, pages
  1895--1904, 2017.

\bibitem{kornowski2021oracle}
G.~Kornowski and O.~Shamir.
\newblock Oracle complexity in nonsmooth nonconvex optimization.
\newblock {\em arXiv preprint arXiv:2104.06763}, 2021.

\bibitem{kovalev2020don}
D.~Kovalev, S.~Horv{\'a}th, and P.~Richt{\'a}rik.
\newblock Don’t jump through hoops and remove those loops: Svrg and katyusha
  are better without the outer loop.
\newblock In {\em Algorithmic Learning Theory}, pages 451--467, 2020.

\bibitem{krizhevsky2009learning}
A.~Krizhevsky, G.~Hinton, et~al.
\newblock Learning multiple layers of features from tiny images.
\newblock 2009.

\bibitem{kuderov2021planning}
P.~Kuderov. and A.~Panov.
\newblock Planning with hierarchical temporal memory for deterministic markov
  decision problem.
\newblock In {\em Proceedings of the 13th International Conference on Agents
  and Artificial Intelligence - Volume 2: ICAART,}, pages 1073--1081. INSTICC,
  SciTePress, 2021.

\bibitem{lacroix2018canonical}
T.~Lacroix, N.~Usunier, and G.~Obozinski.
\newblock Canonical tensor decomposition for knowledge base completion.
\newblock In {\em International Conference on Machine Learning}, pages
  2863--2872, 2018.

\bibitem{lan2020first}
G.~Lan.
\newblock {\em First-order and Stochastic Optimization Methods for Machine
  Learning}.
\newblock Springer, 2020.

\bibitem{lan2019accelerated}
G.~Lan and Y.~Yang.
\newblock Accelerated stochastic algorithms for nonconvex finite-sum and
  multiblock optimization.
\newblock {\em SIAM Journal on Optimization}, 29(4):2753--2784, 2019.

\bibitem{Larson_2019}
J.~Larson, M.~Menickelly, and S.~M. Wild.
\newblock Derivative-free optimization methods.
\newblock {\em Acta Numerica}, 28:287–404, May 2019.

\bibitem{lee2016optimizing}
J.~C.~H. {Lee} and P.~{Valiant}.
\newblock Optimizing star-convex functions.
\newblock In {\em 2016 IEEE 57th Annual Symposium on Foundations of Computer
  Science (FOCS)}, pages 603--614, 2016.

\bibitem{lei2019stochastic}
Y.~Lei, T.~Hu, G.~Li, and K.~Tang.
\newblock Stochastic gradient descent for nonconvex learning without bounded
  gradient assumptions.
\newblock {\em IEEE Transactions on Neural Networks and Learning Systems},
  2019.

\bibitem{levy2016power}
K.~Y. Levy.
\newblock The power of normalization: Faster evasion of saddle points.
\newblock {\em arXiv preprint arXiv:1611.04831}, 2016.

\bibitem{li2018over}
D.~Li, T.~Ding, and R.~Sun.
\newblock Over-parameterized deep neural networks have no strict local minima
  for any continuous activations.
\newblock {\em arXiv preprint arXiv:1812.11039}, 2018.

\bibitem{li2016identifiability}
Y.~Li, K.~Lee, and Y.~Bresler.
\newblock Identifiability in blind deconvolution with subspace or sparsity
  constraints.
\newblock {\em IEEE Transactions on information Theory}, 62(7):4266--4275,
  2016.

\bibitem{li2020page}
Z.~Li, H.~Bao, X.~Zhang, and P.~Richt{\'a}rik.
\newblock Page: A simple and optimal probabilistic gradient estimator for
  nonconvex optimization.
\newblock {\em arXiv preprint arXiv:2008.10898}, 2020.

\bibitem{li2020unified}
Z.~Li and P.~Richt{\'a}rik.
\newblock A unified analysis of stochastic gradient methods for nonconvex
  federated optimization.
\newblock {\em arXiv preprint arXiv:2006.07013}, 2020.

\bibitem{liang2018understanding}
S.~Liang, R.~Sun, Y.~Li, and R.~Srikant.
\newblock Understanding the loss surface of neural networks for binary
  classification.
\newblock In {\em International Conference on Machine Learning}, pages
  2835--2843, 2018.

\bibitem{liu2018zerothorder}
S.~Liu, B.~Kailkhura, P.-Y. Chen, P.~Ting, S.~Chang, and L.~Amini.
\newblock Zeroth-order stochastic variance reduction for nonconvex
  optimization, 2018.

\bibitem{livni2014computational}
R.~Livni, S.~Shalev-Shwartz, and O.~Shamir.
\newblock On the computational efficiency of training neural networks.
\newblock In {\em Advances in neural information processing systems}, pages
  855--863, 2014.

\bibitem{loizou2020stochastic}
N.~Loizou, S.~Vaswani, I.~H. Laradji, and S.~Lacoste-Julien.
\newblock Stochastic polyak step-size for sgd: An adaptive learning rate for
  fast convergence.
\newblock In {\em International Conference on Artificial Intelligence and
  Statistics}, pages 1306--1314. PMLR, 2021.

\bibitem{lojasiewicz1963topological}
S.~Lojasiewicz.
\newblock A topological property of real analytic subsets.
\newblock {\em Coll. du CNRS, Les {\'e}quations aux d{\'e}riv{\'e}es
  partielles}, 117:87--89, 1963.

\bibitem{loshchilov2016sgdr}
I.~Loshchilov and F.~Hutter.
\newblock Sgdr: Stochastic gradient descent with warm restarts.
\newblock {\em arXiv preprint arXiv:1608.03983}, 2016.

\bibitem{lucchi2019stochastic}
A.~Lucchi and J.~Kohler.
\newblock A stochastic tensor method for non-convex optimization.
\newblock {\em arXiv preprint arXiv:1911.10367}, 2019.

\bibitem{ma2018implicit}
C.~Ma, K.~Wang, Y.~Chi, and Y.~Chen.
\newblock Implicit regularization in nonconvex statistical estimation: Gradient
  descent converges linearly for phase retrieval and matrix completion.
\newblock In {\em International Conference on Machine Learning}, pages
  3345--3354. PMLR, 2018.

\bibitem{ma2018power}
S.~Ma, R.~Bassily, and M.~Belkin.
\newblock The power of interpolation: Understanding the effectiveness of sgd in
  modern over-parametrized learning.
\newblock In {\em International Conference on Machine Learning}, pages
  3325--3334. PMLR, 2018.

\bibitem{mairal2015incremental}
J.~Mairal.
\newblock Incremental majorization-minimization optimization with application
  to large-scale machine learning.
\newblock {\em SIAM Journal on Optimization}, 25(2):829--855, 2015.

\bibitem{malik2020derivativefree}
D.~Malik, A.~Pananjady, K.~Bhatia, K.~Khamaru, P.~L. Bartlett, and M.~J.
  Wainwright.
\newblock Derivative-free methods for policy optimization: Guarantees for
  linear quadratic systems, 2020.

\bibitem{martens2010deep}
J.~Martens.
\newblock Deep learning via hessian-free optimization.
\newblock In {\em International Conference on Machine Learning}, volume~27,
  pages 735--742, 2010.

\bibitem{mikolov2012statistical}
T.~Mikolov.
\newblock Statistical language models based on neural networks.
\newblock {\em Presentation at Google, Mountain View, 2nd April}, 80, 2012.

\bibitem{mishchenko2019distributed}
K.~Mishchenko, E.~Gorbunov, M.~Tak{\'a}{\v{c}}, and P.~Richt{\'a}rik.
\newblock Distributed learning with compressed gradient differences.
\newblock {\em arXiv preprint arXiv:1901.09269}, 2019.

\bibitem{mishchenko2020random}
K.~Mishchenko, A.~Khaled, and P.~Richt{\'a}rik.
\newblock Random reshuffling: Simple analysis with vast improvements.
\newblock {\em arXiv preprint arXiv:2006.05988}, 2020.

\bibitem{murty1987some}
K.~G. Murty and S.~N. Kabadi.
\newblock Some np-complete problems in quadratic and nonlinear programming.
\newblock {\em Mathematical Programming}, 39(2):117--129, Jun 1987.

\bibitem{nelder1965simplex}
J.~A. Nelder and R.~Mead.
\newblock A simplex method for function minimization.
\newblock {\em The computer journal}, 7(4):308--313, 1965.

\bibitem{nemirovski1982orth}
A.~Nemirovski.
\newblock Orth-method for smooth convex optimization.
\newblock {\em Izvestia AN SSSR, Transl.: Eng. Cybern. Soviet J. Comput. Syst.
  Sci}, 2:937--947, 1982.

\bibitem{nesterov1983method}
Y.~Nesterov.
\newblock A method of solving a convex programming problem with convergence
  rate $o(1/k^2)$.
\newblock {\em Soviet Mathematics Doklady}, 27(2):372--376, 1983.

\bibitem{nesterov2004introduction}
Y.~Nesterov.
\newblock {\em Introductory Lectures on Convex Optimization: a basic course}.
\newblock Kluwer Academic Publishers, Massachusetts, 2004.

\bibitem{nesterov2012make}
Y.~Nesterov.
\newblock How to make the gradients small.
\newblock {\em Optima}, 88:10--11, 2012.

\bibitem{nesterov2018lectures}
Y.~Nesterov.
\newblock {\em Lectures on convex optimization}, volume 137.
\newblock Springer, 2018.

\bibitem{nesterov2020primal-dual}
Y.~Nesterov, A.~Gasnikov, S.~Guminov, and P.~Dvurechensky.
\newblock Primal-dual accelerated gradient methods with small-dimensional
  relaxation oracle.
\newblock {\em Optimization Methods and Software}, pages 1--28, 2020.
\newblock arXiv:1809.05895.

\bibitem{nesterov2006cubic}
Y.~Nesterov and B.~Polyak.
\newblock Cubic regularization of newton method and its global performance.
\newblock {\em Mathematical Programming}, 108(1):177--205, 2006.

\bibitem{nesterov2017random}
Y.~Nesterov and V.~Spokoiny.
\newblock Random gradient-free minimization of convex functions.
\newblock {\em Found. Comput. Math.}, 17(2):527--566, Apr. 2017.
\newblock First appeared in 2011 as CORE discussion paper 2011/16.

\bibitem{neyshabur2017exploring}
B.~Neyshabur, S.~Bhojanapalli, D.~McAllester, and N.~Srebro.
\newblock Exploring generalization in deep learning.
\newblock In {\em Advances in neural information processing systems}, pages
  5947--5956, 2017.

\bibitem{nguyen2017sarah}
L.~M. Nguyen, J.~Liu, K.~Scheinberg, and M.~Tak{\'a}{\v{c}}.
\newblock Sarah: A novel method for machine learning problems using stochastic
  recursive gradient.
\newblock In {\em International Conference on Machine Learning}, pages
  2613--2621, 2017.

\bibitem{nguyen2017stochastic}
L.~M. Nguyen, J.~Liu, K.~Scheinberg, and M.~Tak{\'a}{\v{c}}.
\newblock Stochastic recursive gradient algorithm for nonconvex optimization.
\newblock {\em arXiv preprint arXiv:1705.07261}, 2017.

\bibitem{nguyen2020unified}
L.~M. Nguyen, Q.~Tran-Dinh, D.~T. Phan, P.~H. Nguyen, and M.~van Dijk.
\newblock A unified convergence analysis for shuffling-type gradient methods.
\newblock {\em arXiv preprint arXiv:2002.08246}, 2020.

\bibitem{nguyen2018loss}
Q.~Nguyen, M.~C. Mukkamala, and M.~Hein.
\newblock On the loss landscape of a class of deep neural networks with no bad
  local valleys.
\newblock {\em arXiv preprint arXiv:1809.10749}, 2018.

\bibitem{nocedal2006numerical}
J.~Nocedal and S.~Wright.
\newblock {\em Numerical optimization}.
\newblock Springer Science \& Business Media, 2006.

\bibitem{osawa2018second}
K.~Osawa, Y.~Tsuji, Y.~Ueno, A.~Naruse, R.~Yokota, and S.~Matsuoka.
\newblock Second-order optimization method for large mini-batch: Training
  resnet-50 on imagenet in 35 epochs.
\newblock {\em arXiv preprint arXiv:1811.12019}, 1:2, 2018.

\bibitem{papernot2017practical}
N.~Papernot, P.~McDaniel, I.~Goodfellow, S.~Jha, Z.~B. Celik, and A.~Swami.
\newblock Practical black-box attacks against machine learning, 2017.

\bibitem{papyan2017convolutional}
V.~Papyan, Y.~Romano, J.~Sulam, and M.~Elad.
\newblock Convolutional dictionary learning via local processing.
\newblock In {\em Proceedings of the IEEE International Conference on Computer
  Vision}, pages 5296--5304, 2017.

\bibitem{park2020combining}
S.~Park, S.~H. Jung, and P.~M. Pardalos.
\newblock Combining stochastic adaptive cubic regularization with negative
  curvature for nonconvex optimization.
\newblock {\em Journal of Optimization Theory and Applications},
  184(3):953--971, 2020.

\bibitem{pascanu2013difficulty}
R.~Pascanu, T.~Mikolov, and Y.~Bengio.
\newblock On the difficulty of training recurrent neural networks.
\newblock In {\em International conference on machine learning}, pages
  1310--1318, 2013.

\bibitem{polyak1963gradient}
B.~Polyak.
\newblock Gradient methods for the minimisation of functionals.
\newblock {\em USSR Computational Mathematics and Mathematical Physics},
  3(4):864 -- 878, 1963.

\bibitem{polyak1987introduction}
B.~Polyak.
\newblock {\em Introduction to Optimization}.
\newblock New York, Optimization Software, 1987.

\bibitem{polyak1964some}
B.~T. Polyak.
\newblock Some methods of speeding up the convergence of iteration methods.
\newblock {\em USSR Computational Mathematics and Mathematical Physics},
  4(5):1--17, 1964.

\bibitem{qu2019nonconvex}
Q.~Qu, X.~Li, and Z.~Zhu.
\newblock A nonconvex approach for exact and efficient multichannel sparse
  blind deconvolution.
\newblock In {\em Advances in Neural Information Processing Systems}, pages
  4015--4026, 2019.

\bibitem{rajput2020closing}
S.~Rajput, A.~Gupta, and D.~Papailiopoulos.
\newblock Closing the convergence gap of sgd without replacement.
\newblock {\em arXiv preprint arXiv:2002.10400}, 2020.

\bibitem{reddi2016stochastic}
S.~J. Reddi, A.~Hefny, S.~Sra, B.~Poczos, and A.~Smola.
\newblock Stochastic variance reduction for nonconvex optimization.
\newblock In {\em International conference on machine learning}, pages
  314--323, 2016.

\bibitem{reddi2019convergence}
S.~J. Reddi, S.~Kale, and S.~Kumar.
\newblock On the convergence of adam and beyond.
\newblock {\em arXiv preprint arXiv:1904.09237}, 2019.

\bibitem{reddi2016proximal}
S.~J. Reddi, S.~Sra, B.~Poczos, and A.~J. Smola.
\newblock Proximal stochastic methods for nonsmooth nonconvex finite-sum
  optimization.
\newblock In {\em Advances in Neural Information Processing Systems}, pages
  1145--1153, 2016.

\bibitem{rezanov2021deep}
A.~Rezanov and D.~Yudin.
\newblock Deep neural networks for ortophoto-based vehicle localization.
\newblock In B.~Kryzhanovsky, W.~Dunin-Barkowski, V.~Redko, and Y.~Tiumentsev,
  editors, {\em Advances in Neural Computation, Machine Learning, and Cognitive
  Research IV}, pages 167--174, Cham, 2021. Springer International Publishing.

\bibitem{Risteski2016AlgorithmsAM}
A.~Risteski and Y.~Li.
\newblock Algorithms and matching lower bounds for approximately-convex
  optimization.
\newblock In {\em NIPS}, 2016.

\bibitem{roy2020escaping}
A.~Roy, K.~Balasubramanian, S.~Ghadimi, and P.~Mohapatra.
\newblock Escaping saddle-point faster under interpolation-like conditions.
\newblock {\em Advances in Neural Information Processing Systems}, 33, 2020.

\bibitem{royer2018complexity}
C.~W. Royer and S.~J. Wright.
\newblock Complexity analysis of second-order line-search algorithms for smooth
  nonconvex optimization.
\newblock {\em SIAM Journal on Optimization}, 28(2):1448--1477, 2018.

\bibitem{safran2018spurious}
I.~Safran and O.~Shamir.
\newblock Spurious local minima are common in two-layer relu neural networks.
\newblock In {\em International Conference on Machine Learning}, pages
  4433--4441. PMLR, 2018.

\bibitem{sankararaman2019impact}
K.~A. Sankararaman, S.~De, Z.~Xu, W.~R. Huang, and T.~Goldstein.
\newblock The impact of neural network overparameterization on gradient
  confusion and stochastic gradient descent.
\newblock {\em arXiv preprint arXiv:1904.06963}, 2019.

\bibitem{schmidt2017minimizing}
M.~Schmidt, N.~Le~Roux, and F.~Bach.
\newblock Minimizing finite sums with the stochastic average gradient.
\newblock {\em Mathematical Programming}, 162(1-2):83--112, 2017.

\bibitem{schmidt2013fast}
M.~Schmidt and N.~L. Roux.
\newblock Fast convergence of stochastic gradient descent under a strong growth
  condition.
\newblock {\em arXiv preprint arXiv:1308.6370}, 2013.

\bibitem{Schumer1968}
M.~Schumer and K.~Steiglitz.
\newblock Adaptive step size random search.
\newblock {\em {IEEE} Transactions on Automatic Control}, 13(3):270--276, June
  1968.

\bibitem{sebbouh2020convergence}
O.~Sebbouh, R.~M. Gower, and A.~Defazio.
\newblock On the convergence of the stochastic heavy ball method.
\newblock {\em arXiv preprint arXiv:2006.07867}, 2020.

\bibitem{Sener2020Learning}
O.~Sener and V.~Koltun.
\newblock Learning to guide random search.
\newblock In {\em International Conference on Learning Representations}, 2020.

\bibitem{shalev2016sdca}
S.~Shalev-Shwartz.
\newblock Sdca without duality, regularization, and individual convexity.
\newblock In {\em International Conference on Machine Learning}, pages
  747--754, 2016.

\bibitem{shechtman2015phase}
Y.~Shechtman, Y.~C. Eldar, O.~Cohen, H.~N. Chapman, J.~Miao, and M.~Segev.
\newblock Phase retrieval with application to optical imaging: a contemporary
  overview.
\newblock {\em IEEE signal processing magazine}, 32(3):87--109, 2015.

\bibitem{shen2019stochastic}
Z.~Shen, P.~Zhou, C.~Fang, and A.~Ribeiro.
\newblock A stochastic trust region method for non-convex minimization.
\newblock {\em arXiv preprint arXiv:1903.01540}, 2019.

\bibitem{shi2020learning}
B.~Shi, W.~J. Su, and M.~I. Jordan.
\newblock On learning rates and schr$\backslash$" odinger operators.
\newblock {\em arXiv preprint arXiv:2004.06977}, 2020.

\bibitem{shi2020manifold}
L.~Shi and Y.~Chi.
\newblock Manifold gradient descent solves multi-channel sparse blind
  deconvolution provably and efficiently.
\newblock In {\em ICASSP 2020-2020 IEEE International Conference on Acoustics,
  Speech and Signal Processing (ICASSP)}, pages 5730--5734. IEEE, 2020.

\bibitem{shi2020rmsprop}
N.~Shi, D.~Li, M.~Hong, and R.~Sun.
\newblock Rmsprop converges with proper hyper-parameter.
\newblock In {\em International Conference on Learning Representations}, 2021.

\bibitem{shibaev2021zeroth-order}
I.~Shibaev, P.~Dvurechensky, and A.~Gasnikov.
\newblock Zeroth-order methods for noisy {H}\"older-gradient functions.
\newblock {\em Optimization Letters}, 2021.
\newblock (accepted), arXiv:2006.11857, doi:10.1007/s11590-021-01742-z.

\bibitem{shin2019effects}
Y.~Shin.
\newblock Effects of depth, width, and initialization: A convergence analysis
  of layer-wise training for deep linear neural networks.
\newblock {\em arXiv preprint arXiv:1910.05874}, 2019.

\bibitem{shor1967generalized}
N.~Z. Shor.
\newblock Generalized gradient descent with application to block programming.
\newblock {\em Kibernetika}, 3(3):53--55, 1967.

\bibitem{skrynnik2021forgetful}
A.~Skrynnik, A.~Staroverov, E.~Aitygulov, K.~Aksenov, V.~Davydov, and A.~I.
  Panov.
\newblock Forgetful experience replay in hierarchical reinforcement learning
  from expert demonstrations.
\newblock {\em Knowledge-Based Systems}, 218:106844, 2021.

\bibitem{smith2017cyclical}
L.~N. Smith.
\newblock Cyclical learning rates for training neural networks.
\newblock In {\em 2017 IEEE Winter Conference on Applications of Computer
  Vision (WACV)}, pages 464--472. IEEE, 2017.

\bibitem{solodov1998incremental}
M.~V. Solodov.
\newblock Incremental gradient algorithms with stepsizes bounded away from
  zero.
\newblock {\em Computational Optimization and Applications}, 11(1):23--35,
  1998.

\bibitem{soltanolkotabi2018theoretical}
M.~Soltanolkotabi, A.~Javanmard, and J.~D. Lee.
\newblock Theoretical insights into the optimization landscape of
  over-parameterized shallow neural networks.
\newblock {\em IEEE Transactions on Information Theory}, 65(2):742--769, 2018.

\bibitem{spokoiny2012parametric}
V.~Spokoiny et~al.
\newblock Parametric estimation. finite sample theory.
\newblock {\em The Annals of Statistics}, 40(6):2877--2909, 2012.

\bibitem{stonyakin2020inexact}
F.~Stonyakin, A.~Tyurin, A.~Gasnikov, P.~Dvurechensky, A.~Agafonov,
  D.~Dvinskikh, M.~Alkousa, D.~Pasechnyuk, S.~Artamonov, and V.~Piskunova.
\newblock Inexact model: A framework for optimization and variational
  inequalities.
\newblock {\em Optimization Methods and Software}, 2021.
\newblock (accepted), WIAS Preprint No. 2709, arXiv:2001.09013,
  arXiv:1902.00990, doi:10.1080/10556788.2021.1924714.

\bibitem{stonyakin2019gradient}
F.~S. Stonyakin, D.~Dvinskikh, P.~Dvurechensky, A.~Kroshnin, O.~Kuznetsova,
  A.~Agafonov, A.~Gasnikov, A.~Tyurin, C.~A. Uribe, D.~Pasechnyuk, and
  S.~Artamonov.
\newblock Gradient methods for problems with inexact model of the objective.
\newblock In M.~Khachay, Y.~Kochetov, and P.~Pardalos, editors, {\em
  Mathematical Optimization Theory and Operations Research}, pages 97--114,
  Cham, 2019. Springer International Publishing.
\newblock arXiv:1902.09001.

\bibitem{sun2019optimization}
R.~Sun.
\newblock Optimization for deep learning: theory and algorithms.
\newblock {\em arXiv preprint arXiv:1912.08957}, 2019.

\bibitem{surazhevsky2021noise-assisted}
I.~Surazhevsky, V.~Demin, A.~Ilyasov, A.~Emelyanov, K.~Nikiruy, V.~Rylkov,
  S.~Shchanikov, I.~Bordanov, S.~Gerasimova, D.~Guseinov, N.~Malekhonova,
  D.~Pavlov, A.~Belov, A.~Mikhaylov, V.~Kazantsev, D.~Valenti, B.~Spagnolo, and
  M.~Kovalchuk.
\newblock Noise-assisted persistence and recovery of memory state in a
  memristive spiking neuromorphic network.
\newblock {\em Chaos, Solitons \& Fractals}, 146:110890, 2021.

\bibitem{sutskever2013importance}
I.~Sutskever, J.~Martens, G.~Dahl, and G.~Hinton.
\newblock On the importance of initialization and momentum in deep learning.
\newblock In {\em International conference on machine learning}, pages
  1139--1147, 2013.

\bibitem{swirszcz2016local}
G.~Swirszcz, W.~M. Czarnecki, and R.~Pascanu.
\newblock Local minima in training of deep networks.
\newblock 2016.

\bibitem{tan2019online}
Y.~S. Tan and R.~Vershynin.
\newblock Online stochastic gradient descent with arbitrary initialization
  solves non-smooth, non-convex phase retrieval.
\newblock {\em arXiv preprint arXiv:1910.12837}, 2019.

\bibitem{tao2018primal}
W.~Tao, Z.~Pan, G.~Wu, and Q.~Tao.
\newblock Primal averaging: A new gradient evaluation step to attain the
  optimal individual convergence.
\newblock {\em IEEE transactions on cybernetics}, 50(2):835--845, 2018.

\bibitem{taylor2019stochastic}
A.~Taylor and F.~Bach.
\newblock Stochastic first-order methods: non-asymptotic and computer-aided
  analyses via potential functions.
\newblock In {\em Conference on Learning Theory}, pages 2934--2992, 2019.

\bibitem{tieleman2012lecture}
T.~Tieleman and G.~Hinton.
\newblock Lecture 6.5-rmsprop: Divide the gradient by a running average of its
  recent magnitude.
\newblock {\em COURSERA: Neural networks for machine learning}, 4(2):26--31,
  2012.

\bibitem{tripuraneni2018stochastic}
N.~Tripuraneni, M.~Stern, C.~Jin, J.~Regier, and M.~I. Jordan.
\newblock Stochastic cubic regularization for fast nonconvex optimization.
\newblock In {\em Advances in neural information processing systems}, pages
  2899--2908, 2018.

\bibitem{tseng1998incremental}
P.~Tseng.
\newblock An incremental gradient (-projection) method with momentum term and
  adaptive stepsize rule.
\newblock {\em SIAM Journal on Optimization}, 8(2):506--531, 1998.

\bibitem{usmanova2017master}
I.~Usmanova.
\newblock Robust solutions to stochastic optimization problems.
\newblock {\em Master Thesis (MSIAM); Institut Polytechnique de Grenoble
  ENSIMAG, Laboratoire Jean Kuntzmann}, 2017.

\bibitem{vaswani2017attention}
A.~Vaswani, N.~Shazeer, N.~Parmar, J.~Uszkoreit, L.~Jones, A.~N. Gomez,
  {\L}.~Kaiser, and I.~Polosukhin.
\newblock Attention is all you need.
\newblock In {\em Advances in neural information processing systems}, pages
  5998--6008, 2017.

\bibitem{vaswani2018fast}
S.~Vaswani, F.~Bach, and M.~Schmidt.
\newblock Fast and faster convergence of sgd for over-parameterized models and
  an accelerated perceptron.
\newblock In {\em The 22nd International Conference on Artificial Intelligence
  and Statistics}, pages 1195--1204. PMLR, 2019.

\bibitem{vaswani2019painless}
S.~Vaswani, A.~Mishkin, I.~Laradji, M.~Schmidt, G.~Gidel, and
  S.~Lacoste-Julien.
\newblock Painless stochastic gradient: Interpolation, line-search, and
  convergence rates.
\newblock In {\em Advances in Neural Information Processing Systems}, pages
  3732--3745, 2019.

\bibitem{vavasis1993black}
S.~A. Vavasis.
\newblock Black-box complexity of local minimization.
\newblock {\em SIAM Journal on Optimization}, 3(1):60--80, 1993.

\bibitem{vidal2017mathematics}
R.~Vidal, J.~Bruna, R.~Giryes, and S.~Soatto.
\newblock Mathematics of deep learning.
\newblock {\em arXiv preprint arXiv:1712.04741}, 2017.

\bibitem{wang2018spiderboost}
Z.~Wang, K.~Ji, Y.~Zhou, Y.~Liang, and V.~Tarokh.
\newblock Spiderboost: A class of faster variance-reduced algorithms for
  nonconvex optimization.
\newblock {\em arXiv preprint arXiv:1810.10690}, 2018.

\bibitem{wang2019spiderboost}
Z.~Wang, K.~Ji, Y.~Zhou, Y.~Liang, and V.~Tarokh.
\newblock Spiderboost and momentum: Faster variance reduction algorithms.
\newblock In {\em Advances in Neural Information Processing Systems}, pages
  2403--2413, 2019.

\bibitem{wang2019stochastic}
Z.~Wang, Y.~Zhou, Y.~Liang, and G.~Lan.
\newblock Stochastic variance-reduced cubic regularization for nonconvex
  optimization.
\newblock In {\em The 22nd International Conference on Artificial Intelligence
  and Statistics}, pages 2731--2740. PMLR, 2019.

\bibitem{wang2020cubic}
Z.~Wang, Y.~Zhou, Y.~Liang, and G.~Lan.
\newblock Cubic regularization with momentum for nonconvex optimization.
\newblock In {\em Uncertainty in Artificial Intelligence}, pages 313--322.
  PMLR, 2020.

\bibitem{ward2019adagrad}
R.~Ward, X.~Wu, and L.~Bottou.
\newblock Adagrad stepsizes: Sharp convergence over nonconvex landscapes.
\newblock In {\em International Conference on Machine Learning}, pages
  6677--6686. PMLR, 2019.

\bibitem{wilson2017marginal}
A.~C. Wilson, R.~Roelofs, M.~Stern, N.~Srebro, and B.~Recht.
\newblock The marginal value of adaptive gradient methods in machine learning.
\newblock In {\em Advances in neural information processing systems}, pages
  4148--4158, 2017.

\bibitem{wright2018optimization}
S.~J. Wright.
\newblock Optimization algorithms for data analysis.
\newblock {\em The Mathematics of Data}, 25:49, 2018.

\bibitem{wu2020hadamard}
F.~Wu and P.~Rebeschini.
\newblock Hadamard wirtinger flow for sparse phase retrieval.
\newblock {\em arXiv preprint arXiv:2006.01065}, 2020.

\bibitem{xie2019general}
G.~Xie, L.~Luo, and Z.~Zhang.
\newblock A general analysis framework of lower complexity bounds for
  finite-sum optimization.
\newblock {\em arXiv preprint arXiv:1908.08394}, 2019.

\bibitem{xu2018global}
P.~Xu, J.~Chen, D.~Zou, and Q.~Gu.
\newblock Global convergence of langevin dynamics based algorithms for
  nonconvex optimization.
\newblock In {\em Advances in Neural Information Processing Systems}, pages
  3122--3133, 2018.

\bibitem{xu2020newton}
P.~Xu, F.~Roosta, and M.~W. Mahoney.
\newblock Newton-type methods for non-convex optimization under inexact hessian
  information.
\newblock {\em Mathematical Programming}, 184(1):35--70, 2020.

\bibitem{xu2020second}
P.~Xu, F.~Roosta, and M.~W. Mahoney.
\newblock Second-order optimization for non-convex machine learning: An
  empirical study.
\newblock In {\em Proceedings of the 2020 SIAM International Conference on Data
  Mining}, pages 199--207. SIAM, 2020.

\bibitem{xu2020momentum}
Y.~Xu.
\newblock Momentum-based variance-reduced proximal stochastic gradient method
  for composite nonconvex stochastic optimization.
\newblock {\em arXiv preprint arXiv:2006.00425}, 2020.

\bibitem{xu2018first}
Y.~Xu, R.~Jin, and T.~Yang.
\newblock First-order stochastic algorithms for escaping from saddle points in
  almost linear time.
\newblock In {\em Advances in Neural Information Processing Systems}, pages
  5530--5540, 2018.

\bibitem{yan2018unified}
Y.~Yan, T.~Yang, Z.~Li, Q.~Lin, and Y.~Yang.
\newblock A unified analysis of stochastic momentum methods for deep learning.
\newblock In {\em Proceedings of the 27th International Joint Conference on
  Artificial Intelligence}, pages 2955--2961, 2018.

\bibitem{yang2019misspecified}
Z.~Yang, L.~F. Yang, E.~X. Fang, T.~Zhao, Z.~Wang, and M.~Neykov.
\newblock Misspecified nonconvex statistical optimization for sparse phase
  retrieval.
\newblock {\em Mathematical Programming}, 176(1-2):545--571, 2019.

\bibitem{yun2018small}
C.~Yun, S.~Sra, and A.~Jadbabaie.
\newblock Small nonlinearities in activation functions create bad local minima
  in neural networks.
\newblock {\em arXiv preprint arXiv:1802.03487}, 2018.

\bibitem{yun2020general}
J.~Yun, A.~C. Lozano, and E.~Yang.
\newblock A general family of stochastic proximal gradient methods for deep
  learning.
\newblock {\em arXiv preprint arXiv:2007.07484}, 2020.

\bibitem{zaheer2018adaptive}
M.~Zaheer, S.~Reddi, D.~Sachan, S.~Kale, and S.~Kumar.
\newblock Adaptive methods for nonconvex optimization.
\newblock In {\em Advances in neural information processing systems}, pages
  9793--9803, 2018.

\bibitem{zhang2020improved}
B.~Zhang, J.~Jin, C.~Fang, and L.~Wang.
\newblock Improved analysis of clipping algorithms for non-convex optimization.
\newblock {\em Advances in Neural Information Processing Systems}, 33, 2020.

\bibitem{zhang2016understanding}
C.~Zhang, S.~Bengio, M.~Hardt, B.~Recht, and O.~Vinyals.
\newblock Understanding deep learning (still) requires rethinking
  generalization.
\newblock {\em Communications of the ACM}, 64(3):107--115, 2021.

\bibitem{zhang2021general}
H.~Zhang, Y.~Bi, and J.~Lavaei.
\newblock General low-rank matrix optimization: Geometric analysis and sharper
  bounds.
\newblock {\em arXiv preprint arXiv:2104.10356}, 2021.

\bibitem{zhang2019gradient}
J.~Zhang, T.~He, S.~Sra, and A.~Jadbabaie.
\newblock Why gradient clipping accelerates training: A theoretical
  justification for adaptivity.
\newblock In {\em International Conference on Learning Representations}, 2020.

\bibitem{zhang2019adam}
J.~Zhang, S.~P. Karimireddy, A.~Veit, S.~Kim, S.~Reddi, S.~Kumar, and S.~Sra.
\newblock Why are adaptive methods good for attention models?
\newblock {\em Advances in Neural Information Processing Systems}, 33, 2020.

\bibitem{zhang2020stochastic}
J.~Zhang and L.~Xiao.
\newblock Stochastic variance-reduced prox-linear algorithms for nonconvex
  composite optimization.
\newblock {\em arXiv preprint arXiv:2004.04357}, 2020.

\bibitem{zhang2018adaptive}
J.~Zhang, L.~Xiao, and S.~Zhang.
\newblock Adaptive stochastic variance reduction for subsampled newton method
  with cubic regularization.
\newblock {\em arXiv preprint arXiv:1811.11637}, 2018.

\bibitem{zhang2021sharp}
R.~Y. Zhang.
\newblock Sharp global guarantees for nonconvex low-rank matrix recovery in the
  overparameterized regime.
\newblock {\em arXiv preprint arXiv:2104.10790}, 2021.

\bibitem{zhang2020from}
Y.~Zhang, Q.~Qu, and J.~Wright.
\newblock From symmetry to geometry: Tractable nonconvex problems.
\newblock {\em arXiv preprint arXiv:2007.06753}, 2020.

\bibitem{zhang2020symmetry}
Y.~Zhang, Q.~Qu, and J.~Wright.
\newblock From symmetry to geometry: Tractable nonconvex problems.
\newblock {\em arXiv preprint arXiv:2007.06753}, 2020.

\bibitem{zhang2020boosting}
Y.~Zhang, Y.~Zhou, K.~Ji, and M.~M. Zavlanos.
\newblock Boosting one-point derivative-free online optimization via residual
  feedback, 2020.

\bibitem{zhigljavsky2007stochastic}
A.~Zhigljavsky and A.~Zilinskas.
\newblock {\em Stochastic global optimization}, volume~9.
\newblock Springer Science \& Business Media, 2007.

\bibitem{zhou2019lower}
D.~Zhou and Q.~Gu.
\newblock Lower bounds for smooth nonconvex finite-sum optimization.
\newblock In {\em International Conference on Machine Learning}, pages
  7574--7583, 2019.

\bibitem{zhou2020stochastic}
D.~Zhou and Q.~Gu.
\newblock Stochastic recursive variance-reduced cubic regularization methods.
\newblock In {\em International Conference on Artificial Intelligence and
  Statistics}, pages 3980--3990. PMLR, 2020.

\bibitem{zhou2018convergence}
D.~Zhou, Y.~Tang, Z.~Yang, Y.~Cao, and Q.~Gu.
\newblock On the convergence of adaptive gradient methods for nonconvex
  optimization.
\newblock {\em arXiv preprint arXiv:1808.05671}, 2018.

\bibitem{zhou2018stochastic}
D.~Zhou, P.~Xu, and Q.~Gu.
\newblock Stochastic nested variance reduced gradient descent for nonconvex
  optimization.
\newblock {\em Advances in neural information processing systems}, 2018.

\bibitem{zhou2019stochastic}
D.~Zhou, P.~Xu, and Q.~Gu.
\newblock Stochastic variance-reduced cubic regularization methods.
\newblock {\em Journal of Machine Learning Research}, 20(134):1--47, 2019.

\bibitem{zhou2019stochastic_jmlr}
D.~Zhou, P.~Xu, and Q.~Gu.
\newblock Stochastic variance-reduced cubic regularization methods.
\newblock {\em Journal of Machine Learning Research}, 20(134):1--47, 2019.

\bibitem{zhu2020adaptive}
X.~Zhu, J.~Han, and B.~Jiang.
\newblock An adaptive high order method for finding third-order critical points
  of nonconvex optimization.
\newblock {\em arXiv preprint arXiv:2008.04191}, 2020.

\end{thebibliography}

\end{document}